\documentclass[12pt]{article}
\usepackage{amssymb,amsmath}
\usepackage{amscd, graphics}
\usepackage{mathrsfs}
\usepackage{ascmac}
\usepackage{txfonts}
\usepackage{color}%SLIM

\allowdisplaybreaks[1]
\newcommand{\T}{\mathbb{T}}
\newcommand{\Z}{\mathbb{Z}}

\newcommand{\C}{\mathbb{C}}
\newcommand{\N}{\mathbb{N}}
\newcommand{\R}{\mathbb{R}}

\newcommand{\dd}{\,{\rm d}}

\renewcommand{\Re}{\mathop{\mathrm{Re}}}
\renewcommand{\Im}{\mathop{\mathrm{Im}}}

\numberwithin{equation}{section}
\newtheorem{thm}{Theorem}[section]

\newtheorem{prop}[thm]{Proposition}
\newtheorem{lem}[thm]{Lemma}
\newtheorem{rem}[thm]{Remark}
\newtheorem{cor}[thm]{Corollary}
\newtheorem{assum}{Assumption}

\begin{document}

%\author{   \footnote{Department of Mathematics and Statistics, University of Victoria, PO BOX 1700 STN CSC, Victoria BC, V8W 2Y2\,.\,   Email: ibrahims@uvic.ca} Slim Ibrahim \and \footnote{Department of Mathematics, Graduate School of Science, Kyoto University, Kitashirakawa Oiwake-cho, Sakyo-ku, Kyoto 606-8502, Japan\,.\,  Email: maekawa@math.kyoto-u.ac.jp} Yasunori Maekawa \and 
%\footnote{Courant Institute of Mathematical Sciences, 251 Mercer Street, New York, NY 10012, USA\,.\,  Email: masmoudi@courant.nyu.edu} Nader Masmoudi}

%\author{S. Ibrahim,  Y. Maekawa \& N. Masmoudi}
% It is required to enter 2010 MSC.
%%%\subjclass{Primary: 35K50, 35B40; Secondary: 35K55, 35K57.}
% Please provide minimum  5 keywords.
%% \keywords{Turbulent flows, Navier-Stokes and Euler equations, Pseudo-spectrum, enhanced dissipation, Metastability}

% Email address of each of all authors is required.
% You may list email addresses of all other authors, separately.
%\email{S. Ibrahim}{ibrahim@math.uvic.ca}
% \email[Y. Maekawa]{maekawa@math.kyoto-u.ac.jp}
% \email[N. Masmoudi]{masmoudi@courant.nyu.edu}

\title{On pseudospectral bound for non-selfadjoint operators and its application to stability of Kolmogorov flows}

%\author[1, 4]{\small Slim Ibrahim}
%\author[2]{\small Yasunori Maekawa}
%\author[3, 4]{\small Nader Masmoudi}

%\affil[1]{\footnotesize Department of Mathematics and Statistics, University of Victoria, PO BOX 1700 STN CSC, Victoria BC, Canada}
%\affil[2]{\footnotesize Department of Mathematics, Graduate School of Science, Kyoto University, Kitashirakawa Oiwake-cho, Sakyo-ku, Kyoto 606-8502, Japan}
%\affil[3]{\footnotesize Courant Institute of Mathematical Sciences, 251 Mercer Street, New York, NY 10012, USA}
%\affil[4]{\footnotesize Department of mathematics, New York University in Abu Dhabi, Saadyiat Island, Abu Dhabi, UAE}

\date{}

\author{Slim Ibrahim\thanks{Department of Mathematics and Statistics, University of Victoria, PO BOX 1700 STN CSC, Victoria BC, Canada. E-mail: \texttt{ibrahims@uvic.ca}}  \and Yasunori Maekawa\thanks{Department of Mathematics, Kyoto University, Kitashirakawa Oiwakecho, Sakyo-ku, Kyoto 606-8502, Japan. E-mail: \texttt{maekawa@math.kyoto-u.ac.jp}}  \and Nader Masmoudi\thanks{Department of mathematics, New York University in Abu Dhabi, Saadyiat Island, Abu Dhabi, UAE. E-mail: \texttt{masmoudi@cims.nyu.edu}}}

%\date{}

\maketitle

\noindent {\bf Abstract} 
We study the stability of the Kolmogorov flows which are stationary solutions to the two-dimensional Navier-Stokes equations in the presence of the shear  external force. We establish the linear stability estimate when the viscosity coefficient $\nu$ is sufficiently small, where the enhanced dissipation is rigorously verified in the time scale $O(\nu^{-\frac12})$ for solutions to the linearized problem, which has been numerically conjectured and is much shorter than the usual viscous time scale $O(\nu^{-1})$. Our approach is based on the detailed analysis for the resolvent problem. We also provide the abstract framework which is applicable to the resolvent estimate for the Kolmogorov flows.

\vspace{0.3cm}


\vspace{0.3cm}


\section{Introduction}

For nearly-inviscid fluids, turbulent phenomena often occur at transient time scales that are much smaller than the viscous time scale. Describing the fluid, by means of simple solutions, for such long transient times helps to understand turbulence. This is of course of great interest both physically and mathematically. But finding such solutions and estimating their basin of attraction are in general not easy tasks both experimentally and theoretically. To investigate this phenomena let us consider the two dimensional incompressible Navier-Stokes equations in the domain $\mathbb M=\mathbb T^2$ or $\mathbb M=\mathbb R^2$,
\begin{equation}
\label{2dNS}
\partial_t U +(U\cdot\nabla)U+\nabla P=\nu\Delta U + F\,, \qquad t>0\,, \quad  (x,y) \in \mathbb M\,.
\end{equation}
Here  $U=(U_1,U_2):\mathbb M^2\times(0,\infty)\to\mathbb R^2$ is the velocity field of a fluid, $P:\mathbb M^2\times(0,\infty)\to \R$ is the pressure field, and $\nu>0$ is the viscosity coefficient. The vector field $F$ describes a given external force. Setting the vorticity $\Omega$ as $\Omega={\rm rot}\, U=\partial_x U_2-\partial_y U_1$, one can rewrite \eqref{2dNS} in the vorticity form
\begin{equation}
\label{2dVort}
\partial_t\Omega+(U\cdot\nabla)\Omega=\nu\Delta\Omega + {\rm rot}\, F\,.
\end{equation}
Recall that the velocity field can be formally recovered from its vorticity using the Biot-Savart law:
\begin{equation}
\label{BioSav}
U=K_{BS}*\Omega\,.
\end{equation}
Here the kernel $K_{BS}$ is given by $K_{BS}(x,y) =\frac1{2\pi}\frac{(-y,x)}{x^2+y^2}$ when $\mathbb{M}=\R^2$, and $*$ denotes the convolution with respect to the spatial variables. In the sequel, we will review two important examples of  solutions to \eqref{2dVort}, the Kolmogorov flow and the Lamb-Oseen vortex, and explain how the study of their stability is related to spectral problems for non-self adjoint operators.

The Kolmogorov flow, which is the main object of this paper, is an explicit stationary solution to \eqref{2dNS} with a shear sourcing term $F=(a \nu \sin  y, 0)$, $a\in \R$, and is given by 
\begin{align}\label{intro.def.Kol}
U^{a} (x,y)= a (\sin y,0)\,, \qquad \Omega^{a}(x,y)= - a \cos y\,. 
\end{align}
By Iudovich \cite{Iu} these solutions are known to be globally stable for initial perturbations in Sobolev class with zero mean condition for the streamfunctions; see also Marchioro \cite{Mar}. 
By changing the length of the periodicity (e.g., for $x$) the detailed bifurcation analysis has also been done, and there are a lot of important works in this direction; see, for example, \cite{MeSi,Iu,MaMi,AB,Ya,OS}. As a closely related subject of this paper, there are also explicit solutions having the similar forms to \eqref{intro.def.Kol} when $F=0$, but instead, the initial data is chosen as in \eqref{intro.def.Kol}. Indeed, in this case one can check that $U^{a} (x,y,t)= a e^{-\nu t}  (\sin y,0)$ solves \eqref{2dNS} with $F=0$. These solutions 
describe a quasi-steady state of the fluid, and are exact steady solutions to the Euler equations when $\nu=0$. These quasi-steady solutions are known as ``bar-states" or also as the Kolmogorov flows, and they qualitatively match the maximum entropy solutions found in \cite{BS,MSMOM,YMC}. 
Both  numerical and experimental evidences \cite{YMC} claim that solutions to \eqref{2dNS} rapidly approach bar-states on time scale $\mathcal O(\frac1{\sqrt\nu})$ for high Reynolds number. Note that the time scale $\mathcal O(\frac1{\sqrt\nu})$ is much  shorter than the scale $\mathcal{O}(\frac{1}{\nu})$ which is the scale for the linear Stokes equation (and thus, heat equation in this problem) with the viscosity $\nu$. 

The aim of this paper is to study  this enhanced dissipation in view of the stability analysis of the steady Kolmogorov flows \eqref{intro.def.Kol}.  Expanding solutions to \eqref{2dVort} around \eqref{intro.def.Kol} yields
\begin{equation}
\label{2dLinNS}
\partial_t \omega = \mathcal L^{\nu, a}  \omega + \text{nonlinear term}
\end{equation}
where we have set $\Omega=\Omega^{a}+\omega$, and the linearized operator $\mathcal{L}^{\nu, a}$ is given by 
\begin{align}\label{def.L^nu}
\mathcal L^{\nu, a} \omega = \nu\Delta \omega - a  \sin y \, \partial_x (I+\Delta^{-1}) \omega \,.
\end{align}
We note that the linearized operator around the bar-state has the similar form but becomes time-dependent as
\begin{align}\label{def.L^nu'}
\mathcal L^{\nu, a} (t) \omega = \nu\Delta \omega -a  e^{-\nu t} \sin y \, \partial_x (I+\Delta^{-1}) \omega \,.
\end{align}
Showing that the solution $\omega$ to \eqref{2dLinNS} decays rapidly within a nontrivial time scale $t\ll O(\frac{1}{\nu})$ is a challenging mathematical problem, even in the linear case.
In studying the flows generated by \eqref{def.L^nu} or \eqref{def.L^nu'} the main difficulty comes from the presence of the non-local term in these linearized operators. In \cite{BW}, Beck and Wayne proved the stability  and enhanced dissipation of the bar-states for the model linear problem by removing the nonlocal term $\Delta^{-1}$ from \eqref{def.L^nu'}. Their method is based on hypocoercivity arguments developed by Villani \cite{Vi}, and provide the decay in the time scale $\mathcal{O}(\frac{1}{\sqrt{\nu}})$ for solutions to the model linear problem in a suitable invariant subspace. 
However, it is not clear how to extend their argument in the presence of the nonlocal term $\Delta^{-1}$.
Moreover, beside the nonlocality,  the presence of $\Delta^{-1}$ in \eqref{def.L^nu} or \eqref{def.L^nu'} leads to an additional difficulty in view of the symmetry of the operator. Indeed, although the operator $\sin y \partial_x$ is antisymmetric in the standard $L^2$ space,  $\sin y \partial_x (I+ \Delta^{-1})$ is not.
Very recently, the full evolution operator \eqref{def.L^nu'} and the corresponding nonlinear problem were studied in details by Lin and Xu \cite{LX},
and the enhanced dissipation is verified at some time scale $o(\frac{1}{\nu})$ for $\mathbb{Q}\omega$, called the non-shear part of the solution $\omega$ in \cite{LX}, where $\mathbb{Q}$ is the projection to the orthogonal complement of the kernel of  $-\sin y \partial_x ( I + \Delta^{-1})$ in $L^2$.  The core idea in \cite{LX} is to use the Hamiltonian structure of the operator $-\sin y \partial_x (I + \Delta^{-1})=J L$ with $J= -\sin y \partial_x$ and $L=I+\Delta^{-1}$ that naturally leads to the use of the weighted $L^2$ space $\langle L \cdot, \cdot\rangle_{L^2}$ in which $JL$ becomes antisymmetric. Then one can apply the RAGE Theorem for the estimate of the group $e^{t J L}$ and the argument of Constantin, Kiselev, Ryzhik, and Zlato${\rm \check{s}}$ \cite{CKRZ}, see also Zlato${\rm \check{s}}$ \cite{Z}, which study the enhanced dissipation for the advection-diffusion equations. 
Note that the inner product $\langle L \cdot, \cdot \rangle_{L^2}$ was a key tool also in obtaining the global stability of the Kolmogorov flows with arbitrary amplitude $a$; cf. \cite{Iu}.

The argument and the result of \cite{LX} are verified without any change also for the stability problem of the steady Kolmogorov flows \eqref{2dLinNS} - \eqref{def.L^nu}. 
However, for a deeper quantitative point of view, the spectral property of $\mathcal{L}^{\nu,a}$ requires a further study.
Indeed, the argument in \cite{LX} provides little information on the required smallness of $\nu$ to achieve the smallness of  $\frac{\| \mathbb{Q} \omega (t/\nu)\|_{L^2}}{\|\omega (0)\|_{L^2}}$, which depends on $t$ in an implicit way, even for the linear solution $\omega (t) = e^{t\mathcal{L}^{\nu,a}} \omega (0)$. 
In particular,  the question whether or not the smallness of  $\frac{\| \mathbb{Q} \omega (t/\nu^\beta)\|_{L^2}}{\|\omega (0)\|_{L^2}}$ holds for some $\beta\in (0,1)$, as solved in \cite{BW} with $\beta=\frac12$ for the model problem to \eqref{def.L^nu'}, has been a challenging problem; see Remark \ref{rem.LWZ} below. 

In this paper we will establish some resolvent estimate on the imaginary axis of the resolvent parameters for the linearization \eqref{def.L^nu} around the steady Kolmogorov flows. 
Our resolvent estimate is related to the pseudospectrum as in the work by Gallagher, Gallay, and Nier \cite{GaGaNi} of  the spectral analysis for large skew-symmetric perturbations of the Harmonic oscillator.
As a main result, we will verify the enhanced dissipation in a time scale $\mathcal{O}(\frac{1}{\sqrt{\nu}})$ for the linear flow $e^{t \mathcal{L}^{\nu,a}} \omega_0$; see \eqref{intro.main.Kol} below.
In particular, our result gives an affirmative answer to the problem numerically conjectured in \cite{BW}.
We expect that the similar enhanced dissipation will be true also for the linear flow generated by the evolution operator $\mathcal{L}^{\nu,a} (t)$ in \eqref{def.L^nu'}, which is still under investigation due to an obstacle from the time dependence of the operator.
The nonlinear problem \eqref{2dLinNS} can also be handled based on the linear estimates of this paper, but here we focus only on the linear problem.

By rescaling time as $t\mapsto \nu t$, one can rewrite the evolution problem $\partial_t \omega = \mathcal{L}^{\nu,a} \omega $ as,  by relabeling the variable and the unknown again as $t$ and $\omega$ respectively, 
\begin{align}\label{eq.rescaled.linear}
\partial_t \omega = \Delta \omega  - \frac a \nu \sin y \partial_x (I+\Delta^{-1}) \omega \,.
\end{align}
This problem is viewed  in the more abstract form 
\begin{equation}
\label{abstract.eq}
\partial_t \omega =(A - \alpha \Lambda) \omega\,,
\end{equation} 
where $\alpha>0$ is a large positive parameter, $A$ is a dissipative operator, and $\Lambda$ has a Hamiltonian structure. It will be worthwhile investigating the spectral property for such operators in the abstract level, which is handled in Section \ref{sec.general}.
The problem \eqref{eq.rescaled.linear} for the Kolmogorov flows is discussed in Section \ref{sec.Kol}, and we will  show the key estimate for the resolvent {\it with a rate on $\alpha=\frac{a}{\nu}$} (Theorem \ref{thm.Lambda.Kol.rate.2}), and then for the semigroup (Theorem \ref{thm.semigroup.Kol}).
In the original variables, our result in particular provides the bound for the semigroup  $e^{t\mathcal{L}^{\nu,a}}$ such as
\begin{align}\label{intro.main.Kol}
\| \mathbb{Q} e^{t \mathcal{L}^{\nu,a}} \omega_0\|_{L^2}\leq C e^{-c \sqrt{a \nu}  \, t} \| \mathbb{Q} \omega_0 \|_{L^2}\,, \qquad t \geq \frac{1}{\sqrt{a \nu}}\,,
\end{align}
see Corollary \ref{cor.thm.semigroup.Kol}. 
Here $C$ and $c$ are positive constants independent of $t$, $a$, $\omega_0$, and sufficiently small $\nu>0$.
This implies the enhanced dissipation in the time scale $\mathcal{O} (\nu^{-\frac12})$ for $0<\nu\ll 1$, that is much shorter than $o(\frac{1}{\nu})$ and provides a decay of $\frac{\|\mathbb{Q} \omega (t/\nu^\beta)\|_{L^2}}{\| \omega (0)\|_{L^2}}$ for $\beta=\frac12$ in the linear problem.
\begin{rem}{\rm It should be emphasized that the semigroup estimate \eqref{intro.main.Kol} is in fact new even for the model problem considered by \cite{BW} in which the nonlocal operator $\Delta^{-1}$ is dropped, though our result does not handle the time-dependent operator as in \eqref{def.L^nu'}. More precisely, the argument used in \cite{BW} provides the semigroup bound of the model problem in a weighted $H^1$ space whose norm has a dependence in $\nu$. In particular, the norm introduced in \cite{BW} involves the term $\nu^{-\frac14} \| M_{\cos y} f\|_{L^2}$, in addition to the usual $L^2$ norm. 
As a result, the estimate obtained in \cite{BW} shows that $\| M_{\cos y} \omega (t) \|_{L^2}$ becomes small  in a time scale $O(\nu^{-\frac12})$ for the model problem, while in order to achieve the dissipation in the usual $L^2$ norm it seems that one needs to take a slightly longer  time scale, e.g., $O(\nu^{-\frac12} |\log \nu|)$. Our result \eqref{intro.main.Kol} gives the dissipation in the $L^2$ norm exactly in a time scale $O(\nu^{-\frac12})$. Our proof is based on the detailed resolvent analysis and is very different from the approach in \cite{BW}.}
\end{rem}

Next let us briefly mention  a topic which is  closely related to the present work as another example of \eqref{abstract.eq}: the asymptotic stability problem of the Lamb-Oseen vortex. By working on $\mathbb M=\mathbb R^2$, it is known that there exists a family of self-similar solutions to the vorticity equation \eqref{2dVort} given by 
\begin{align}\label{def.Oseen}
\Omega(x,y,t)=\frac{\gamma}{\nu t}G(\frac{x}{\sqrt{\nu t}}, \frac{y}{\sqrt{\nu t}}) ,\quad\mbox{and }\quad U(x,y,t)=\frac{\gamma}{\sqrt{\nu t}}V^G(\frac{x}{\sqrt{\nu t}}, \frac{y}{\sqrt{\nu t}})\,,
\end{align}
where the profiles are  $G(\xi)=\frac1{4\pi}e^{-|\xi|^2/4}$ and  $V^G(\xi)=\frac1{2\pi}\frac{\xi^\perp}{|\xi|^2}(1-e^{-|\xi|^2/4})$. The constant $\gamma=\int_{\R^2}\Omega(x,y,t)\, \dd x \dd y$ is the circulation at infinity of the flow. By the significant work of Gallay and Wayne \cite{GW} it is known that this solution is the only forward self-similar solution to \eqref{2dNS} in $\R^2$ with an integrable vorticity. This solution is called the Lamb-Oseen vortex. 
It is well known that, through a suitable similarity transformation, the asymptotic stability of the Lamb-Oseen vortex is equivalent with the two dimensional stability of the Burgers vortex, which is a stationary solution to the three dimensional Navier-Stokes equations in the presence of the axisymmeric linear strain.
The reader is referred to a recent review article \cite{GaMa} by Gallay and the second author of this paper  about the research on the stability of the Burgers vortex.
The two dimensional linearized problem for the Burgers vortex with circulation $\alpha$ is given by    
\begin{equation}
\label{LambO1}
\partial_\tau \omega = (\mathcal L-\alpha\Lambda)\omega \,, \qquad \tau>0\,, \quad \xi \in \R^2\,,
\end{equation}
where $\mathcal L \omega =\Delta \omega  +\frac12\xi \cdot \nabla \omega + \omega $, and $\Lambda \omega = (V^G\cdot\nabla)\omega +(K_{BS}*\omega\cdot\nabla)G$. Here $\Delta$ and $\nabla$ are now about the variables $\xi = (\xi_1,\xi_2)$.
In the weighted $L^2$ space $L^2(\R^2; \frac{\dd \xi}{G})$, the operator $-\mathcal{L}$ is nonnegative self-adjoint with compact resolvent, and $\Lambda$ becomes antisymmetric as proved in \cite{GW}. 
Hence the linear analysis falls into the analysis of the operator of the form \eqref{abstract.eq}. 
In the space $L^2(\mathbb R^2,G^{-1}d\xi)$ with zero mass condition, we have $-\mathcal{L}\geq \frac12$, and thus, the antisymmetry of $\Lambda$ provides $\frac12$ spectral gap for $\mathcal{L}-\alpha \Lambda$ for any $\alpha$. This yields the linear stability with a uniform estimate in $\alpha$. However, this simple argument does not provide further informations for the fast rotation case $|\alpha|\gg 1$, at the time when numerical and experimental evidence suggest that the basin of attraction should be $\alpha$-dependent, at least ``away" from the kernel of the operator $\Lambda$. In \cite{Ma} the second author of this paper verifies a behavior of the pseudospectral bound but without the information on the rate about $\alpha$. 
On the other hand, in \cite{GaGaNi} and Deng \cite{De1} simplified model operators are studied in details, where the main simplification is dropping the nonlocal term $(K_{BS}*w, \nabla )G$, and the optimal dependence on $\alpha$ of the pseudospectral bound that decays like $|\alpha|^\frac13$ is obtained for these model operators. The same result is proved for the full linearized operator $\mathcal{L}-\alpha \Lambda$ in Deng \cite{De2} but in a smaller subspace than the orthogonal complement of ${\rm Ker}\, \Lambda$.
Very recently, Li, Wei, and Zhang \cite{LWZ} gave a sharp pseudospectral bound as well as the spectral bound of $L-\alpha \Lambda$ in the orthogonal complement of  ${\rm Ker}\, \Lambda$, and this result is applied to the nonlinear problem by Gallay \cite{Ga2017}. 

\begin{rem}\label{rem.LWZ}{\rm In Li, Wei, and Zhang \cite{LWZ} the key elegant idea is to introduce the wave operator which converts the original skew-symmetric operator $\Lambda$, containing a nonlocal term that leads to an essential difficulty, into a skew-symmetric operator for which the nonlocal operator is removed and hence the approach of \cite{GaGaNi} is applied. As announced in \cite{LWZ}, it is recently shown by Wei, Zhang, and Zhao \cite{WZZ} (see also Li, Wei, and Zhang \cite{LWZ2} for $3$D problem) that this approach for the Lamb-Oseen operator can be applied also for the estimate of the enhanced dissipation around the Kolmogorov flows and the optimal enhanced dissipation as in \eqref{intro.main.Kol} is obtained together with the algebraic dissipation in the time scale $O(\nu^{-\frac13})$ for the velocity field. We note that our approach for \eqref{intro.main.Kol} or Theorem \ref{thm.improve.Kol.intro} below is different and independent of \cite{LWZ, WZZ}, and in particular, does not rely on the construction of the wave operator.
}
\end{rem}

To summarize, the above two examples of the Kolmogorov flows and the Lamb-Oseen vortex show that to measure the basin of attraction, it is important to obtain a pseudospectral bound as sharp as possible for the operator in the abstract form given in \eqref{abstract.eq}.  We also note that the enhanced dissipation is one of the important subjects in fluid mechanics, and recently,  significant progress has been achieved around some class of simple flows such as the Couette flow; see, e.g., \cite{BeGerMa1, BeGerMa2, BeMaVi, BeGerMa3}.

This paper consists of two parts. The first one is an abstract result, in which the spectral properties of some class of non self-adjoint operators are established.
The other one is the application of the abstract result to the linearized operator for the Kolmogorov flows. As for the abstract part, we consider the operator in a Hilbert space $X$ of the form
\begin{align}\label{def.L_alpha.intro}
L_\alpha  =  A - \alpha \Lambda\,,
\end{align}
where $-A$ is positive self-adjoint with compact resolvent, $\alpha\in \R$, and $\Lambda$ is a densely defined closed linear operator relatively compact to $A$.  
For later use we set $\hat{\Lambda}$ by the relation 
\begin{align}
\Lambda = i \hat{\Lambda}\,.
\end{align}
We denote by $D_X (A)$ the domain of $A$ in $X$. 
We are interested in the spectral property of $L_\alpha$ for large $|\alpha|$. 
Since the effect of $\alpha$ is absent for functions in ${\rm Ker} \,\hat{\Lambda}$
it is natural to introduce the orthogonal projections 
\begin{align}\label{def.Q}
\mathbb{Q} : X \rightarrow Y :=  \big ({\rm Ker}\, \hat{\Lambda} \big ) ^\bot\,,
\end{align}
where $K^\bot$ denotes the orthogonal complement space in $X$ for a given closed subspace $K$.
We are interested in the estimate of  $\mathbb{Q}e^{t L_\alpha}$ for large $\alpha$. 
Since the semigroup $e^{t L_\alpha}$ is expressed  in terms of the resolvent of $L_\alpha$ the problem is reduced  to the estimate of  $\mathbb{Q} (i \lambda- L_\alpha)^{-1}$ when $i\lambda$ belongs to the resolvent set of $L_\alpha$.
When the invariance $\mathbb{Q}A\subset A\mathbb{Q}$ holds, which will be assumed in this paper, 
the estimate of $\mathbb{Q} (i \lambda- L_\alpha)^{-1}$ is reduced to the resolvent analysis of the operator $\mathbb{Q}L_\alpha$ in $Y$,
%Our particular interest is the estimates for the quantity
%\begin{align}
%\sup_{\lambda\in \R} \| \mathbb{Q} (i\lambda -  L_\alpha )^{-1} \|_{X\rightarrow X} \label{def.pseudo.1}
%\end{align}
%for large $|\alpha|$, when $i\lambda$ belongs to the resolvent set of $L_\alpha$ in $X$. 
%Closely related to \eqref{def.pseudo.1} it is useful to introduce the operator $\mathbb{Q}L_\alpha$ in $Y$
which is realized as 
\begin{align}
\begin{split}
D_Y (\mathbb{Q}L_\alpha) & = D_Y (\mathbb{Q}A) := D_X(A) \cap Y\,, \\
\mathbb{Q}L_\alpha u & =  \mathbb{Q} A u - i \alpha \mathbb{Q}  \hat{\Lambda} u\,, \qquad u \in  D_Y (\mathbb{Q}L_\alpha) \,.
\end{split}
\end{align}
Indeed, we have $\mathbb{Q} (i\lambda - L_\alpha)^{-1} f =(i\lambda - \mathbb{Q} L_\alpha)^{-1}f$ for $f\in Y$ when $\mathbb{Q}A\subset A \mathbb{Q}$.
In order to obtain the estimate of $\mathbb{Q}e^{t L_\alpha}$ or $e^{t\mathbb{Q}L_\alpha}$,
the following quantity plays an essential role: 
\begin{align}\label{def.pseudo}
\Psi_Y (\alpha; \mathbb{Q}L_\alpha) = \bigg ( \sup_{\lambda\in \R} \| (i\lambda - \mathbb{Q} L_\alpha )^{-1} \|_{Y\rightarrow Y} \bigg )^{-1}= \bigg ( \sup_{\lambda\in \R} \| (i\lambda + \mathbb{Q} L_\alpha )^{-1} \|_{Y\rightarrow Y} \bigg )^{-1}\,.
\end{align}
The quantity \eqref{def.pseudo} was introduced in \cite{GaGaNi}, where the basic pseudospectral property and the relation with the semigroup estimate are also presented. 
For convenience we call \eqref{def.pseudo} the pseudospectral bound of $\mathbb{Q}L_\alpha$. 
In our framework the operator $\Lambda$ is not necessarily antisymmetric, but instead, is assumed to possess a Hamiltonian structure; see Assumption \ref{assum.Lambda} in Section \ref{sec.general}. This structural assumption is of course motivated by the application to the Kolmogorov flows.
There are two theorems in the abstract part. The first one is a pseudospectral bound without  a concrete dependence on  $\alpha$  (Theorem \ref{thm.abstract.1}).
The argument in Theorem \ref{thm.abstract.1} shares some common features with the argument in \cite{LX}. 
 While the result of \cite{LX}  is based on the RAGE theorem, 
  our argument  is much  more elementary, though Theorem \ref{thm.abstract.1} does not necessarily give a stronger result than \cite{LX}. 
The second  result  of the abstract part has   a concrete dependence on $\alpha$  (see Theorem \ref{thm.abstract.2}),  under additional assumptions on $\hat{\Lambda}$. 
The key additional condition is Assumption \ref{assum.Lambda.rate} which imposes some coercive estimate for $\mu - \hat{\Lambda}$ with $\mu \in \R$ by allowing a presence of the term yielding a ``loss of derivative'' but with a small factor in front. This derivative loss with a small prefactor is controlled by the smoothing effect of $A$ at the end, and this balance determines the rate in $\alpha$ for the pseudospectral bound.
In Assumption \ref{assum.Lambda.rate} another key condition is imposed on the cross term $\Im \langle A u, \hat{\Lambda} u \rangle_X$, which is useful in achieving the resolvent estimate with a sharper dependence on $\alpha$. This type of condition fits with the case when $A$ is of the form $A=-T^*T$ and is related with the commutator $[T,\hat{\Lambda}]$, and thus, our approach is highly motivated by the work of \cite{Vi, GaGaNi}.

As an application of the abstract result,  we study in Section \ref{sec.Kol} the rescaled version of  the linear operator \eqref{def.L^nu}, i.e., the problem \eqref{eq.rescaled.linear}.
By taking the Fourier series in $x$, the key is to analyze the operator only in the $y$ variable in the space $L^2 (\T)$:
\begin{align*}
L_{\alpha,l} = A_l - i \alpha l \hat{\Lambda}_l\,, \qquad \hat{\Lambda}_l = M_{\sin y} (I + A_l^{-1})\,.
\end{align*} 
Here $A_l=\partial_y^2 - l^2$, $\alpha = \frac{a}{\nu}$, $M_{\sin y} f = \sin y\, f$, and $l\in \Z\setminus\{0\}$.
The operator $\hat{\Lambda}_l$ has a nontrivial kernel only when $l=\pm 1$ which is spanned by the constant functions, and thus, the projection $\mathbb{Q}_l: L^2 (\T) \rightarrow L^2(\T)$ is defined by 
\begin{align*}
\mathbb{Q}_l f=f \quad \text{for}~~|l|\geq 2\,, \qquad \mathbb{Q}_l f = f - \frac{1}{2\pi} \int_0^{2\pi} f \dd y \quad \text{for}~~|l|=1\,.
\end{align*}
The main effort is to check the coercive estimates described in Assumption \ref{assum.Lambda.rate} for $\mu-\hat{\Lambda}_l$ with $\mu\in \R$ which is essential to achieve the pseudospectral bound with a rate in $\alpha$.  We shall verify Assumption \ref{assum.Lambda.rate} by analyzing the ODE corresponding to the operator $\mu-\hat{\Lambda}_l$. The main difficulty here is the presence of the term $A_l^{-1}$ in $\hat{\Lambda}_l$, which makes the problem nonlocal and also leads to some  lack of invariance, namely the fact that  $(I-\mathbb{Q}_l)\hat{\Lambda}_l\ne 0$ when $l=\pm 1$. This loss of invariance is due to the absence of the symmetry of $\hat{\Lambda}_l$, and gives rise to an additional nonlocality coming from the projection $\mathbb{Q}_l$.
Therefore, we have to deal with two nonlocalities; the one in  $A_l^{-1}$ and the one in  $\mathbb{Q}_l$. 
For a given $\mu\in \R$ the point $y\in \T$ satisfying $\sin y = \mu$ is called a critical point of this problem. 
The difficulty coming from the nonlocality of $A_l^{-1}$ is significant when the critical points are degenerate, and this corresponds to the case when $|\mu|$ is around $1$ in the analysis of $\mu-\hat{\Lambda}_l$. The core part of the analysis is Lemma \ref{lem.Lambda.Kol.rate.2} which deals with this singularity. The key idea is to use a contradiction argument, which enables us to focus on the functions concentrating around the critical points, for which the nonlocal term essentially becomes a small order since the operator $A_l^{-1}$ has a smoothing effect.
On the other hand, the influence of the projection $\mathbb{Q}_l$ becomes relevant only when $\mu$ is close to $0$ in the analysis of $\mu-\hat{\Lambda}_l$, for $I-\mathbb{Q}_l$ is the projection to the kernel of $\hat{\Lambda}_l$. 
As a result, these two kinds of difficulty related to $A_l^{-1}$ and to $\mathbb{Q}_l$ appear in different parameter regimes of $\mu$, and thus we can handle them   separately. 
%The  two conditions \eqref{est.assum.Lambda.rate'} and \eqref{est.assum.Lambda.rate} in Assumption \ref{assum.Lambda.rate} in fact reflect such a situation for the Kolmogorov flows. 
After establishing the key coercive bounds of $\mu-\hat{\Lambda}_l$, which are stated in Proposition \ref{thm.Lambda.Kol.rate.1}, the resolvent estimate for $L_{\alpha,l}$ is obtained in Theorem \ref{thm.Lambda.Kol.rate.2} by applying the abstract result in Section \ref{sec.general} and also by using a specific property of the trigonometric functions.
For convenience,  it will be worthwhile stating our resolvent estimate for the Kolmogorov flow 
 in this introductory section: 
\begin{thm}\label{thm.improve.Kol.intro} There exist $C, \alpha_0>0$ such that the following statement holds for all $\alpha\in \R$ with $|\alpha|\geq \alpha_0$. Let $\lambda\in \R$ and $l\in \Z\setminus \{0\}$. Then
\begin{equation}\label{est.thm.improve.Kol.intro.1}
\| (i\lambda + \mathbb{Q}_l L_{\alpha,l})^{-1}\|_{Y_l \rightarrow Y_l}   \leq 
\begin{cases}
& \displaystyle \frac{C}{|\alpha l| \, (|\frac{\lambda}{\alpha l}|-1)} ~~\qquad {\rm if}~~|\frac{\lambda}{\alpha l}|>1+\frac{1}{|\alpha l|^\frac12}\,,\\
& \displaystyle \frac{C}{|\alpha l |^\frac12} \qquad \qquad \qquad {\rm if}~~1-\frac{1}{|\alpha l|^\frac12} < |\frac{\lambda}{\alpha l}| \leq 1 + \frac{1}{|\alpha l|^\frac12}\,,\\
& \displaystyle \frac{C}{|\alpha l|^\frac23 (1-|\frac{\lambda}{\alpha l}|)^\frac13} \qquad {\rm if}~~ |\frac{\lambda}{\alpha l}| \leq 1 - \frac{1}{|\alpha l|^\frac12}\,.
\end{cases}
\end{equation}
Here $Y_l = \mathbb{Q}_l L^2 (\T)$.
\end{thm}
The estimate \eqref{est.thm.improve.Kol.intro.1} actually gives more detailed information on the spectrum of $\mathbb{Q}L_{\alpha,l}$ than the pseudospectral bound defined by \eqref{def.pseudo}, and seems to be considerably sharp in view of the degeneracy of the critical points.  
In fact, we observe that the critical points become degenerate when $|\frac{\lambda}{\alpha l}|\sim 1$, and \eqref{est.thm.improve.Kol.intro.1} claims that the rate is $O(|\alpha l|^{-\frac12})$ around this case. When $|\frac{\lambda}{\alpha l}|$ is less than $1$, the critical points are nondegenerate and the rate is improved as $O(|\alpha l|^{-\frac23})$. Note that these rates, $O(|\alpha l|^{-\frac12})$ and $O(|\alpha l|^{-\frac23})$, depending on the degeneracy of the critical points, are compatible with the result in \cite{GaGaNi} and hence they are optimal if the nonlocal term $M_{\sin y} A_l^{-1}$ is dropped from $\hat{\Lambda}_l$. Additional remark for the local operator $A_l-i\alpha l M_{\sin y}$ is that, near the critical point $\sin y \sim \mu$, the operator is modeled by the complex Airy operator $\partial_y^2 + i y$ when $\mu$ is away from $\pm 1$ (nondegenerate case) that is responsible for the exponent $2/3$, while it is modeled by $\partial_y^2 \pm i y^2$ when $\mu$ is close to $\pm 1$ (degenerate case), resulting the exponent $1/2$.
Finally, if $|\frac{\lambda}{\alpha l}|$ is larger than $1$,  the critical points are no longer present, resulting in the rate $O(|\alpha l|^{-1})$.  

This paper is organized as follows. In Section \ref{sec.general} we discuss the problem in an abstract framework.
Section \ref{sec.Kol} is devoted to the study of the linearized problem for the Kolmogorov flows.
The main results in Section \ref{sec.Kol} are Theorem \ref{thm.semigroup.Kol} and its Corollary \ref{cor.thm.semigroup.Kol} for the estimate of the semigroup $\{e^{t\mathcal{L}^{\nu,a}} \}_{t\geq 0}$.
In Section 4 we also consider the application to the Lamb-Oseen vortices by omitting some details of the proof since the argument is similar to the case of the Kolmogorov flow. Section 4 provides alternative approach for the result of \cite{LWZ}.

\section{Abstract result}\label{sec.general}

In this section we establish the abstract result in obtaining the resolvent estimate for the operator \eqref{def.L_alpha.intro}, by taking into account the application to the stability of the Kolmogorov flows.
In fact, to prove the estimate stated in Theorem \ref{thm.improve.Kol.intro} requires a rather complicated and long argument, and thus, the abstract result is useful in understanding the basic strategy.
First we state the basic assumption on $A$.
\begin{assum}\label{assum.A} The operator $A: D_X(A)\subset X \rightarrow X$ is self-adjoint in $X$ with compact resolvent, and $-A$ is positive and satisfies 
\begin{align}\label{est.assum.A}
\langle -A u, u\rangle _X\geq \|u\|_X^2 \,,  \qquad u\in D_X(A)\,.
\end{align}
\end{assum}

\begin{rem}{\rm One can extend the result of this section to more general class of $A$ such that $-A$ is  {\it m-sectorial} satisfying the positivity $\Re \langle -Au, u \rangle_X\geq \| u\|_X^2 + C |\Im \langle Au, u \rangle_X|$ with compact resolvent, by slightly modifying the assumption on $\hat{\Lambda}$. But for simplicity we focus on the case when $A$ is self-adjoint.
}
\end{rem}

Next we state the conditions on the relation between $A$ and $\hat{\Lambda}$. 
\begin{assum}\label{assum.Lambda} {\rm (i)} $\hat{\Lambda}$ is a densely defined closed operator and is relatively compact to $-A$ in $X$. 

\noindent {\rm (ii)} Set $Y=({\rm Ker}\, \hat{\Lambda})^\bot$, the orthogonal complement space of ${\rm Ker}\, \hat{\Lambda}$ in $X$, and let $\mathbb{Q}: X\rightarrow Y$ be the orthogonal projection. Then $\mathbb{Q} A \subset A \mathbb{Q}$, $D_X(A) \subset D_X(\hat{\Lambda}^*)$, and there exists a positive constant $C$ such that 
\begin{align*}
|\langle \hat{\Lambda}u\,, u\rangle_X | + |\Im \langle - A u,  \hat{\Lambda} u\rangle_X| \leq C \langle -Au, u\rangle_X\,, \qquad u\in D_X(A)\,.
\end{align*}

\noindent {\rm (iii)} There exist closed symmetric operators $B_1$, $B_2$, and positive constants $c_1$ and $C$ such that $\hat{\Lambda}=B_1 B_2$, ${\rm Ker}\, \hat{\Lambda} = {\rm Ker}\, B_2$, and 
\begin{align}
\| B_2 u \|_X^2 & \leq  C \| (-A)^\frac12 u\|_X^2 \,, \qquad u\in D_X(A) \,,\label{est.assum.Lambda.0}\\
\langle u\,,  B_2 u \rangle_X & \geq c_1 \| u\|_X^2   \,,  \qquad \qquad u\in D_X(A)\cap Y\,,\label{est.assum.Lambda.1}\\
\Re \langle - A u, B_2 u \rangle_X   & \geq  c_1 \| (- A)^\frac12 u \|_X^2  \qquad u \in D_X (A) \cap Y \,.\label{est.assum.Lambda.2}
\end{align} 
\end{assum}

\begin{rem}{\rm Assumption \ref{assum.Lambda} (iii) states that $\hat{\Lambda}$ has a property similar to a closed symmetric operator. Indeed, if $\hat{\Lambda}$ is closed symmetric then it suffices to take $B_1=\hat{\Lambda}$ and $B_2=\mathbb{Q}$. 
}
\end{rem}

Let us denote by ${\rm Ran} \, \hat{\Lambda}$ the range of $\hat{\Lambda}$, i.e, 
${\rm Ran}\, \hat{\Lambda} = \{f\in X~|~ f=\hat{\Lambda} g ~\text{ for some }~ g\in D_X(\hat{\Lambda})\}$. 

\begin{assum}\label{assum.Lambda'} {\rm (i)} ${\rm Ker}\, \hat{\Lambda} \subset D_X(\hat{\Lambda}^*)$.

\noindent {\rm (ii)} ${\rm Ker}\, \hat{\Lambda} \cap {\rm Ran}\, \hat{\Lambda}=\{0\}$.

\noindent {\rm (iii)} $\hat{\Lambda}$ does not possess eigenvalues in $\R\setminus\{0\}$.
\end{assum}

\begin{rem}{\rm (1) Assumption \ref{assum.Lambda'} (i) is imposed in order to justify the formal computation.

\noindent (2) If $\hat{\Lambda}$ is closed symmetric then Assumption \ref{assum.Lambda'} (ii) holds. Indeed, it suffices to use the orthogonal decomposition $X={\rm Ker}\, \hat{\Lambda}^* \oplus \overline{{\rm Ran}\, \hat{\Lambda}^{**}} = {\rm Ker}\, \hat{\Lambda}^* \oplus \overline{{\rm Ran}\, \hat{\Lambda}}$; then for any $f$ we have the corresponding decomposition $f=\varphi + \psi$, and $\psi$ is approximated by $\{\psi_n\}$ such that $\psi_n = \hat{\Lambda}\phi_n$. Then for any $u\in {\rm Ker}\, \hat{\Lambda} \cap {\rm Ran}\, \hat{\Lambda}$ we have 
\begin{align*}
\langle u, f \rangle_X = \langle u, \varphi\rangle_X + \langle u, \psi_n \rangle_X + \langle u, \psi-\psi_n \rangle_X = \langle u, \psi-\psi_n \rangle_X \rightarrow 0\qquad n\rightarrow \infty.
\end{align*}
Hence $u=0$.
}
\end{rem}

First we state the abstract result of the spectral behavior of $\mathbb{Q} L_\alpha$ with $\alpha \in \R$ in the limit $|\alpha|\rightarrow \infty$,
but without any information on the rate of convergence.
\begin{thm}\label{thm.abstract.1} Suppose that Assumptions \ref{assum.A}, \ref{assum.Lambda}, and \ref{assum.Lambda'} hold. Let $\sigma_Y (\mathbb{Q}L_\alpha)$ be the set of the spectrum of $\mathbb{Q}L_\alpha$, $\alpha \in \R$, in $Y$. Then we have 
\begin{align}\label{thm.abstract.1.1} 
\lim_{|\alpha| \rightarrow \infty} \sup_{\mu \in \sigma_Y (\mathbb{Q}L_\alpha)} \Re \mu =-\infty\,,
\end{align}
and 
\begin{align}
\lim_{|\alpha| \rightarrow \infty} \sup_{\lambda\in \R} \| (i\lambda -\mathbb{Q}L_\alpha)^{-1}\|_{Y\rightarrow Y} & =0\,.\label{thm.abstract.1.2} 
\end{align}
Moreover, for sufficiently large $|\alpha|$ the set $\{\zeta\in \C~|~\Re \zeta >-1 \}$ is contained in the resolvent set of $L_\alpha$ in $X$, and we have 
\begin{align}
\lim_{|\alpha| \rightarrow \infty} \sup_{\lambda\in \R} \| \mathbb{Q} (i\lambda - L_\alpha)^{-1}\|_{X\rightarrow X} & =0\,.\label{thm.abstract.1.3} 
\end{align} 
\end{thm}

\noindent {\it Proof.} The proof consists of several steps. Without loss of generality we may assume that $\alpha>0$.

\noindent Step 1: {\it The operator $\mathbb{Q}A$ in $Y$ is a closed linear operator with compact resolvent.} 

\noindent This follows directly from the invariance $\mathbb{Q} A \subset A \mathbb{Q}$ and Assumption \ref{assum.A}. We denote by $A|_{Y}$ the restriction of $A$ to $Y$ with the domain $D_Y(A|_Y)=D_X(A) \cap Y$.

\

\noindent Step 2: {\it $\sigma_Y (\mathbb{Q}L_\alpha) \subset \{\mu \in \C~|~\Re \mu <0 \}$.}

\noindent Let $c_1>0$ be the number in \eqref{est.assum.Lambda.2}. We have already seen that  $-A|_Y$ is bounded from below in $Y$ and has a compact resolvent in $Y$. Then, since $\hat{\Lambda}$ is relatively compact to $A$ by the assumption, we see that $\mathbb{Q}L_\alpha = A|_Y - i \alpha \mathbb{Q} \hat{\Lambda}$ is also a closed operator with compact resolvent in $Y$.
Thus, the spectrum of $\mathbb{Q}L_\alpha$ consists of isolated eigenvalues with finite multiplicities.
Let $\mu\in \C$ be an eigenvalue of $\mathbb{Q}L_\alpha$ in $Y$ and let $u\in D_X(A)\cap Y$ be an associated eigenfunction such that $\|u\|_X=1$. Note that $-\mu u = -\mathbb{Q}L_\alpha u$ holds. By taking the inner product with $B_2u$, we have 
\begin{align}\label{proof.thm.abstract.1.2} 
\begin{split}
-\mu \langle u\,, B_2 u \rangle_X & = \langle -\mathbb{Q} A  u\,,  B_2 u\rangle _X + i \alpha \langle \mathbb{Q} \hat{\Lambda} u\,, B_2 u\rangle _X\\
& =  \langle -\mathbb{Q} A  u\,,  B_2 u\rangle _X + i \alpha \langle B_1 B_2 u\,, \mathbb{Q} B_2 u\rangle _X\,.
\end{split}
\end{align}
Here we have used $\hat{\Lambda}=B_1B_2$ by Assumption \ref{assum.Lambda} (ii). 
Moreover, we verify that $(I-\mathbb{Q}) B_2 = 0$ since $B_2$ is closed symmetric and ${\rm Ker}\, \hat{\Lambda}={\rm Ker}\, B_2$. Hence we have $\mathbb{Q}B_2=B_2$, which implies $\langle B_1 B_2 u\,, \mathbb{Q} B_2 u\rangle _X=\langle B_1 B_2 u\,, B_2 u\rangle _X$ and $\langle -\mathbb{Q} A  u\,,  B_2 u\rangle _X = \langle - A  u\,,  \mathbb{Q} B_2 u\rangle _X =\langle -A  u\,,  B_2 u\rangle _X$.
Therefore, since $B_1$ is closed symmetric in $X$, the real part of \eqref{proof.thm.abstract.1.2} yields 
\begin{align}\label{proof.thm.abstract.1.3} 
-(\Re \mu ) \langle u\,, B_2 u \rangle_X & = \Re \langle -Au\,, B_2 u \rangle _X \nonumber \\
& \geq c_1 \| (-A)^\frac12 u\|_X^2  \geq c_1\,.
\end{align} 
Here we have used the assumptions \eqref{est.assum.Lambda.2} and \eqref{est.assum.A} with $\|u\|_X=1$.
Hence, \eqref{est.assum.Lambda.1} yields
\begin{align}\label{proof.thm.abstract.1.4} 
\Re (\mu)  \leq -\frac{c_1}{\langle u, B_2 u\rangle_X} <0\,.
\end{align}
Thus, $\sigma_Y (\mathbb{Q}L_\alpha) \subset \{\mu \in \C~|~\Re \mu < 0\}$ for all $\alpha \geq 0$.

\

\noindent  Step 3: {\it The spectral limits \eqref{thm.abstract.1.1} and \eqref{thm.abstract.1.2} hold.}

\noindent By Step 2 it suffices to show 
\begin{align}\label{proof.thm.abstract.1.5} 
\lim_{\alpha\rightarrow \infty} \sup_{\lambda\in \R} \| (i\lambda  -\mathbb{Q}L_\alpha)^{-1}\|_{Y\rightarrow Y} =0\,.
\end{align}
Suppose that \eqref{proof.thm.abstract.1.5} does not hold. Then there exist $\delta>0$, $\{\alpha_n\}\subset \R_+$, $\{\lambda_n\}\subset \R$, and $f_n\in Y$ with $\|f_n \|_X=1$ such that $\alpha_n\rightarrow \infty$, and $\|(i\lambda_n  - \mathbb{Q}L_\alpha)^{-1} f_n \|_X\geq \delta$. Set $u_n = (i\lambda_n  - \mathbb{Q}L_\alpha)^{-1} f_n$, which solves 
\begin{align}\label{proof.thm.abstract.1.6} 
i\lambda_n  u_n - \mathbb{Q}A u_n + i\alpha_n \mathbb{Q} \hat{\Lambda} u_n =f_n\,.
\end{align} 
By taking the inner product with $B_2 u_n$ in the above equation, we obtain 
\begin{align}\label{proof.thm.abstract.1.7} 
i\lambda_n  \langle u_n\,, B_2 u_n\rangle_X +\langle -A u_n \,, B_2u_n \rangle_X + i\alpha_n \langle B_1B_2 u_n\,, B_2u_n \rangle_X = \langle f_n\,, B_2 u_n \rangle_X\,.
\end{align} 
Here we have used $\mathbb{Q} B_2 = B_2$ as observed in Step 2.
Since $B_1$ and $B_2$ are symmetric, the real part of \eqref{proof.thm.abstract.1.7} yields
\begin{align*}
\Re \langle -A u_n\,, B_2 u_n \rangle_X = \Re \langle f_n\,, B_2 u_n \rangle_X\,,
\end{align*}
and then the assumptions  \eqref{est.assum.Lambda.0}, \eqref{est.assum.Lambda.1}, and \eqref{est.assum.Lambda.2} imply
\begin{align*}
c_1 \langle -A u_n\,, u_n \rangle_X  \leq \Re \langle f_n\,, B_2 u_n \rangle_X 
& \leq \|f_n \|_X \| B_2 u_n\|_X \leq C \| (-A)^\frac12 u_n\|_X\,.
\end{align*}
Thus we obtain the uniform bound 
\begin{align}\label{proof.thm.abstract.1.8} 
\sup_n \| (-A)^\frac12 u_n \|_X <\infty\,.
\end{align}
Now we recall that, since $-A$ is positive self-adjoint with compact resolvent, $(-A)^\frac12$ also has a compact resolvent ( \cite[Theorem V-3.49]{Ka}). Since \eqref{proof.thm.abstract.1.8} implies the uniform bound of $\| (-A)^\frac12 u_n \|_X$, $\{u_n\}$ is compact in $X$, and thus, in $Y$.
Then there exists a subsequence of $\{u_n\}$, denoted again by $\{u_n\}$, which strongly converges to some $u_\infty\in Y$ and satisfies $\| (-A)^\frac12 u_\infty \|_X \leq \sup_n \| (-A)^\frac12 u_n\|_X<\infty$. 
By the strong convergence we have $\|u_\infty\|_X\geq \delta$, so $u_\infty\in Y$ is nontrivial. 
Let us go back to \eqref{proof.thm.abstract.1.6}, and take the inner product with $u_n$. Then we have 
\begin{align*}
i\frac{\lambda_n}{\alpha_n} \| u_n \|_X^2  + \frac{1}{\alpha_n} \langle -A u_n\,, u_n\rangle_X + i \langle \hat{\Lambda}u_n\,, u_n \rangle_X = \frac{1}{\alpha_n} \langle f_n\,, u_n\rangle_X
\end{align*}
and the imaginary part of this identity yields the bound 
\begin{align*}
|\frac{\lambda_n}{\alpha_n}| \, \|u_n\|_X^2 \leq  |\langle \hat{\Lambda} u_n\,, u_n \rangle_X|  + |\Im \langle f_n\,, u_n\rangle_X | 
& \leq C \| (-A)^\frac12 u_n\|_X^2   + \| f_n\|_X \| u_n\|_X \nonumber \\
& \leq C \big (\sup_n  \| (-A)^\frac12 u_n\|_X^2 +1 \big )<\infty\,.
\end{align*}
Since $\|u_n\|_X\geq \delta$ we have the uniform bound 
\begin{align*}
\sup_n |\frac{\lambda_n}{\alpha_n}|<\infty\,.
\end{align*}
Set $\mu_n = \frac{\lambda_n}{\alpha_n}$. By taking a suitable subsequence we may assume that $\mu_n$ converges to some $\mu_\infty \in \R$. For any $\varphi\in D_X(A)$ we have from \eqref{proof.thm.abstract.1.6},
\begin{align*}
i\mu_n \langle u_n\,, \varphi\rangle_X  - \frac{i}{\alpha_n} \langle u_n\,, A \mathbb{Q}\varphi\rangle_X  + i \langle u_n\,, \hat{\Lambda}^* \mathbb{Q} \varphi \rangle_X = \frac{1}{\alpha_n} \langle f_n\,, \varphi \rangle_X\,,
\end{align*}
and by taking the limit $n\rightarrow \infty$, we have 
\begin{align*}
i\mu_\infty \langle u_\infty\,, \varphi\rangle_X + i\langle u_\infty\,, \hat{\Lambda}^* \mathbb{Q} \varphi \rangle_X =0\,,
\end{align*}
and thus, 
\begin{align*}
\langle u_\infty\,, \hat{\Lambda}^* \varphi \rangle_X =- \mu_\infty \langle u_\infty\,, \varphi\rangle_X + \langle u_\infty\,, \hat{\Lambda}^* (I-\mathbb{Q}) \varphi \rangle_X\,.
\end{align*}
Note that $(I-\mathbb{Q}): X \rightarrow {\rm Ker}\, \hat{\Lambda} \subset  D_X(\hat{\Lambda}^*)$ by the assumption, and thus, $\hat{\Lambda}^*(I-\mathbb{Q})$ defines a bounded operator in $X$ by the closed mapping theorem.
Then, since $D_X(A)$ is dense in $X$ this identity holds for all $\varphi \in D_X(\hat{\Lambda}^*)$, which implies $u_\infty\in D_X(\hat{\Lambda})$ and 
\begin{align*}
\hat{\Lambda} u_\infty = - \mu_\infty u_\infty + (I-\mathbb{Q}) \hat{\Lambda} u_\infty \,.
\end{align*}
This shows $\hat{\Lambda}u_\infty \in D_X(\hat{\Lambda})$ and 
\begin{align}
\label{proof.thm.abstract.1.9} 
\hat{\Lambda}^2 u_\infty = - \mu_\infty \hat{\Lambda} u_\infty \,.
\end{align}
Since we have shown that $u_\infty \in Y$ and $u_\infty \ne 0$, we conclude that $\hat{\Lambda} u_\infty \ne 0$.
Thus, $-\mu_\infty$ is an eigenvalue of $\hat{\Lambda}$ in $X$.
Then, by the assumption of the theorem $\mu_\infty$ must be $0$, which implies $\hat{\Lambda} u_\infty\in {\rm Ker}\, \hat{\Lambda}$. Thus, $\hat{\Lambda} u_\infty \in {\rm Ran}\, \hat{\Lambda} \cap {\rm Ker}\, \hat{\Lambda}$, and by the assumption we conclude that $\hat{\Lambda} u_\infty=0$, i.e., $u_\infty \in {\rm Ker}\, \hat{\Lambda}$. On the other hand, we have also seen $u_\infty \in Y=\big ({\rm Ker}\, \hat{\Lambda} \big )^\bot$, and hence, $u_\infty=0$. This is a contradiction, and \eqref{proof.thm.abstract.1.5} must hold. 

\

\noindent Step 4: {\it $\sigma (L_\alpha) \subset \{\zeta\in \C~|~\Re \zeta\leq -1\}$ and \eqref{thm.abstract.1.3} holds.}

\noindent 
Let $\zeta \in \C$ satisfy $\Re \zeta >-1$ and let $f\in X$. Let $v\in D_Y(\mathbb{Q}L_\alpha)$ be the unique solution to 
\begin{align}\label{proof.thm.abstract.1'.1}
(\zeta - \mathbb{Q}L_\alpha) v = \mathbb{Q} f\,,
\end{align}
which is well-defined for all sufficiently large $|\alpha|$ by \eqref{thm.abstract.1.1}. Let $w\in D_X(A)$ be the solution to 
\begin{align*}
(\zeta -A )w = - \mathbb{P} (\zeta- L_\alpha) v + \mathbb{P}f\,, \qquad \mathbb{P}=I-\mathbb{Q}\,,
\end{align*}
where the term $- \mathbb{P} (\zeta- L_\alpha) v$ coincides with $-i\alpha \mathbb{P} \hat{\Lambda} v$ by the invariance  $\mathbb{Q}A \subset A \mathbb{Q}$ and $v\in Y$.
That is, we have the following  formula for $w$:
\begin{align}\label{proof.thm.abstract.1'.3}
w = -i \alpha (\zeta-A)^{-1}\mathbb{P} \hat{\Lambda} (\zeta- \mathbb{Q}L_\alpha)^{-1} \mathbb{Q} f + (\zeta-A)^{-1} \mathbb{P} f\,.
\end{align}
Note that \eqref{est.assum.A} implies that $\zeta\in \C$ satisfying $\Re \zeta>-1$ belongs to the resolvent set of $A$,
and thus the above formula is well-defined.
Moreover,  from $\mathbb{Q}A\subset A\mathbb{Q}$ and  the above equation, we get: 
\begin{align*}
0=\mathbb{Q} (\zeta - A ) w = (\zeta - A )\mathbb{Q} w\,.
\end{align*}
Hence we have $\mathbb{Q}w=0$, that is, $w\in {\rm Ker}\, \hat{\Lambda}$.
Then $u=v+w$ solves 
\begin{align*}
(\zeta - L_\alpha) u = (\zeta - L_\alpha ) v + (\zeta - A) w & = \mathbb{Q}f + \mathbb{P} (\zeta - L_\alpha) v + (\zeta - A) w \\
& = \mathbb{Q} f + \mathbb{P} f=f\,.
\end{align*}
Hence $u\in D_X(L_\alpha)$ solves the resolvent problem, and the above construction also implies the uniqueness. Moreover, we have from the construction that 
\begin{align}\label{proof.thm.abstract.1'.2}
\mathbb{Q} (\zeta - L_\alpha)^{-1} f = (\zeta - \mathbb{Q}L_\alpha)^{-1} \mathbb{Q} f\,, \qquad f\in X\,.
\end{align}
Hence, \eqref{thm.abstract.1.3} holds by \eqref{thm.abstract.1.2}. 
The proof is complete.

\

\begin{rem}\label{rem.thm.abstract.1'}{\rm From \eqref{proof.thm.abstract.1'.3} and \eqref{proof.thm.abstract.1'.2} we have the formula
\begin{align}\label{eq.rem.thm.abstract.1'}
\begin{split}
\mathbb{Q} (\zeta - L_\alpha)^{-1}  & = (\zeta - \mathbb{Q}L_\alpha)^{-1} \mathbb{Q} \,,\\
\mathbb{P} (\zeta - L_\alpha )^{-1} & = -i \alpha (\zeta-A)^{-1}\mathbb{P} \hat{\Lambda} (\zeta- \mathbb{Q}L_\alpha)^{-1} \mathbb{Q} + (\zeta-A)^{-1} \mathbb{P} \,.
\end{split}
\end{align}
}
\end{rem}

\

Theorem \ref{thm.abstract.1} and its proof do not give any information on the rate of convergence for $|\alpha|\rightarrow \infty$. To obtain a rate we make further assumptions as follows.

\begin{assum}\label{assum.Lambda.rate} There exist $C>0$, $\tau\in (0,\infty]$, $m_0\geq 1$, and bounded nonnegative functions $h_j: [m_0,\infty) \times \R\rightarrow [0,\infty)$, $j=1,2$, satisfying $\displaystyle \lim_{m\rightarrow \infty} \sup_{\mu\in \R} h_j (m,\mu)=0$, such that the following statements hold.

\noindent {\rm (i)} ${\rm Ker}~\hat{\Lambda}\subset D_X(A)$.

\noindent  {\rm (ii)} {\rm (a)} It follows that
\begin{align}
\| B_2 u \|_X \leq C \| u\|_X\,, \qquad u \in X\,.\label{est.assum.Lambda.rate.0'}
\end{align}

\noindent {\rm (b)} For all $\mu\in \R$ and $m\geq m_0$ it follows that
\begin{align}
\| u\|_X^2 &\leq C \bigg (  m^2 \|  (\mu - \hat{\Lambda} )  u \|_X^2 +  h_1^2 (m,\mu) \| (-A)^\frac12 u \|_X^2 \bigg )  ~~ \qquad \text{ if }~~|\mu|\geq \tau ~~ \text{and} ~~ u\in D_X(A)\,,\label{est.assum.Lambda.rate'}\\
\begin{split}
\| u \|_X^2 &\leq C\bigg (  m^2 \| \mathbb{Q} (\mu - \hat{\Lambda} )  u \|_X^2 +  h_1^2 (m,\mu)  \| (-A)^\frac12 u \|_X^2  \bigg )  \quad  \text{ if} ~~ |\mu| < \tau ~~ \text{and} ~~ u\in D_X(A) \cap Y\,.\label{est.assum.Lambda.rate}
\end{split}
\end{align}

\noindent {\rm (iii)} There exists a densely defined closed operator $B_3: D_X(B_3) \rightarrow X$ such that 

\noindent {\rm (a)} $D_X((-A)^\frac12) \subset D_X (B_3)$ and 
\begin{align}
| \Im \langle A u, \hat{\Lambda} u \rangle_X | \leq C \| (-A)^\frac12 u \|_X \| B_3 u \|_X  \,,\qquad u \in D_X (A)\,,\label{est.assum.Lambda.rate.0''}
\end{align}
and 

\noindent {\rm (b)} for all $\mu\in \R$ and $m\geq m_0$ it follows that
\begin{align}
\| B_3 u \|_X^2 &\leq C \bigg (  m^2 \| (\mu - \hat{\Lambda} )  u \|_X^2 +  h_2^2 (m, \mu) \| (-A)^\frac12 u \|_X^2 \bigg ) ~~ \qquad \text{ if }~~|\mu|\geq \tau ~~ \text{and} ~~ u\in D_X(A)\,,\label{est.assum.Lambda.rate''}\\
\begin{split}
\| B_3 u \|_X^2 &\leq C\bigg (  m^2 \| \mathbb{Q} (\mu - \hat{\Lambda} )  u \|_X^2 +  h_2^2 (m,\mu)  \| (-A)^\frac12 u \|_X^2  \bigg )  \quad \text{ if} ~~ |\mu| < \tau ~~ \text{and} ~~ u\in D_X(A) \cap Y\,.\label{est.assum.Lambda.rate'''}
\end{split}
\end{align} 
\end{assum}

\begin{rem}\label{rem.assum.Lambda.rate.1}{\rm 
(1) The condition (ii) (b)  states quantitatively the absence of  eigenvalues of $\hat{\Lambda}$ in $\R\setminus \{0\}$. To be precise let us compare the conditions \eqref{est.assum.Lambda.rate'}-\eqref{est.assum.Lambda.rate} with Assumption \ref{assum.Lambda'}. Assume that ${\rm Ker}\, \hat{\Lambda} \subset D_X(A)$ holds. Then \eqref{est.assum.Lambda.rate'}-\eqref{est.assum.Lambda.rate} imply that $\hat{\Lambda}$ does not possess eigenvalues in $\R\setminus\{0\}$, otherwise one can take the limit $m\rightarrow \infty$ in \eqref{est.assum.Lambda.rate'}-\eqref{est.assum.Lambda.rate} for the eigenvalue $\mu$  and the eigenfunction $u$, which leads to a contradiction.  
Moreover, if in addition ${\rm Ker}\, \hat{\Lambda} ^2 \subset D_X(A)$, then for any $u\in {\rm Ker}\, \hat{\Lambda} \cap {\rm Ran} \hat{\Lambda}$ there exists $f\in D_X(A)\cap Y$ such that $u=\hat{\Lambda} f$. Then \eqref{est.assum.Lambda.rate} with $\mu=0$ implies $f=0$, which leads to $u=0$.  
Thus, in this case we also have ${\rm Ker}\, \hat{\Lambda} \cap {\rm Ran}\, \hat{\Lambda}=\{0\}$.
As a conclusion, under Assumption \ref{assum.Lambda.rate} (ii) and the condition ${\rm Ker}\, \hat{\Lambda}^2\subset D_X(A)$,  Assumption \ref{assum.Lambda'} is automatically satisfied. 

\noindent (2) The case $\tau=\infty$ means that the conditions \eqref{est.assum.Lambda.rate} and \eqref{est.assum.Lambda.rate'''} hold for all $\mu\in \R$. But since $\mathbb{Q}$ is nonlocal in actual applications,  in the case $\mu$ is away from $0$ the conditions \eqref{est.assum.Lambda.rate'} and \eqref{est.assum.Lambda.rate''} will be easier to check. If $\hat{\Lambda}$ is closed symmetric then $\mathbb{Q}$ in \eqref{est.assum.Lambda.rate} and \eqref{est.assum.Lambda.rate'''} is automatically dropped since $\mathbb{P}\hat{\Lambda}=0$ in this case. 
}
\end{rem}

\begin{rem}\label{rem.assum.Lambda.rate.2}{\rm  In fact, one can obtain the pseudospectral bound with some rate without assuming (iii) (b) of Assumption \ref{assum.Lambda.rate}. However,  the condition (iii) (b) is useful in obtaining the pseudospectral bound with a better rate in $\alpha$. Indeed, when $A$ is of the form $-A=T^* T$ for some densely defined closed operator $T$ as discussed in the work of \cite{Vi}, a natural candidate of $B_3$ is $B_1T(B_2-I) + [T, B_1] B_2$; formally we can compute as,  in virtue of the symmetry of $B_1$, 
\begin{align}
|\Im \langle Au, \hat{\Lambda} u \rangle_X| = |\Im \langle Tu, T B_1 B_2 u \rangle_X|  & = |\Im \langle Tu, B_1 T B_2 u + [T, B_1]B_2 u \rangle_X | \nonumber \\
& = |\Im \langle Tu, B_1 T (B_2-I) u + [T, B_1]B_2 u \rangle_X|\,.
\end{align}
Thus, the estimate \eqref{est.assum.Lambda.rate.0''} is valid with $B_3 = B_1T(B_2-I)+ [T,B_1] B_2$.
 When $B_2$ is a {\it smooth enough} perturbation from the identity operator, then one can even expect to take $B_3$ just as $[T, B_1]$. 
}
\end{rem}

\begin{rem}\label{rem.assum.Lambda.rate.3}{\rm The key idea of Assumption \ref{assum.Lambda.rate} is that we reduce the whole analysis to several coercive estimates of $\hat{\Lambda}$.
This idea is useful in actual applications. Indeed, when $A$ and $\hat{\Lambda}$ are (pseudo)differential operators,  the {\it order} of $\hat{\Lambda}$ is lower than $A-i\alpha \hat{\Lambda}$, and hence, the analysis of  $\hat{\Lambda}$ itself is expected to be simpler than the combined operator $A-i \alpha \hat{\Lambda}$.  In principle, the operator $A$ plays a role of recovering the regularity which was lost in the coercive estimates for $\hat{\Lambda}$. The functions $h_j$ in the assumption describe the degeneracy of the operator $\mu - \hat{\Lambda}$, which leads to an essential influence to the resolvent estimate for the full operator $L_\alpha$.
}
\end{rem}

The next theorem provides the information on the convergence rate, once we know the behavior of the functions $h_j$.
\begin{thm}\label{thm.abstract.2}  Suppose that Assumptions \ref{assum.A}, \ref{assum.Lambda},  and \ref{assum.Lambda.rate} hold. 
Then for any $\alpha \in \R$ the set $\{\zeta\in \C~|~\Re \zeta \geq 0 \}$ is contained in the resolvent set of $L_\alpha$ in $X$ and of $\mathbb{Q}L_\alpha$ in $Y$. Moreover, there exists a large number $M_0>0$ such that if $|\alpha|\geq M_0$ then the following resolvent estimate holds for any $\lambda\in \R$. 
\begin{align}
\| \mathbb{Q} (i\lambda + L_\alpha)^{-1}\|_{X\rightarrow X}  =\|  (i\lambda + \mathbb{Q}L_\alpha)^{-1}\|_{Y\rightarrow Y} \leq C F(\alpha,\frac{\lambda}{\alpha}) \,, \label{est.thm.abstract.2.1}
\end{align}
where
\begin{align}\label{def.F.main}
F (\alpha,\mu) =  \inf_{m_1,m_2\geq m_0} \bigg (  \frac{m_1}{|\alpha|}   + \frac{m_1^2 m_2^2}{\alpha^2}  + \frac{m_1^2 h_2(m_2,\mu)}{|\alpha|} + h_1^2(m_1,\mu) \bigg )\,
\end{align}
and $m_0$ is given in Assumption 4.
Here $C$ is independent of $\alpha$ and $\lambda$.
\end{thm}

\begin{rem}{\rm To evaluate $F(\alpha,\mu)$ we first choose $m_2$ so that $m_2^{2} = |\alpha| h_2 (m_2,\mu)$ holds, which gives the balance $\frac{m_1^2m_2^2}{\alpha^2}=\frac{m_1^2 h_2(m_2,\mu)}{|\alpha|}$ for any $m_1$. With this choice of $m_2$, the number $m_1$ is chosen so that $\frac{m_1}{|\alpha|} + \frac{m_1^2m_2^2}{\alpha^2}  + h_1^2(m_1,\mu)$ is minimized.
}
\end{rem}

\noindent {\it Proof of Theorem \ref{thm.abstract.2}.} 
It suffices to consider the case $\alpha>0$. Set $\mu=\frac{\lambda}{\alpha}$ and let $\tau\in (0,\infty]$ be the number in Assumption \ref{assum.Lambda.rate}.
For simplicity of notations we write $h_j$ instead of $h_j (m_j,\mu)$, $j=1,2$.
In the following argument any constant which is independent of $\alpha$, $\lambda$, $m_1$, and $m_2$ will be denoted by $C$, and thus, the constant $C$ can change from line to line.

\noindent (i) The case $|\mu|<\tau$. From the definition, we have for $u\in D_X(A) \cap Y$,
\begin{align*}
\mathbb{Q} (L_\alpha + i\lambda ) = \mathbb{Q} A - i \alpha \mathbb{Q} (\hat{\Lambda} - \mu )
\end{align*}
and 
\begin{align}\label{proof.thm.reduction.1.0} 
\langle  (\mathbb{Q} L_\alpha + i\lambda ) u, \big ( \hat{\Lambda} - \mu \big )  u \rangle_X & = \langle \mathbb{Q} A u, \big ( \hat{\Lambda}  - \mu \big ) u \rangle_X - i \alpha \| \mathbb{Q} (\hat{\Lambda} - \mu )  u \|_X^2\,.
\end{align}
By taking the imaginary parts of both sides, we obtain 
\begin{align*}
\alpha \| \mathbb{Q} \big ( \hat{\Lambda}  - \mu \big ) u \|_X^2  & \leq |\Im \langle  (\mathbb{Q} L_\alpha + i\lambda ) u, \mathbb{Q} \big ( \hat{\Lambda}  - \mu \big ) u \rangle_X| + | \Im \langle  -\mathbb{Q}A u,  \big ( \hat{\Lambda}  - \mu \big ) u \rangle_X | \\
& \leq \| (\mathbb{Q} L_\alpha + i\lambda )u\|_X \| \mathbb{Q} ( \hat{\Lambda}  - \mu \big ) u \|_X +  |\Im \langle -A u, \hat{\Lambda} u\rangle_X|
\end{align*}
Here we have used that for $u\in D_X(A) \cap Y$, we have  $\mathbb{Q} Au = A\mathbb{Q} u = Au$ and the self-adjointness of $A$, which gives $\Im \langle -\mathbb{Q}Au, \mu u \rangle_X=0$. For the mixing term $\Im \langle -Au, \hat{\Lambda}u\rangle_X$, we use Assumption \ref{assum.Lambda.rate} (iii) (a) and arrive at the estimate
\begin{align}\label{proof.thm.reduction.1.2} 
\alpha \| \mathbb{Q} \big ( \mu - \hat{\Lambda}\big ) u \|_X^2  
& \leq \frac{C}{\alpha} \| (\mathbb{Q} L_\alpha + i\lambda )u\|_X^2   +  C \| (-A)^\frac12 u \|_X  \| B_3 u \|_X\,.
\end{align}
Let us estimate $\| B_3 u \|_X$. Fix $m_2\geq m_0$.
We have from \eqref{est.assum.Lambda.rate'''} and \eqref{proof.thm.reduction.1.2},
\begin{align}\label{proof.thm.reduction.1.3} 
\|  B_3 u \|_X^2 & \leq C m_2^2 \| \mathbb{Q} ( \mu - \hat{\Lambda}) u \|_X^2 + Ch_2^2  \| (-A)^\frac12 u \|_X^2 \nonumber \\
\begin{split}
& \leq \frac{Cm_2^2}{\alpha} \bigg ( \frac{1}{\alpha} \| (\mathbb{Q} L_\alpha + i\lambda )u\|_X^2   +  \| (-A)^\frac12 u \|_X  \| B_3 u \|_X \bigg ) + h_2^2 \| (-A)^\frac12 u \|_X^2
\end{split}\nonumber \\
\begin{split}
& \leq \frac{Cm_2^2}{\alpha^2} \| (\mathbb{Q} L_\alpha + i\lambda )u\|_X^2  + C \bigg ( \frac{m_2^4}{\alpha^2} + h_2^2 \bigg ) \| (-A)^\frac12 u \|_X^2\,.
\end{split}
\end{align}
Then \eqref{proof.thm.reduction.1.2} and \eqref{proof.thm.reduction.1.3} imply 
\begin{align}\label{proof.thm.reduction.1.4}
\| \mathbb{Q} \big ( \mu - \hat{\Lambda}\big ) u \|_X^2 & \leq \frac{C}{\alpha^2} \|(\mathbb{Q} (L_\alpha + i\lambda) u \|_X^2  \nonumber \\
&   + \frac{C}{\alpha} \| (-A)^\frac12 u \|_X  \Big (  \frac{m_2^2}{\alpha^2} \| (\mathbb{Q} L_\alpha + i\lambda )u\|_X^2   + \big ( \frac{m_2^4}{\alpha^2} + h_2^2  \big ) \| (-A)^\frac12 u \|_X^2 \Big )^\frac12 \nonumber \\
& \leq  \frac{C}{\alpha^2} \|(\mathbb{Q} (L_\alpha + i\lambda) u \|_X^2 + C \big (\frac{m_2^2}{\alpha^2} + \frac{h_2}{\alpha} \big )   \| (-A)^\frac12 u \|_X^2\,.
\end{align}
Then, by combining \eqref{proof.thm.reduction.1.4} with the assumption \eqref{est.assum.Lambda.rate}, we see 
\begin{align}\label{proof.thm.reduction.1.5}
\| u \|_X^2 & \leq C m_1^2 \| \mathbb{Q} (\mu - \hat{\Lambda}) u \|_X^2 + Ch_1^2  \| (-A)^\frac12 u \|_X^2 \nonumber \\
& \leq \frac{C m_1^2}{\alpha^2} \| (\mathbb{Q} L_\alpha + i\lambda )u \|_X^2 +  C \big ( \frac{m_1^2m_2^2}{\alpha^2} + \frac{m_1^2h_2}{\alpha} + h_1^2 \big ) \| (-A)^\frac12 u \|_X^2\,.
\end{align}
Next we estimate $\| (-A)^\frac12 u \|_X^2$. We observe that the following identity holds:
\begin{align*}
\Re \langle (\mathbb{Q} L_\alpha + i\lambda )u, B_2 u \rangle_X = \Re \langle (L_\alpha + i\lambda )u, \mathbb{Q} B_2 u \rangle_X = \Re \langle Au\,, B_2u\rangle_X\,.
\end{align*}
Here we have used  $\mathbb{Q} B_2 = B_2$, $\Im \langle \hat{\Lambda} u\,, B_2 u\rangle_X=\Im \langle B_1 B_2 u\,, B_2 u\rangle_X=0$, and $\Im \langle u\,, B_2 u\rangle_X=0$ since $B_1$ and $B_2$ are symmetric.
By \eqref{est.assum.Lambda.2} we have $\Re \langle - A u\,, B_2 u\rangle_X \geq c_1 \langle -Au\,, u\rangle_X$, which gives  
\begin{align}\label{proof.thm.reduction.1.6} 
\| (-A)^\frac12 u \|_X^2  & \leq C \|( \mathbb{Q} L_\alpha + i\lambda )u\|_X \| B_2 u\|_X  \,, \qquad u\in D_X(A) \cap Y\,.
\end{align}
Thus we conclude from \eqref{est.assum.Lambda.rate.0'}, \eqref{proof.thm.reduction.1.5}, and \eqref{proof.thm.reduction.1.6} that
\begin{align}\label{proof.thm.reduction.1.7} 
\| u\|_X^2 & \leq C \bigg ( \frac{m_1^2}{\alpha^2}   + \frac{m_1^4m_2^4}{\alpha^4} + \frac{m_1^4 h_2^2}{\alpha^2}  + h_1^4    \bigg ) \,  \| (\mathbb{Q} L_\alpha + i\lambda )u\|_X^2 \,, \qquad u\in D_X(A) \cap Y \,.
\end{align}
Note that \eqref{proof.thm.reduction.1.7} is valid for any $m_1,m_2\geq m_0$ and $\alpha>0$.

\noindent (ii) The case $|\mu|\geq \tau$. In this case we drop the projection $\mathbb{Q}$ in \eqref{proof.thm.reduction.1.0} and use the identity 
\begin{align*}%\label{proof.thm.reduction.1.0'} 
\langle  ( L_\alpha + i\lambda ) u, \big ( \hat{\Lambda} - \mu \big )  u \rangle_X & = \langle A u, \big ( \hat{\Lambda}  - \mu \big ) u \rangle_X - i \alpha \| (\hat{\Lambda} - \mu )  u \|_X^2\,, \qquad u \in D_X(A)\,,
%& =\langle A u, \big ( \hat{\Lambda}  + \mu \big ) \mathbb{Q} u \rangle_X + i \alpha \langle \big ( \hat{\Lambda} + \mu \big ) u, (\hat{\Lambda} + \mu ) \mathbb{Q} u \rangle_X^2 \,.
\end{align*}
which gives by taking the imaginary parts, 
\begin{align*}
\alpha \| \big ( \hat{\Lambda}  - \mu \big ) u \|_X^2 
& \leq \| (L_\alpha + i\lambda )u\|_X \|  ( \hat{\Lambda}  - \mu \big ) u \|_X +  |\Im \langle -A u, \hat{\Lambda} u\rangle_X| \,.
%+ |\mu \Im \langle -A u, u\rangle_X | \\
%& \leq \| (L_\alpha - i\lambda )u\|_X \|  ( \hat{\Lambda}  + \mu \big ) u \|_X +  C \| B_2 u \|_X \| (-A)^{-\frac12} u \|_X + C \| (-A)^\frac12 u \|_X \| [T, B_1] B_2 u \|_X\,.
\end{align*}
Thus, as in \eqref{proof.thm.reduction.1.2}, we obtain
\begin{align}\label{proof.thm.reduction.1.8} 
\begin{split}
\alpha \| \big ( \hat{\Lambda}  - \mu \big ) u \|_X^2 & \leq \frac{C}{\alpha} \| ( L_\alpha + i\lambda )u\|_X^2  +  \| (-A)^\frac12  u \|_X \| B_3 u \|_X\,, \qquad u\in D_X(A) \,.
\end{split}
\end{align}
Then the estimates of the terms $B_3 u$ and $(-A)^{-\frac12}u$ are obtained in the same manner as in the case (i), and we have 
\begin{align}\label{proof.thm.reduction.1.9} 
\begin{split}
\| B_3 u \|_X^2 & \leq \frac{Cm_2^2}{\alpha^2} \| (L_\alpha + i\lambda )u\|_X^2   + C \bigg ( \frac{m_2^4}{\alpha^2} + h_2^2 \bigg ) \| (-A)^\frac12 u \|_X^2\,, \\
\| \big ( \mu - \hat{\Lambda}\big ) u \|_X^2 
& \leq  \frac{C}{\alpha^2} \| (L_\alpha + i\lambda) u \|_X^2 + C \big (\frac{m_2^2}{\alpha^2} + \frac{h_2}{\alpha} \big )   \| (-A)^\frac12 u \|_X^2\,.
\end{split}
\end{align}
Then $u$ is estimated as follows, by arguing as in \eqref{proof.thm.reduction.1.5}:
\begin{align}\label{proof.thm.reduction.1.11}
\begin{split} 
\| u \|_X^2 
& \leq \frac{C m_1^2}{\alpha^2} \| (L_\alpha + i\lambda )u \|_X^2 +  C_1 \big ( \frac{m_1^2m_2^2}{\alpha^2} + \frac{m_1^2h_2}{\alpha} + h_1^2 \big ) \| (-A)^\frac12 u \|_X^2\,. 
\end{split}
\end{align}
Next we observe that $A\mathbb{P}$, where $\mathbb{P}=I-\mathbb{Q}$, is a bounded operator by the assumption ${\rm Ker}\, \hat{\Lambda}\subset D_X(A)$, which gives  from $\mathbb{Q}A\subset A\mathbb{Q}$,
$$ \| (-A)^\frac 12 u \|_X^2 = \langle -Au, u\rangle_X =  \langle -A \mathbb{P}u, \mathbb{P} u \rangle_X +  \langle -A \mathbb{Q}u, \mathbb{Q} u\rangle_X\leq C_2 \| \mathbb{P} u \|_X^2 + \langle -A \mathbb{Q}u, \mathbb{Q} u \rangle_X\,.$$
Hence, for sufficiently large $\alpha$ and $m_j$ so that 
\begin{align}\label{proof.thm.reduction.1.11'}
C_1 C_2 \bigg (  \frac{m_1^2m_2^2}{\alpha^2} + \frac{m_1^2 h_2}{\alpha}  + h_1^2   \bigg ) \leq \frac14\,,
\end{align}
the estimate \eqref{proof.thm.reduction.1.11} yields
\begin{align}\label{proof.thm.reduction.1.12} 
\begin{split}
\| u\|_X^2 & \leq  \frac{Cm_1^2}{\alpha^2}  \| (L_\alpha + i\lambda )u\|_X^2 + C \bigg (  \frac{m_1^2m_2^2}{\alpha^2} + \frac{m_1^2 h_2}{\alpha} + h_1^2 \bigg )  \| (-A)^\frac12 \mathbb{Q} u \|_X^2\,.
\end{split}
\end{align}
Since $\langle -A \mathbb{Q} u, \mathbb{Q} u\rangle_X$ is estimated as \eqref{proof.thm.reduction.1.6} but with $u$ replaced by $\mathbb{Q}u$, we obtain the estimate like \eqref{proof.thm.reduction.1.7} such as 
\begin{align}\label{proof.thm.reduction.1.13} 
\| u\|_X^2 & \leq C\frac{m_1^2}{\alpha^2}  \| (L_\alpha + i\lambda )u\|_X^2  + C \bigg (  \frac{m_1^4m_2^4}{\alpha^4} + \frac{m_1^4 h_2^2}{\alpha^2}  + h_1^2    \bigg ) \,  \| \mathbb{Q} (L_\alpha + i\lambda )u\|_X^2 \nonumber \\
\begin{split}
& \leq C  \bigg (  \frac{m_1^2}{\alpha^2} + \frac{m_1^4m_2^4}{\alpha^4} + \frac{m_1^4 h_2^2}{\alpha^2}  + h_1^4   \bigg ) \,  \| (L_\alpha + i\lambda )u\|_X^2  \,, \qquad u\in D_X(A) \,.
\end{split}
\end{align}
This is the desired estimate in the case $|\mu|\geq \tau$. 
Note that \eqref{proof.thm.reduction.1.13} is valid for any $m_1,m_2\geq m_0$ and $\alpha>0$  satisfying  \eqref{proof.thm.reduction.1.11'}. Such a set is not empty  since $\displaystyle \lim_{m_j\rightarrow \infty} \sup_{\mu\in \R} h_j (m_j,\mu)=0$, and in particular, we can find a positive constant $M_0$ uniformly in $\lambda$ so that if $\alpha\geq M_0$ then there exists $(m_1,m_2)$ satisfying \eqref{proof.thm.reduction.1.11'}.

Now we recall that $\sigma (L_\alpha)$, $\sigma_Y (\mathbb{Q}L_\alpha)$, and $\sigma_{\mathbb{P}X} (A|_{\mathbb{P}X})$ consist  only of  discrete eigenvalues.
Moreover, we have from $\mathbb{Q} A\subset A\mathbb{Q}$, 
\begin{align}\label{proof.thm.reduction.1.14} 
\sigma (L_\alpha ) = \sigma_Y (\mathbb{Q}L_\alpha) \cup \sigma_{\mathbb{P}X} (A|_{\mathbb{P}X})\,,
\end{align}
and the formula \eqref{eq.rem.thm.abstract.1'} holds for all $\zeta \in \C\setminus \big (\sigma (L_\alpha)\big )$. 
Indeed, let $\zeta\in \C$ be  in the  resolvent of $\mathbb{Q}L_\alpha$ in $Y$ and also of $A|_{\mathbb{P}X}$in $\mathbb{P}X$. Let us show that $\zeta-L_\alpha$ is injective in $X$, which shows that $\zeta$ belongs to the 
  resolvent of $L_\alpha$ since $\sigma (L_\alpha)$ in $X$ consists of eigenvalues. 
If $u\in D_X(A)$ satisfies $(\zeta - L_\alpha )u=0$, then $(\zeta - \mathbb{Q}L_\alpha)\mathbb{Q} u=0$ by the invariance $\mathbb{Q}A\subset A\mathbb{Q}$ and $\hat{\Lambda} u = \hat{\Lambda}\mathbb{Q}u$.
Then $\mathbb{Q}u=0$ by the assumption, and therefore, $(\zeta-A) \mathbb{P}u=0$, which also gives $\mathbb{P}u=0$ by the assumption. Thus $u=0$, and we have shown the inclusion
\begin{align*}
 \sigma (L_\alpha ) \subset \sigma_Y (\mathbb{Q}L_\alpha) \cup \sigma_{\mathbb{P}X} (A|_{\mathbb{P}X}) \,.
\end{align*}
On the other hand, let $\zeta$ belong to the resolvent set of $L_\alpha$ in $X$.
Since $\sigma_{\mathbb{P}X} (A|_{\mathbb{P}X})\subset \sigma (L_\alpha)$ holds by the assumptions $\mathbb{Q} A \subset A \mathbb{Q}$ and ${\rm Ker}\, \hat{\Lambda}\subset D_X(A)$, $\zeta$ is also a resolvent of $A|_{\mathbb{P}X}$ in $\mathbb{P}X$.
If $u\in D_X(A)\cap Y$ satisfies $(\zeta -  \mathbb{Q}L_\alpha) u=0$ then  by setting $v\in D_X(A) \cap \mathbb{P}X$ as $v= - i \alpha (\zeta-A|_{\mathbb{P}X})^{-1} \mathbb{P} \hat{\Lambda}u$ we see that $w=u+v$ solves from $\mathbb{Q}A\subset A\mathbb{Q}$,
\begin{align*}
(\zeta - L_\alpha ) w = (\zeta - L_\alpha ) u + (\zeta - L_\alpha ) v  = (\zeta- \mathbb{Q}L_\alpha ) u  + i\alpha \mathbb{P}\hat{\Lambda} u -i \alpha \mathbb{P}\hat{\Lambda} u =0\,.
\end{align*}
Since $\zeta$ is a resolvent of $L_\alpha$ in $X$, $w=0$. This implies $u=\mathbb{Q}w=0$.
Hence $\zeta$ is a resolvent of $\mathbb{Q}L_\alpha$ in $Y$. Thus \eqref{proof.thm.reduction.1.14} holds.

By arguing as in Step 2 of the proof of Theorem \ref{thm.abstract.1}, we can show that
$\sigma_Y (\mathbb{Q}L_\alpha) \subset \{\zeta\in \C~|~\Re \zeta <0 \}$ for all $\alpha$. Therefore, we observe from $\sigma(A)\subset \{\zeta\in \C~|~\Re \zeta \leq -1\}$ and \eqref{proof.thm.reduction.1.14} that $\sigma (L_\alpha) \subset \{\zeta\in \C~|~\Re \zeta <0 \}$ for all $\alpha$.
In particular, $i\R$ belongs to the resolvent set of $L_\alpha$ and also of $\mathbb{Q}L_\alpha$.

Let $\lambda\in \R$.  If $|\frac{\lambda}{\alpha}|<\tau$ then \eqref{proof.thm.reduction.1.7} gives the estimate for the resolvent $(i\lambda + \mathbb{Q} L_\alpha)^{-1}$ in $Y$.
If $|\frac{\lambda}{\alpha}|\geq \tau$ then \eqref{proof.thm.reduction.1.13} yields the estimate of the resolvent $(i\lambda + L_\alpha)^{-1}$ in $X$ as long as \eqref{proof.thm.reduction.1.11'} is satisfied.
The estimate of the resolvent $(i\lambda + \mathbb{Q}L_\alpha)^{-1}$ in $Y$ is then obtained from the formula
\begin{align*}
(i\lambda + \mathbb{Q}L_\alpha)^{-1}f  = \mathbb{Q} (i\lambda + L_\alpha)^{-1} f \,, \qquad f\in Y\,,
\end{align*}
and by using the inequality $\| \mathbb{Q} u \|_X\leq \| u\|_X$. 
As a summary, there exists $M_0>0$ such that we have
\begin{align}\label{proof.thm.reduction.1.15} 
\begin{split}
\| (i\lambda + \mathbb{Q}L_\alpha)^{-1}\|_{Y\rightarrow Y}  & \leq C \bigg (  \frac{m_1}{|\alpha|} + \frac{m_1^2 m_2^2}{\alpha^2} + \frac{m_1^2 h_2}{|\alpha|}  + h_1^2   \bigg )  \,,
\end{split}
\end{align}
as long as $\alpha\geq M_0$. The proof is complete.

\begin{rem}{\rm From \eqref{proof.thm.reduction.1.6} in the proof of Theorem \ref{thm.abstract.2}, we observe that 
\begin{align}
\| (-A)^\frac12 (i\lambda + \mathbb{Q}L_\alpha)^{-1}  \|_{Y\rightarrow Y} \leq C F(\alpha,\frac{\lambda}{\alpha})^\frac12\,.
\end{align}
Moreover, when $m_1$ and $m_2$ are the numbers such that the infimum in the definition of $F$ is evaluated, we have 
\begin{align}
\| B_3 (i\lambda + \mathbb{Q}L_\alpha)^{-1} \|_{Y\rightarrow X} & \leq C \Big ( \frac{m_2}{|\alpha|} + \big ( \frac{m_2^2}{|\alpha|} + h_2 (m_2,\mu) \big ) F(\alpha,\frac{\lambda}{\alpha})^\frac12 \Big ) \,,\\
\| \mathbb{Q} (\frac{\lambda}{\alpha} - \hat \Lambda) (i\lambda + \mathbb{Q}L_\alpha)^{-1} \|_{Y\rightarrow Y} & \leq   C \Big ( \frac{1}{|\alpha|}  + \big ( \frac{m_2}{|\alpha|} + \frac{h_2(m_2,\mu)^\frac12}{|\alpha|^\frac12}\big ) F(\alpha,\frac{\lambda}{\alpha})^\frac12 \Big )\,, \quad |\frac{\lambda}{\alpha}|<\tau\,,\\
\| (\frac{\lambda}{\alpha} - \hat \Lambda) (i\lambda + L_\alpha)^{-1} \|_{X\rightarrow X} & \leq   C \Big ( \frac{1}{|\alpha|}  + \big ( \frac{m_2}{|\alpha|} + \frac{h_2(m_2,\mu)^\frac12}{|\alpha|^\frac12}\big ) F(\alpha,\frac{\lambda}{\alpha})^\frac12 \Big )\,, \quad |\frac{\lambda}{\alpha}|\geq \tau\,.
\end{align}
}
\end{rem}

\subsection{On the proof of Assumption 4 (ii) (b) and (iii) (b) in actual applications}

In actual applications to the Kolmogorov flow or the Lamb-Oseen vortex, the most trivial part is to verify the interpolation inequalities \eqref{est.assum.Lambda.rate'}, \eqref{est.assum.Lambda.rate}, and \eqref{est.assum.Lambda.rate''}, \eqref{est.assum.Lambda.rate'''}. To find appropriate $B_3$ itself is not a difficult task in these examples, by recalling Remark \ref{rem.assum.Lambda.rate.2}.
We will show these interpolation inequalities by a contradiction argument. The approach using a contradiction argument is standard, and one can go back to the very abstract and classical result as follows: Suppose that the triple of Banach spaces $(X,Y,Z)$ satisfies the embedding property $Z\hookrightarrow \hookrightarrow X$ and $Z\hookrightarrow Y$, in particular, $Z$ is compactly embedded in $X$ (note that, for example in \eqref{est.assum.Lambda.rate'}, we can consider $\|u\|_Y=\|(\mu-\hat{\Lambda})u\|_X$ and $\|u\|_Z=\|(-A)^\frac12 u\|_X$ when it is assumed that $\mu-\hat{\Lambda}$ is injective and that $(-A)^{-\frac12}$ is compact). Then for any $\epsilon>0$ there exists $C_\epsilon>0$ such that $\|u\|_X \leq C_\epsilon \| u \|_Y + \epsilon \|u\|_Z$ for any $u\in Z$. Indeed, one can easily prove this inequality, depending on $\epsilon$, by a contradiction argument. Since the assumption of this abstract result is too general, we do not know how $C_\epsilon$ depends on $\epsilon$. 
But for concrete applications, we expect to be able to estimate the degeneracy of $\mu-\hat{\Lambda}$ that gives an information about the concrete dependence of $C_\epsilon$ on $\epsilon$ when $\|u\|_Y=\|(\mu-\hat{\Lambda})u \|_X$ (this leads to a candidate of $h_1(m,\mu)$ in \eqref{est.assum.Lambda.rate'}). Then we can try to prove the interpolation-type inequality by a contradiction but together with the presence of $h_1(m,\mu)$.

In the application to the Kolmogorov flow and the Lamb-Oseen vortex, the operator $\hat{\Lambda}$ is essentially of the form $\hat{\Lambda}_1+\hat{\Lambda}_2$, where $\hat{\Lambda}_1$ is a simple local operator  (multiplication operator) and $\hat{\Lambda}_2$ is a nonlocal compact operator. Then we expect that the degeneracy of $\mu-\hat{\Lambda}$ is dominated by $\mu-\hat{\Lambda}_1$, and thus, the interpolation-type inequality such as \eqref{est.assum.Lambda.rate'} has a close relation to the similar inequality but with $\mu-\hat{\Lambda}$ replaced by $\mu-\hat{\Lambda}_1$. This is indeed shown to be true in the above two examples, though the whole proof requires a long argument.  
It will be useful to point out that the interpolation-type inequalities in Assumption \ref{assum.Lambda.rate} are also related to the estimate of the limiting absorption principle (LAP) typically stated as $\lim_{\epsilon\downarrow} \| (-A)^{-\frac12} (\mu\pm i\epsilon - \hat{\Lambda})^{-1} f \|_X \leq C \|(-A)^\frac12 f \|_X$, for which a contradiction argument is a familiar tool in the proof. Indeed, our argument in Lemma \ref{lem.Lambda.Kol.rate.2} for the Kolmogorov flow share a common feature with the proof of the limiting absorption principle around the shear flows obtained by Wei, Zhang, and Zhao \cite{WZZ2}. But on the other hand, there are some differences in technical details between the proof of \eqref{est.assum.Lambda.rate'} and the proof of LAP in \cite{WZZ2}, mainly due to the difference of the regularity condition on $f=(\mu-\hat{\Lambda})u$; in \eqref{est.assum.Lambda.rate'} we impose $f\in X$, while in LAP it is $f\in D_X((-A)^\frac12)$. 
In fact, the lower regularity condition $(\mu-\hat{\Lambda}) u\in X$ in \eqref{est.assum.Lambda.rate'} makes the argument more technical in the analysis around the critical points.

\section{Application to Kolmogorov flow}\label{sec.Kol}

In this section we study the spectral property of the operator \eqref{eq.rescaled.linear} related to the linearization of the Kolmogorov flow.  Set 
\begin{align*}
X=L^2_0 (\T^2) = \{ \omega \in L^2 (\T^2)~|~\int_0^{2\pi} \omega (x,y) \dd x=0~~a.e.~y \in \T\}\,.
\end{align*}
Let $A$ be the realization of $\Delta=\partial_x^2 + \partial_y^2$ in $L_0^2(\T^2)$, i.e., 
\begin{align*}
D(A)=W^{2,2}(\T^2)\cap L^2_0 (\T^2)\,, \qquad A \omega = \Delta \omega\,, \quad \omega \in D(A)\,.
\end{align*}
Next let us denote by $M_g$ the multiplication operator with the multiplier $g$, i.e., $M_g f = g f$.
We denote by $\hat{\Lambda}$ the realization of $-i \partial_x M_{\sin y}  (I + A^{-1})$ in $L^2_0(\T^2)$ which is given by 
\begin{align*}
D(\hat{\Lambda})=\{\omega \in L^2_0 (\T^2)~|~\partial_x M_{\sin y} \omega \in L^2_0 (\T^2)\}\,, \qquad \hat{\Lambda} \omega = -i \partial_x M_{\sin y}  (I + A^{-1}) \omega \,, \quad \omega\in D(\hat{\Lambda})\,.
\end{align*}
Since $i\partial_x$ is realized as a self-adjoint operator in $L^2_0 (\T^2)$ and $M_{\sin y} (I + A^{-1})$ is bounded in $L^2_0(\T^2)$, the operator $\hat{\Lambda}$ is a closed operator. Moreover, since $H^1 (\T^2) \cap L^2_0 (\T^2) \subset D(\hat{\Lambda})$, it is densely defined in $L^2_0(\T^2)$. We are interested in the spectral property of 
\begin{align}\label{def.L_alpha.Kol}
L_\alpha = A- i \alpha \hat{\Lambda}\,, \qquad D(L_\alpha) = D(A)\,.
\end{align} 
Let us denote by $Y$ the closed subspace of $L^2_0 (\T^2)$ defined by 
\begin{align}
Y=\{\omega \in L^2_0(\T^2)~|~\int_0^{2\pi} (\mathcal{P}_l \omega )(x,y) \dd y =0\,,~ |l|=1\,, ~\text{for all}~  x\in \T \}\,,
\end{align}
where 
\begin{align*}
(\mathcal{P}_l \omega )(x,y) = \frac{1}{2\pi} \int_0^{2\pi} \omega (s,y) e^{-i l s} \dd s \ e^{i l x}\,, \qquad l\in \Z\,.
\end{align*} 
The orthogonal projection from $L^2_0(\T^2)$ to $Y$ is denoted by $\mathbb{Q}$. 
We observe that $Y$ is an invariant space under the action of $A$ and that 
\begin{align*}
Y ~\text{ is the orthogonal complement of }~~\{a\cos x + b\sin x~|~a,b\in \C \} ~~\text{in}~~X\,. 
\end{align*}
The spaces $L^2_0 (\T^2)$ and $Y$ are diagonalized as  
\begin{align}
L^2_0(\T^2) = \oplus_{l\in \Z\setminus\{0\}} \mathcal{P}_l L^2_0 (\T^2)\,, \qquad 
Y = \oplus_{l\in \Z\setminus\{0\}} \mathcal{P}_l Y\,,
\end{align}
and each of $\mathcal{P}_l L^2_0 (\T^2)$ and $\mathcal{P}_l Y$ is identified with $L^2(\T)$ and $Y_l$ respectively, where 
\begin{align}
Y_l =
\begin{cases}
L^2(\T) \qquad  \qquad \text{ if } ~ l\ne \pm 1\,,\\
\{ f \in L^2(\T)~|~\int_0^{2\pi} f \dd y =0\} \qquad \text{ if } ~ l=\pm 1\,.
\end{cases}
\end{align}
The orthogonal projection from $L^2(\T)$ to $Y_l$ is denoted by $\mathbb{Q}_l$.
Since $\mathcal{P}_l L^2_0 (\T^2)$ and $\mathcal{P}_l Y$ are invariant spaces for $L_\alpha$, the operator $L_\alpha$ is also diagonalized as 
\begin{align}
L_\alpha = \oplus_{l\in \Z\setminus\{0\}} L_\alpha|_{\mathcal{P}_l L^2_0 (\T^2)}\,,
\end{align}
where $L_\alpha|_{\mathcal{P}_l L^2_0 (\T^2)}$ is the restriction of $L_\alpha$  to  the invariant subspace  $\mathcal{P}_l L^2_0 (\T^2)$, which is identified with $L_{\alpha,l}$ in $L^2(\T)$ defined as follows:
\begin{align}\label{def.Kol.2}
L_{\alpha,l} =  A_l - i \alpha l \hat{\Lambda}_{l}  \,,\qquad D(L_{\alpha,l}) = W^{2,2}(\T)\,,
\end{align}
where
\begin{align}\label{def.Kol.1}
\begin{split}
A_l & = \partial_y^2 -l^2\,, \qquad D (A_l) = W^{2,2}(\T)\,,\\
\hat{\Lambda}_{l}    & = M_{\sin y} \, (I + A_l^{-1} ) \,, \qquad D(\hat{\Lambda}_l) = L^2 (\T)\,.
\end{split}
\end{align}

It is straightforward to see the following result.
\begin{prop}\label{prop.Kol.assum} Let $|l|\geq 1$. {\bf (i)} $-A_l$ is positive self-adjoint in $L^2(\T)$ and satisfies Assumption \ref{assum.A} in $L^2(\T)$. Moreover, the invariance $\mathbb{Q}_l A_l \subset A_l \mathbb{Q}_l$ holds.

\noindent {\rm (ii)} $({\rm Ker}\, \hat{\Lambda}_l )^\bot = Y_l$, ${\rm Ran}\, \hat{\Lambda}_l \cap {\rm Ker}\, \hat{\Lambda}_l=\{0\}$, and $\hat{\Lambda}_l$ satisfies Assumption \ref{assum.Lambda} (i).
\end{prop}

\noindent {\it Proof.} We give a proof only for the statement ${\rm Ran}\, \hat{\Lambda}_l \cap {\rm Ker}\, \hat{\Lambda}_l=\{0\}$ with $l=\pm 1$, since the other statements are easy to check.
Let $f\in {\rm Ran}\, \hat{\Lambda}_l \cap {\rm Ker}\, \hat{\Lambda}_l$ with $l=\pm 1$. Then, since ${\rm Ker}\, \hat{\Lambda}_{\pm 1}=\{Const.\}$, there exists a constant $c$ and a function $g\in L^2(\T)$ such that $f=c=\hat{\Lambda}_l g$. By the definition of $\hat{\Lambda}_l$, we have $\sin y \, \big ( I + A_l^{-1}\big ) g =c$.
However, $\displaystyle \big ( I + A_l^{-1}\big ) g=\frac{c}{\sin y}$ cannot belong to $L^2(\T)$ if $c\ne 0$.
Hence, we must have $c=0$, that is, $f=0$. The proof is complete.

\

The following corollary immediately follows from the above proposition.
\begin{cor} {\rm (i)} $-A$ is positive self-adjoint in $L^2_0(\T^2)$ and satisfies Assumption \ref{assum.A} in $L^2_0(\T^2)$. Moreover, the invariance $\mathbb{Q}A\subset A \mathbb{Q}$ holds.

\noindent {\rm (ii)} ${\rm Ker}\, \hat{\Lambda} = \{ f=a \sin x + b \cos x\,, ~a,b\in \C \}$, $({\rm Ker}\, \hat{\Lambda})^\bot=Y$, ${\rm Ran}\, \hat{\Lambda} \cap {\rm Ker}\, \hat{\Lambda} =\{0\}$, and $\hat{\Lambda}$ satisfies Assumption \ref{assum.Lambda} (i).
\end{cor}

\subsection{Estimate without rate}\label{subsec.Kol.norate}
In this subsection we aim to apply Theorem \ref{thm.abstract.1}.
Let us first check Assumption \ref{assum.Lambda} (ii) for $L_{\alpha,l}$ in $L^2(\T)$.
We observe that
\begin{align}\label{def.Kol.B_j}
\hat{\Lambda}_l = M_{\sin y} B_{2,l}\,,  \qquad B_{2,l} = \big ( I + A_l^{-1} \big )\,,
\end{align}  
and $B_{2,l}$ is bounded self-adjoint in $L^2(\T)$. 
We can also check that ${\rm Ker}\, \hat{\Lambda}_l = {\rm Ker}\, B_{2,l}$ without difficulty.
The operator $B_{2,l}$ is positive in $Y_l$. To see this we set $\phi=A_l^{-1} f$ for $f\in Y_l$, which satisfies 
\begin{align}
\| \partial_y \phi \|_{L^2}^2 + l^2 \|\phi \|_{L^2}^2 = - \langle f \,, \phi\rangle_{L^2} \leq \| f\|_{L^2} \| \phi\|_{L^2}\,.
\end{align}
Since $f,\phi\in Y_l$ we see $\|\partial_y \phi\|_{L^2}^2\geq \| \phi\|_{L^2}^2$ if $l=\pm 1$, and $\| \partial_ y \phi\|_{L^2}^2\geq 0$ if $|l|\geq 2$. Thus we have 
\begin{align}
\| \phi \|_{L^2}\leq 
\begin{cases}
& \displaystyle 2^{-1} \| f \|_{L^2} \qquad \text{ if }~l=\pm 1\,,\\
& l^{-2} \| f \|_{L^2} \qquad \text{ if }~|l|\geq 2\,.
\end{cases}
\end{align}
Hence, for $|l|\geq 1$ and $f \in Y_l$,
\begin{align}\label{coercive.B_2}
\langle f\,, B_{2,l} f \rangle_{L^2} & = \|f\|_{L^2}^2 + \langle f\,, \phi \rangle_{L^2}  \geq  \| f \|_{L^2}^2 -\| f\|_{L^2} \| \phi\|_{L^2}  \geq \| f \|_{L^2}^2 - \frac12 \| f\|_{L^2}^2  = \frac12 \| f \|_{L^2}^2\,.
\end{align}
Since \eqref{est.assum.Lambda.2} is also not difficult to check, Assumptions \ref{assum.A} and \ref{assum.Lambda} are satisfied by the above operators.
To apply Theorem \ref{thm.abstract.1} it remains to show 

\begin{prop}\label{prop.kernel.Kol} $\sigma (\hat{\Lambda}_l)=[-1,1]$.  Moreover, $\hat{\Lambda}_l$ in $L^2(\T)$ does not have eigenvalues in $\C\setminus \{0\}$.
\end{prop}

\begin{rem}{\rm Proposition \ref{prop.kernel.Kol} and its Corollary \ref{cor.prop.kernel.Kol} below are not new, and have been proved in a more general framework; see Lin and Xu \cite[Lemmas 2.4, 5.1]{LX}. We give a proof here just for  the convenience of  the reader.
}
\end{rem}

\noindent {\it Proof of Proposition \ref{prop.kernel.Kol}.} We first observe that the spectrum of $M_{\sin y}$ in $L^2(\T)$ consists of the essential spectrum and is $[-1,1]$. Since $M_{\sin y} A_l^{-1}$ is compact in $L^2(\T)$, the essential spectrum of $\hat{\Lambda}_{l}$ coincides with the one of $M_{\sin y}$, and thus, is $[-1,1]$.
Hence it suffices to consider the existence of eigenvalues of $\hat{\Lambda}_l$.
Suppose that  $f \in L^2 (\T)$ and $\mu\in \C\setminus\{0\}$ satisfies 
\begin{align}\label{proof.prop.kernel.Kol.1}
\hat{\Lambda}_l f =\mu f\,.
\end{align}
We first consider the case $\mu\in \R\setminus \{0\}$.
Then we have $(\sin y - \mu) f + \sin y\,  A_l^{-1} f=0$, and thus, 
\begin{align}\label{proof.prop.kernel.Kol.1'}
f+\frac{\sin y}{\sin y-\mu} A_l^{-1} f = 0\,, \qquad y \notin S_\mu\,, 
\end{align}
where $S_\mu=\{\theta \in \T~|~\sin \theta =\mu\}$.
Note that $\phi=A_l^{-1} f\in W^{2,2}(\T)$ is a $C^{1+\delta}(\T)$ function for some $\delta>0$ by the Sobolev embedding inequality. From \eqref{proof.prop.kernel.Kol.1'} we see that $f\in C^{1+\delta}(\T\setminus\{S_\mu\})$,
and also  \eqref{proof.prop.kernel.Kol.1'} implies that $\phi (y_\mu)=0$ for $y_\mu\in S_\mu$, otherwise $f$ cannot be in $L^2 (\T)$.
By the bootstrap argument and \eqref{proof.prop.kernel.Kol.1'}, we see that $f$ is smooth in $\T\setminus\{S_\mu\}$. Thus $\phi$ is smooth in $\T\setminus\{S_\mu\}$ and solves the ODE
\begin{align*}
(M_{\sin y} - \mu ) A_l \phi + M_{\sin y} \phi =0\,.
\end{align*}
By the identity $M_{\sin y}=M_{\sin y}-\mu + \mu$, we have $(M_{\sin y}-\mu ) (A_l + 1)\phi + \mu \phi =0$, and thus,
\begin{align}\label{proof.prop.kernel.Kol.2}
-  (A_l +1 )\phi + \frac{\mu}{\mu - \sin y} \phi =0\,, \qquad y\in \T\setminus \{S_\mu \}\,.
\end{align}

\noindent Case (i) $\mu \geq 0$. 
When $0\leq \mu < 1$ let $y_\mu, z_\mu \in S_\mu$ be the points such that $y_\mu\in [\frac{1}{2}\pi,\pi]$ and $z_\mu \in [2\pi, \frac52\pi]$ (they are uniquely determined). When $\mu \geq 1$ we simply take $y_\mu=\frac12\pi$ and $z_\mu=\frac52\pi$.
Then $\mu-\sin y\geq 0$ for $y\in (y_\mu,z_\mu)$, and we obtain 
\begin{align}\label{proof.prop.kernel.Kol.3}
\int_{y_\mu}^{z_\mu}  (-A_l -1 ) \phi \,  \bar{\phi} \dd y + \int_{y_\mu}^{z_\mu} \frac{\mu}{\mu - \sin y} |\phi|^2 \dd y=0\,.
\end{align}
Note that the second integral converges due to the regularity $\phi\in C^{1+\delta} (\T)$ and $\phi(y_\mu)=\phi(z_\mu)=0$ when $0\leq \mu\leq 1$. 
As for the first integral, the integration by parts and the condition $\phi(y_\mu)=\phi(z_\mu)=0$ when $0\leq \mu\leq 1$ and the periodicity of $\phi$ when $\mu>1$ yield
\begin{align*}
\int_{y_\mu}^{z_\mu}  (-A_l -1 ) \phi \,  \bar{\phi} \dd y = \int_{y_\mu}^{z_\mu} |\partial_y \phi|^2 +(l^2-1) |\phi|^2 \dd y\,,
\end{align*}
thus we have 
\begin{align}\label{proof.prop.kernel.Kol.4}
\int_{y_\mu}^{z_\mu} |\partial_y \phi|^2 +(l^2-1) |\phi|^2 \dd y + \int_{y_\mu}^{z_\mu} \frac{\mu}{\mu - \sin y} |\phi|^2 \dd y=0\,.
\end{align}
Hence $\phi=0$ in $[y_\mu, z_\mu]$. When $\mu\geq 1$ we clearly have $\phi=0$ on $[\frac12\pi, \frac52\pi]$. When $0\leq \mu<1$ we see from $\phi\in C^{1+\delta}(\T)$ that $\phi (y_\mu) = \phi'(y_\mu)=0$.
However, since the singularity of $\displaystyle \frac{1}{\mu-\sin y}$ is first order when $0\leq\mu<1$ it is easy to see that any $C^{1+\delta}$ solution $\phi$ to the ODE \eqref{proof.prop.kernel.Kol.2} satisfying $\phi(y_\mu) = \phi'(y_\mu)=0$ must be trivial. Thus, we have $f=0$ in $\T$.

\noindent Case (ii) $\mu<0$. The argument is the same as above, and we omit the details.

\noindent Case (iii) $\mu\notin \R$. Since $-\mu f + M_{\sin y} (I + A_l^{-1}) f=0$, by taking the inner product with $B_{2,l}f$, we have $-\mu \langle f, B_{2,l} f\rangle_{L^2} + \langle M_{\sin y} B_{2,l}f, B_{2,l} f\rangle_{L^2} =0$. The imaginary part of this equality gives $(\Im \mu) \langle f, B_{2,l} f \rangle_{L^2}=0$, and thus, from the definition of $B_{2,l}$, we observe that $f=0$ if $|l|\geq 2$ and $f=constant$ if $|l|=1$.
On the other hand, if $|l|=1$ and $f=constant$ then $(I+A_l^{-1})f=B_{2,l} f=0$, which gives $-\mu f=0$. Thus $f=0$ since $\mu\ne 0$. The proof is complete.

\

The above result for $\hat{\Lambda}_l$ in $L^2(\T)$ is easily translated to $\hat{\Lambda}$ in $L_0^2(\T^2)$.
Indeed, 
\begin{align}
\hat{\Lambda} =B_1 B_2\,, \qquad B_1 = -i\partial_x M_{\sin y}\,, \qquad B_2 = (I+ A^{-1})\,,
\end{align}
and $B_1$ is closed symmetric and $B_2$ is bounded self-adjoint in $L^2_0 (\T^2)$, and ${\rm Ker}\, \hat{\Lambda} = {\rm Ker}\, B_2$. The operator $B_2$ is positive in $Y$, for so is $B_{2,l}$ in $Y_l$ for each $l\in \Z\setminus\{0\}$ with a uniform lower bound in $l$.
Proposition \ref{prop.kernel.Kol} therefore implies
\begin{cor}\label{cor.prop.kernel.Kol} $\sigma (\hat{\Lambda}) = \R$. Moreover, $\hat{\Lambda}$ in $L^2_0(\T^2)$ does not have eigenvalues in $\C\setminus\{0\}$.
\end{cor}

\noindent {\it Proof.} Let $f\in L^2_0 (\T^2)$ and set $f_l (y)=(\mathcal{P}_l f) (x,y) e^{-i l x}$. If $\zeta\in \C$ and $\Im \zeta\ne 0$ then $\zeta-l \hat{\Lambda}_l$ is invertible for any $l\in \Z\setminus\{0\}$ and $\omega_l = (\zeta-l\hat{\Lambda}_l)^{-1} f_l$ satisfies 
\begin{align*}
\|\omega_l\|_{L^2(\T)} \leq \frac{C}{|\Im \zeta|} \| f_l \|_{L^2(\T)} ~~ \text{ if }~ |l|\geq 2\,, \qquad \| \omega_l \|_{L^2(\T)} \leq C(|\Im \zeta|) \| f_l\|_{L^2 (\T)} ~~\text{ if } ~|l|=1\,.
\end{align*}
Here $C$ is independent of $l$ and $C(|\Im \zeta|)$ depends only on $|\Im \zeta|$ (the concrete dependence of $C(|\Im \zeta|)$ on $|\Im \zeta|$ is not needed in the argument below). 
Moreover, from $\zeta \omega_l - l M_{\sin y} B_{2,l} \omega_l=f_l$, we also have 
\begin{align*}
\| l M_{\sin y} \omega_l\|_{L^2 (\T)} \leq  \| l M_{\sin y} A_l^{-1} \omega_l \|_{L^2(\T)} + |\zeta | \, \| \omega_l \|_{L^2(\T)} + \| f_l \|_{L^2 (\T)} \leq C(\zeta) \| f_l \|_{L^2(\T)}
\end{align*}
with $C(\zeta)$ depending only on $\zeta$ and independent of $l$. Hence, $\omega (x,y) =\sum_{l\in \Z\setminus\{0\}} \omega_l (y) e^{i l x}$ satisfies $\omega \in L^2_0 (\T^2)$  and   $\partial_x M_{\sin y} \omega \in L_0^2 (\T^2)$. Clearly $\omega$ solves $(\zeta - \hat{\Lambda}) \omega = f$. The uniqueness is also shown by taking the Fourier series in $x$. Thus $\zeta$ belongs to the resolvent set of $\hat{\Lambda}$ in $L^2_0 (\T^2)$. 
This shows $\sigma (\hat{\Lambda}) \subset \R$. Since $\sigma (l\hat{\Lambda}_l) =[-l,l]$, we conclude that $\sigma (\hat{\Lambda})=\R$.
If $\zeta\in \R$ is an eigenvalue of $\hat{\Lambda}$ and $\omega \in D(\hat{\Lambda})$ is an eigenfunction then there exists $l\in \Z\setminus \{0\}$ such that $\omega_l (y) = (\mathcal{P}_l \omega) (x,y) e^{-ilx}$ is nontrivial. Since $\omega_l$ satisfies $(\zeta - l \hat{\Lambda}_l)\omega_l=0$, the number $\frac{\zeta}{l}$ must be an eigenvalue of $\hat{\Lambda}_l$ in $L^2 (\T)$, which is a contradiction. The proof is complete.

\

We can now apply Theorems \ref{thm.abstract.1}, which yields the following result.

\begin{thm}\label{thm.Kol.1} Let $L_{\alpha}$ be as in \eqref{def.L_alpha.Kol}. Then 
\begin{align}\label{thm.Kol.1.1} 
\lim_{|\alpha|\rightarrow \infty} \sup_{\lambda\in \R} \| (i\lambda -\mathbb{Q}L_{\alpha})^{-1}\|_{Y\rightarrow Y} 
=\lim_{|\alpha|\rightarrow \infty} \sup_{\lambda\in \R} \| \mathbb{Q} (i\lambda - L_{\alpha})^{-1}\|_{X\rightarrow X} =0\,.
\end{align} 
A similar result holds also for $L_{\alpha,l}$ for each $l\in \Z\setminus\{0\}$. 
\end{thm}

\subsection{Estimate with rate}\label{subsec.Kol.rate}

Theorem \ref{thm.Kol.1} does not give any estimates on the convergence rate.
To apply Theorem \ref{thm.abstract.2}  we focus on the study of $\hat{\Lambda}_l = M_{\sin y} (I + A_l^{-1})$ in $L^2 (\T)$. Note that
\begin{align*}
-A_l = T_l^*T_l\,, \qquad T_l=\partial_y - l
\end{align*}
and therefore,
\begin{align*}
[T_l,B_1] = M_{\cos y}\,.
\end{align*}
In particular, it is not difficult to show from $B_{2,l}=I + A_l^{-1}$,
\begin{align}
| \Im \langle A_l \phi, \hat{\Lambda}_l \phi \rangle_{L^2} | = | \Im \langle T_l \phi, T_l B_1 B_{2,l} \phi \rangle_{L^2} | & = | \Im \langle T_l \phi, B_1 T_l B_{2,l} \phi + M_{\cos y} B_{2,l} \phi \rangle_{L^2} | \nonumber \\
& = | \Im \langle T_l \phi, B_1 T_l A_l^{-1} \phi + M_{\cos y} B_{2,l} \phi \rangle_{L^2} | \nonumber \\
& \leq \| T \phi \|_{L^2} \| \big ( B_1 T_l A_l^{-1} + M_{\cos y} B_{2,l} \big ) \phi \|_{L^2}\,,
\end{align}
which ensures Assumption \ref{assum.Lambda.rate} (iii) (a)  with $B_3 = B_1 T_l A_l^{-1}+ M_{\cos y} B_{2,l}$ and $X=X_l=L^2 (\T)$. 
The result for $\hat{\Lambda}$ is obtained by the diagonalization $\hat{\Lambda} = \oplus_{l\in \Z\setminus\{0\}} l \hat{\Lambda}_l$. 
To simplify the notation we use the symbols $u$ and $v$ as scalar functions on $\T$ in this subsection (i.e.,  in this subsection $u$ and $v$ do not mean velocity fields).

Our goal is to show the coercive estimates of $\hat{\Lambda}_l$ as stated in Assumption \ref{assum.Lambda.rate}.
The main difficulty comes from the degeneracy of the critical points, i.e., the case when $(\sin y)'|_{y=y_\mu}=\cos y_\mu$ vanishes for the point $y_\mu\in \sin^{-1} \mu$.
This is the case $\mu=\pm1$. More precisely, the difficulty is to show \eqref{est.assum.Lambda.rate'} uniformly in a neighborhood of $\mu=\pm 1$, rather than on the exact points $\mu=\pm 1$.
Below we divide the regime of $\mu$ into three parts. The first part, discussed in Lemma \ref{lem.Lambda.Kol.rate.1},  corresponds to the case $|\mu|\geq 1$, though we can take $|\mu|$ slightly below  $1$ depending on the value $m$ in \eqref{est.assum.Lambda.rate'}. 
The second part, which is the core part of this section and discussed in Lemma \ref{lem.Lambda.Kol.rate.2}, is $\frac12\leq |\mu|<1$. The last part is $|\mu|<\frac12$ and  will be treated in Lemma \ref{lem.Lambda.Kol.rate.3}, where  the critical point is nondegenerate, while we need to handle the additional nonlocality due to the presence of the projection $\mathbb{Q}_l$ when $|l|=1$.

\begin{lem}\label{lem.Lambda.Kol.rate.1} There exist $\kappa\in (0,1)$ and $C>0$ such that the following statements hold for all $\delta\in (0,1]$ and $l\in \Z\setminus\{0\}$.
If  $\mu\in \R$ satisfies $1-\kappa \delta^2 \leq |\mu|\leq 1+\kappa \delta^2$ then 
\begin{align}\label{est.lem.Lambda.Kol.rate.1}
\delta^2 \| u \|_{L^2}^2  + \| M_{\cos y} B_{2,l} u \|_{L^2}^2 + \|(-A_l)^{-\frac12} u \|_{L^2}^2 \leq C \big ( \delta^{-2} \| (\mu-\hat{\Lambda}_l ) u \|_{L^2}^2 + \delta^4 \| (-A_l)^\frac12 u \|_{L^2}^2 \big )\,, \qquad u\in H^1 (\T)\,,
\end{align}
while if $|\mu|>1$ then 
\begin{align}\label{est.lem.Lambda.Kol.rate.1'}
\begin{split}
& (|\mu|-1)^2 \| B_{2,l} u \|_{L^2}^2  + (|\mu|-1) \| M_{\cos y} B_{2,l} u\|_{L^2}^2 \\
& \qquad \qquad \qquad \qquad +  |\mu| (|\mu|-1)  \|(-A_l)^{-\frac12} u \|_{L^2}^2  \leq C  \| (\mu-\hat{\Lambda}_l ) u \|_{L^2}^2 \,,  \qquad u\in H^1 (\T)\,.
\end{split}
\end{align}
\end{lem}

\noindent {\it Proof.} Set $f=(\mu-\hat{\Lambda}_l) u$ for $u\in H^1 (\T)$, i.e.,
\begin{align}\label{proof.lem.Lambda.Kol.rate.2}
(\mu  -M_{\sin y})  u - M_{\sin y} A_l^{-1} u = f
\end{align}
by the definition of $\hat{\Lambda}_l$. Setting $v= A_l^{-1} u$ and using $M_{\sin y}=M_{\sin y} - \mu  + \mu$ , \eqref{proof.lem.Lambda.Kol.rate.2} is also written as 
\begin{align}\label{proof.lem.Lambda.Kol.rate.3}
(\mu  -M_{\sin y})  (A_l + 1)  v - \mu v = f\,.
\end{align}
Note that $(A_l+1)v = B_{2,l} u$ by the definition. 
Below we fix $\delta\in (0,1)$ and take $\kappa>0$ sufficiently small.
Taking the inner product with $(A_l+1) v$ and by considering the real part, we obtain 
\begin{align}\label{proof.lem.Lambda.Kol.rate.4}
\int_0^{2\pi} \frac{\mu-\sin y}{\mu} |(A_l+1) v|^2 \dd y + \| \partial_y v\|_{L^2}^2 + (l^2-1) \| v\|_{L^2}^2 = \frac{1}{\mu} \Re \langle f\,, (A_l + 1) v\rangle_{L^2}\,.
\end{align}
If $0<1-\kappa\delta^2 \leq \mu\leq 1+\kappa \delta^2$ then we have
\begin{align}\label{proof.lem.Lambda.Kol.rate.5}
& \int_0^{2\pi} \frac{1-\sin y}{\mu} |(A_l+1) v|^2 \dd y + \| \partial_y v\|_{L^2}^2 + (l^2-1)\| v\|_{L^2}^2 \nonumber \\
& \qquad = \frac{1-\mu}{\mu} \int_0^{2\pi}  |(A_l+1) v|^2 \dd y  + \Re  \frac{1}{\mu} \langle f\,, (A_l + 1) v\rangle_{L^2}\nonumber \\
& \qquad \leq \frac{\kappa\delta^2}{\mu} \int_0^{2\pi}  |(A_l+1) v|^2 \dd y  +  \frac{1}{\mu}\|  f\|_{L^2} \| (A_l + 1) v\|_{L^2} \,.
\end{align}
Thus, we have
\begin{align}\label{proof.lem.Lambda.Kol.rate.7}
\begin{split}
& \int_0^{2\pi} (1-\sin y) |(A_l+1) v|^2 \dd y + \| \partial_y v\|_{L^2}^2 + (l^2-1) \| v\|_{L^2}^2  \\
& \qquad \leq C\kappa \delta^2  \| (A_l + 1) v \|_{L^2}^2  +  C\|  f\|_{L^2} \| (A_l+1) v\|_{L^2}\,.
\end{split}
\end{align}
Next we compute 
\begin{align}
\int_0^{2\pi} (1-\sin y) |(A_l+1) v|^2 \dd y & = \int_{|y-\frac12\pi|\geq \delta} \ldots \dd y + \int_{|y-\frac12\pi|<\delta}\ldots \dd y \nonumber \\
& \geq  \int_{|y-\frac12\pi|\geq \delta} \ldots \dd y\nonumber \\
& \geq C \delta^2 \| (A_l+1) v\|_{L^2 (\{|y-\frac12\pi|\geq \delta\})}^2 \nonumber \\
& \geq C \delta^2 \| (A_l + 1) v\|_{L^2}^2 - C \delta^3 \| (A_l +1 ) v\|_{L^\infty}^2\,.
\end{align}
Thus, from $\| (A_l + 1) v\|_{L^\infty} \leq C \| u\|_{L^\infty}$ with $C$ independent of $l$, and \eqref{proof.lem.Lambda.Kol.rate.5}-\eqref{proof.lem.Lambda.Kol.rate.7}, we deduce that
\begin{align*}
& \delta^2 \| (A_l +1 ) v\|_{L^2}^2 + \| \partial_y v\|_{L^2}^2 + (l^2-1) \| v\|_{L^2}^2 \\
& \quad  \leq C \bigg ( \kappa \delta^2 \| (A_l + 1) v \|_{L^2}^2 + \delta^3 \| u\|_{L^\infty}^2 + \| f\|_{L^2} \| (A_l+1) v\|_{L^2}\bigg ) \,,
\end{align*}
and if $\kappa>0$ is small enough but independently of $\mu$ and $\delta$, then
\begin{align*}
\delta^2 \| (A_l +1 ) v\|_{L^2}^2 + \| \partial_y v\|_{L^2}^2 + (l^2-1) \| v\|_{L^2}^2  \leq C \bigg (  \delta^3 \| u\|_{L^\infty}^2 + \delta^{-2} \| f\|_{L^2}^2 \bigg ) \,.
\end{align*}
From $A_l v=u$ we finally obtain 
\begin{align}\label{proof.lem.Lambda.Kol.rate.8}
\delta^2 \| u \|_{L^2}^2 + \| (-A_l)^{-\frac12} u \|_{L^2}^2 \leq C \bigg (\delta^{-2}\| (\mu - \hat{\Lambda}_l ) u \|_{L^2}^2   + \delta^3 \| u \|_{L^\infty}^2   \bigg ) \,.
\end{align}
The argument is the same as above for the case  $-1-\kappa\delta^2 \leq \mu\leq -1 + \kappa \delta^2 <0$, and we have \eqref{proof.lem.Lambda.Kol.rate.8} also in this case. The details are  omitted here.  
%Now we take $f=(\mu + i\epsilon - \hat{\Lambda}_l) u$ for $u\in H^1(\T)$. 
%Then $w=u\in H^1(\T)$ and we obtain from \eqref{proof.lem.Lambda.Kol.rate.8},
%\begin{align}\label{proof.lem.Lambda.Kol.rate.10}
%\delta^2 \| u \|_{L^2}^2 +\| (-A_l)^{-\frac12} u \|_{L^2}^2 \leq C \bigg ( \delta^{-2} \|  (\mu + i\epsilon - \hat{\Lambda}_l) u \|_{L^2}^2 + \delta^3 \| u \|_{L^\infty}^2 \bigg )\,.
%\end{align}
%Note that $C$ is independent of $\epsilon>0$, $|\mu|\in [1-\kappa \delta^2,1+\kappa \delta^2]$, and $\delta$.
%Letting $\epsilon\rightarrow 0$, we have arrived at
%\begin{align}\label{proof.lem.Lambda.Kol.rate.11}
%\delta^2 \| u \|_{L^2}^2 +\| (-A_l)^{-\frac12} u \|_{L^2}^2 \leq C \bigg ( \delta^{-2} \|  (\mu - \hat{\Lambda}_l) u \|_{L^2}^2 + \delta^3 \| u \|_{L^\infty}^2 \bigg )\,, \qquad u\in H^1 (\T)\,,
%\end{align}
%for $1-\kappa \delta^2 \leq |\mu|\leq 1+\kappa \delta^2$ and $\delta\in (0,\frac{1}{100})$.
The estimate of $\| M_{\cos x} B_{2,l} u\|_{L^2}^2$ follows from \eqref{proof.lem.Lambda.Kol.rate.7} and \eqref{proof.lem.Lambda.Kol.rate.8} by the inequality 
$$\cos^2 y =(1-\sin^2 y) \leq 2 (1\pm \sin y)\,.$$ 
Then it suffices to apply the interpolation inequality $\|u\|_\infty^2 \leq C \| u \|_{H^1}\| u \|_{L^2}$ and $\|u\|_{H^1}\leq \| (-A_l)^\frac12 u \|_{L^2}$ to obtain \eqref{est.lem.Lambda.Kol.rate.1}.
Estimate \eqref{est.lem.Lambda.Kol.rate.1'} for the case $|\mu| > 1$ easily follows from  the identity \eqref{proof.lem.Lambda.Kol.rate.4}. Indeed, \eqref{proof.lem.Lambda.Kol.rate.4} is written from $(A_l+1) v= B_{2,l} u$,
\begin{align*}
(1-\frac{1}{|\mu|} ) \| B_{2,l} u \|_{L^2}^2 +  \frac{1}{|\mu|} \int_0^{2\pi} \big (1-\frac{|\mu|}{\mu}  \sin y \big ) \, |B_{2,l} u |^2 \dd y + \| \partial_y v\|_{L^2}^2 + (l^2-1) \| v\|_{L^2}^2 = \frac{1}{\mu} \Re \langle f\,, B_{2,l} u \rangle_{L^2}\,,
\end{align*}
which gives \eqref{est.lem.Lambda.Kol.rate.1'} for $|\mu|>1$. 
The details are omitted here. The proof is complete.

\

The coercive estimate for $|\mu|<1-\kappa \delta^2$ is more delicate,
especially when $\kappa \delta^2 <|\mu\pm 1| \leq o(1)$ as $\delta\rightarrow 0$
due to the degeneracy around the points such that $(\sin y)'=\cos y=0$ and the nonlocality.
To overcome the difficulty we apply a contradiction argument.
The contradiction argument is useful since it enables us to focus on the functions which concentrate around the critical points, by which the nonlocality is reduced since the nonlocal operator has a {\it smoothing} effect and thus becomes a small perturbation of the local operator for such functions. 
The following lemma, which requires a long proof, is the core result of this subsection.
\begin{lem}\label{lem.Lambda.Kol.rate.2}
Let $\kappa\in (0,1)$ be the number in Lemma \ref{lem.Lambda.Kol.rate.1}.
There exists $C>0$  such that if $\delta\in (0,1]$, $l\in \Z\setminus\{0\}$, and $\mu\in \R$ with $\frac12\leq |\mu| < 1-\kappa \delta^2$, then 
\begin{align}\label{est.lem.Lambda.Kol.rate.2}
\delta^2 \| u \|_{L^2}^2 + \| (-A_l)^{-\frac12} u \|_{L^2}^2 + \frac{1}{1-|\mu|} \| A_l^{-1} u\|_{L^2}^2 \leq C \bigg ( \delta^{-2} \| (\mu-\hat{\Lambda}_l) u \|_{L^2}^2 + \frac{\delta^6}{1-|\mu|} \| (-A_l)^\frac12 u\|_{L^2}^2  \bigg ) \,, \quad u\in H^1 (\T) \,,
\end{align}
and 
\begin{align}\label{est.lem.Lambda.Kol.rate.2'}
\| M_{\cos y} B_{2,l} u \|_{L^2}^2 \leq C \bigg ( \delta^{-2} \| (\mu-\hat{\Lambda}_l) u \|_{L^2}^2 + \delta^2 (1-|\mu| ) \, \| (-A_l)^\frac12 u\|_{L^2}^2 \bigg ) \,, \quad u\in H^1 (\T) \,.
\end{align}
\end{lem}

The proof consists of several steps. We first consider \eqref{est.lem.Lambda.Kol.rate.2}.
%\begin{lem}\label{lem.Lambda.Kol.rate.2-1} The inequality \eqref{est.lem.Lambda.Kol.rate.2} holds.
% \end{lem}

\noindent {\it Proof of  \eqref{est.lem.Lambda.Kol.rate.2}.} Since $\mu$ is a real number and $\hat{\Lambda}_l$ and $A_l$ preserve the real valued, without loss of generality it suffices to show \eqref{est.lem.Lambda.Kol.rate.2} for real valued functions. We may also assume that $\mu$ is positive, for the case $\mu<0$ is proved in the same manner. By the density argument it suffices to show the claim for $u\in H^2 (\T)$, rather than $u\in H^1(\T)$. We make use of a contradiction argument.
Suppose that  the estimate 
\begin{align}\label{proof.prop.Lambda.Kol.rate.2.1}
\begin{split}
& \delta^2 \| u \|_{L^2}^2 + \| (-A_l)^{-\frac12} u \|_{L^2}^2 +  \frac{1}{1-\mu} \| A_l^{-1} u\|_{L^2}^2  \leq C \bigg ( \delta^{-2} \| (\mu - \hat{\Lambda}_l) u \|_{L^2}^2 + \frac{\delta^6}{1-\mu} \| (-A_l)^\frac12  u\|_{L^2}^2  \bigg )\,, \\
& \delta\in (0,1]\,, \quad l\in \Z\setminus\{0\}\,, \quad \frac12\leq \mu<1-\kappa |\delta|^2\,, \quad u\in H^2 (\T; \R)
\end{split}
\end{align}
does not hold.
Then there exist $\{\delta_n, l_n, \mu_n\}_{n\in \N}$, $\delta_n\in (0,1]$, $l_n\in \Z\setminus\{0\}$, $\mu_n\in  [\frac12, 1-\kappa \delta_n^2)$, and $\{u_n\}\subset H^2 (\T; \R)$ such that 
\begin{align*}
& \lim_{n\rightarrow \infty} \delta_n=\delta_\infty\in [0,1]\,, \quad \lim_{n\rightarrow \infty} l_n = l_\infty\in \{\pm \infty\}\cup \Z\setminus\{0\}\,,\\
& \lim_{n\rightarrow \infty}\mu_n=\mu_\infty \in [\frac12, 1-\kappa \delta_\infty^2]\,,
\end{align*}
and 
\begin{align}\label{proof.prop.Lambda.Kol.rate.2.2}
\begin{split}
& \delta_n^2 \| u_n \|_{L^2}^2 + \| (-A_{l_n})^{-\frac12} u_n \|_{L^2}^2 + \frac{1}{1-\mu_n} \| A_{l_n}^{-1} u_n \|_{L^2}^2 =1\,, \\
& \lim_{n\rightarrow \infty}  \bigg ( \delta_n^{-2}  \|  (\mu_n - \hat{\Lambda}_{l_n}) u_n \|_{L^2}^2 + \frac{\delta_n^6}{1-\mu_n} \| (-A_{l_n})^\frac12   u_n \|_{L^2}^2 \bigg ) =0\,.
\end{split}
\end{align}
Set
\begin{align}\label{proof.prop.Lambda.Kol.rate.2.3}
f_n = \delta_n^{-1} (\mu_n  - \hat{\Lambda}_{l_n} ) u_n\,, \qquad v_n = A_{l_n}^{-1} u_n\,,
\end{align}
and then $v_n$ satisfies 
\begin{align}\label{proof.prop.Lambda.Kol.rate.2.4}
(\mu_n - M_{\sin y}) \big ( A_{l_n} + 1) v_n - \mu_n  v_n = \delta_n f_n\,.
\end{align}
The normalized condition in \eqref{proof.prop.Lambda.Kol.rate.2.2} implies $\delta_n^2 \| u_n \|_{L^2}^2 + \|(-A_{l_n})^\frac12 v_n \|_{L^2}^2 + \frac{1}{1-\mu_n} \| v_n \|_{L^2}^2 =1$, and thus, from the integration by parts,
\begin{align}\label{proof.prop.Lambda.Kol.rate.2.2'}
\delta_n^2 \| u_n \|_{L^2}^2 + \|\partial_y v_n \|_{L^2}^2 + l_n^2 \| v_n \|_{L^2}^2 +  \frac{1}{1-\mu_n} \| v_n \|_{L^2}^2 =1\,.
\end{align}
Note that we have also from \eqref{proof.prop.Lambda.Kol.rate.2.2'} and the interpolation inequality $\| v_n \|_{L^\infty} \leq C \| v_n \|_{H^1}^\frac12 \| v_n \|_{L^2}^\frac12$ that 
\begin{align}\label{proof.prop.Lambda.Kol.rate.2.2''}
\sup_n \frac{\| v_n \|_{L^\infty}}{|1-\mu_n|^\frac14} <\infty\,.
\end{align}
Since $\sup_n \| v_n\|_{H^1}<\infty$, we may assume that, after taking a suitable subsequence, $\{v_n\}$ converges to $v_\infty$ weakly in $H^1(\T; \R)$, and thus, strongly in $L^2(\T;\R)$. 
First we exclude the possibility $\delta_\infty>0$. Indeed, in this case we have $\frac{\|(-A_{l_n})^\frac12 u_n\|_{L^2}^2}{1-\mu_n}\rightarrow 0$ by \eqref{proof.prop.Lambda.Kol.rate.2.2}, which implies 
\begin{align*}
\delta_n^2 \| u_n \|_{L^2}^2 + \| (-A_{l_n})^\frac12 v_n \|_{L^2}^2 + \frac{1}{|1-\mu_n|} \| v_n \|_{L^2}^2 \leq C \frac{\| (-A_{l_n})^\frac12 u_n \|_{L^2}^2}{1-\mu_n} \rightarrow 0\,.
\end{align*}
This contradicts with the normalized condition in \eqref{proof.prop.Lambda.Kol.rate.2.2'}.
Thus it suffices to consider the case $\delta_\infty=0$.
Let $S_{\mu_n}=\{y\in \T~|~\sin y = \mu_n\}$ be the set of critical points.
Then we have from \eqref{proof.prop.Lambda.Kol.rate.2.4},
\begin{align}\label{proof.prop.Lambda.Kol.rate.2.4'}
\mu_n v_n (y_\mu) + \delta_n f_n (y_\mu)=0\,, \qquad {\rm for~any}~~y_\mu\in S_{\mu_n} \,.
\end{align}
This fact plays an important role in the analysis.

Let us start from the following claim.

\noindent {\bf Step 1: $\displaystyle \lim_{n\rightarrow \infty} \delta_n \| u_n \|_{L^2}=0$.} 

The Taylor expansion implies for any $y_{\mu_n} \in S_{\mu_n}$, 
\begin{align}\label{eq.sin}
\sin y & = \mu_n + \cos y_{\mu_n} \, (y-y_{\mu_n}) + R_n (y) (y-y_{\mu_n})^2 \,, \qquad |R_n (y)|\leq C \,, \qquad y\in [-\frac{\pi}{2}, \frac{3\pi}{2}]\,,
\end{align}
and 
\begin{align}\label{eq.cos}
|1-\mu_n|^\frac12 \leq |\cos y_{\mu_n}| = \sqrt{1-\mu_n^2} \leq 2^\frac12|1-\mu_n|^\frac12\,.
\end{align}
Here $C$ is independent of $n$. 
Let $\kappa_1>0$ be fixed but arbitrary small number. We decompose the interval $[-\frac{\pi}{2},\frac{3\pi}{2}]$ into $I_n$ and $I_n^c=[-\frac{\pi}{2},\frac{3\pi}{2}]\setminus I_n$, where 
\begin{align}
I_n & = \big \{y\in [-\frac{\pi}{2},\frac{3\pi}{2}]~ \big |~{\rm dist}\, (y, S_{\mu_n}) \leq \frac{\delta_n^2}{\kappa_1 |1-\mu_n|^\frac12} \big \}\,.
\end{align}
Then we have from $\sup_n \delta_n \| B_{2,l_n} u_n \|_{L^2}\leq C$ by \eqref{proof.prop.Lambda.Kol.rate.2.2'},
\begin{align*}
\delta_n \| B_{2,l_n} u_n \|_{L^2(I_n)} \leq \delta_n |I_n|^\frac12 \| B_{2,l_n} u_n \|_{L^\infty}  
& \leq \frac{C\delta_n^2}{\kappa_1^\frac12 |1-\mu_n|^\frac14} \| B_{2,l_n} u_n \|_{H^1}^\frac12 \| B_{2,l_n} u_n \|_{L^2}^\frac12 \\
& \leq \frac{C \delta_n^\frac32}{\kappa_1^\frac12 |1-\mu_n|^\frac14} \| B_{2,l_n} u_n \|_{H^1}^\frac12 \\
& \leq C \Big ( \frac{\delta_n^3}{\kappa_1 |1-\mu_n|^\frac12} \| (-A_{l_n})^\frac12 u_n \|_{L^2}\Big )^\frac12  \rightarrow 0 \quad (n\rightarrow \infty\,, ~ {\rm by~\eqref{proof.prop.Lambda.Kol.rate.2.2}}) \,.
\end{align*}
Next, \eqref{proof.prop.Lambda.Kol.rate.2.4} gives 
\begin{align}\label{proof.prop.Lambda.Kol.rate.2.5}
B_{2,l_n} u_n & = \frac{\mu_n v_n+\delta_n f_n}{\mu_n - \sin y}\,.
\end{align}
We decompose $I_n^c$ as $I_n^c = (I_n^c \cap [-\frac{\pi}{2},\frac{\pi}{2}] )\cup (I_n^c \cap [\frac{\pi}{2}, \frac{3\pi}{2}])=:I_{n,1}^c \cup I_{n,2}^c$, and we find that there exists $C>0$ independent of $n$ such that 
\begin{align}\label{est.sin}
\begin{split}
\sum_{j=1,2} \| \frac{1}{\mu_n - \sin y} \|_{L^2 (I_{n,j}^c)} \leq \frac{C\kappa_1^\frac12}{\delta_n (1-\mu_n)^\frac14}\,, \qquad 
\sum_{j=1,2} \| \frac{1}{\mu_n - \sin y} \|_{L^\infty (I_{n,j}^c)} \leq \frac{C}{\delta_n^2}\,.
\end{split}
\end{align}
Then we have from \eqref{proof.prop.Lambda.Kol.rate.2.5} with \eqref{est.sin} and \eqref{proof.prop.Lambda.Kol.rate.2.2''},
\begin{align*}
\delta_n \| B_{2,l_n} u_n \|_{L^2(I_{n,1}^c)} & \leq C \delta_n \|v_n \|_{L^\infty} \| \frac{1}{\mu_n-\sin y} \|_{L^2(I_{n,1}^c)} + C\delta_n^2  \| f_n \|_{L^2} \| \frac{1}{\mu_n-\sin y} \|_{L^\infty (I_{n,1}^c)} \\
& \leq \frac{C\kappa_1^\frac12\| v_n \|_{L^\infty}}{|1-\mu_n|^\frac14} + C \| f_n \|_{L^2}\\
& \leq C \kappa_1^\frac12 + C \| f_n \|_{L^2}\,.
\end{align*}
The same estimate holds for $\delta_n \| B_{2,l_n} u_n \|_{L^2(I_{n,2}^c)}$.
Thus we obtain 
$\delta_n \| B_{2,l_n} u_n \|_{L^2 (I_n^c)} \leq  C \kappa_1^\frac12 + C \| f_n \|_{L^2}$.
Since $\kappa_1>0$ is arbitrary and $\displaystyle \lim_{n\rightarrow \infty} \| f_n \|_{L^2}=0$, we have 
$\displaystyle \limsup_{n\rightarrow \infty} \delta_n \| B_{2,l_n} u_n \|_{L^2} =0$,
which gives $\displaystyle \lim_{n\rightarrow \infty} \delta_n \| u_n \|_{L^2}=0$ by the relation $B_{2,l_n} u_n = u_n + v_n$. 

\

\noindent {\bf Step 2: $\displaystyle \lim_{n\rightarrow \infty} \frac{\mu_n |v_n (y_n)|}{|1-\mu_n|^\frac14} = \lim_{n\rightarrow \infty} \frac{\delta_n |f_n (y_n)|}{|1-\mu_n|^\frac14}=0$ for $y_n\in S_{\mu_n}$.} 

\noindent In virtue of \eqref{proof.prop.Lambda.Kol.rate.2.4'} it suffices to show $\displaystyle \lim_{n\rightarrow \infty} \frac{\delta_n |f_n (y_n)|}{|1-\mu_n|^\frac14}=0$ for $y_n\in S_{\mu_n}$. As in the previous step, let $y_{n,1}$ be the unique point of $S_{\mu_n}$ such that $y_{n,1}\in (0,\frac{\pi}{2})$. 
Let $\kappa_2>0$ be fix but arbitrary small number and set $II_{n,1}=(y_{n,1}, y_{n,1}+\kappa_2 \frac{\delta_n^2}{|1-\mu_n|^\frac12}]$. There exists a point $\tilde y_{n,1}\in II_{n,1}$ such that $ \frac{\kappa_2 \delta_n^2}{|1-\mu_n|^\frac12} |f_n (\tilde y_{n,1})|^2 \leq \| f_n \|_{L^2(II_{n,1})}^2$, and thus, we have $\frac{\delta_n | f_n (\tilde y_{n,1})|}{|1-\mu_n|^\frac14} \leq \frac{\|f_n \|_{L^2}}{\kappa_2^\frac12}$.
Then we see 
\begin{align*}
\frac{\delta_n |f_n (y_{n,1})|}{|1-\mu_n|^\frac14}\leq  \frac{\delta_n |f(y_{n,1})- f(\tilde y_{n,1})|}{|1-\mu_n|^\frac14} + \frac{\delta_n |f(\tilde y_{n,1})|}{|1-\mu_n|^\frac14} & \leq \frac{\delta_n |II_{n,1}|^\frac12}{|1-\mu_n|^\frac14} \| \partial_y f_n \|_{L^2(II_{n,1})} + \frac{\|f_n \|_{L^2}}{\kappa_2^\frac12}\\
& \leq \frac{C \kappa_2^\frac12 \delta_n^2}{|1-\mu_n|^\frac12} \| \partial_y f_n \|_{L^2(II_{n,1})} + \frac{\|f_n \|_{L^2}}{\kappa_2^\frac12}\,.
\end{align*}
To estimate the $H^1$ norm of $f_n$ we observe that $\delta_n \partial_y f_n$ is written from \eqref{proof.prop.Lambda.Kol.rate.2.4} as 
\begin{align}\label{proof.prop.Lambda.Kol.rate.2.6}
\delta_n \partial_y f_n & = (\mu_n  - \sin y ) (\partial_y u_n + \partial_y v_n) - \cos y\, (A_{l_n} + 1) v_n  - \mu_n \partial_y v_n \nonumber \\
\begin{split}
& =  - M_{\cos y} B_{2,l_n} u_n  - M_{\sin y} \, \partial_y v_n +   (\mu_n  - M_{\sin y} ) \partial_y u_n \,.
\end{split}
\end{align}
From $\cos y = \cos y_{n,1} + O((y-y_{n,1})^2)$ and \eqref{eq.cos} we observe that $|\cos y|\leq C |1-\mu_n|^\frac12$ for $y\in II_{n,1}$, and also $|\mu_n-\sin y|\leq C \delta_n^2$ for $y\in II_{n,1}$. Thus we have 
\begin{align}\label{proof.prop.Lambda.Kol.rate.2.7}
\| \delta_n \partial_y f_n \|_{L^2(II_{n,1})} \leq C |1-\mu_n|^\frac12 \| B_{2,l_n} u_n\|_{L^2} + C \| \partial_y v_n \|_{L^2} + C \delta_n^2 \| \partial_y u_n \|_{L^2}
\end{align}
Therefore, we arrive at 
\begin{align*}
\frac{\delta_n |f_n (y_{n,1})|}{|1-\mu_n|^\frac14}\leq C \kappa_2^\frac12 \delta_n \| u_n \|_{L^2} + \frac{C\kappa_2^\frac12\delta_n}{|1-\mu_n|^\frac12} \| \partial_y v_n \|_{L^2} + \frac{C \kappa_2^\frac12 \delta_n^3}{|1-\mu_n|^\frac12} \| \partial_y u_n \|_{L^2} + \frac{\|f_n \|_{L^2}}{\kappa_2^\frac12}\,,
\end{align*}
which shows from Step 1, \eqref{proof.prop.Lambda.Kol.rate.2.2}, and $1-\mu_n \geq \kappa \delta_n^2$ that $\displaystyle \limsup_{n\rightarrow \infty} \frac{\delta_n |f_n(y_{n,1})|}{|1-\mu_n|^\frac14} \leq C (\frac{\kappa_2}{\kappa})^\frac12$.
Since $\kappa_2>0$ is arbitrary, the claim is proved. The estimate for the point $y_{n,2}\in S_{\mu_n}\cap (\frac{\pi}{2},\pi)$ is proved in the same manner.

\

\noindent {\bf Step 3:} {\it Estimate of \, $\| v_n \|_{L^2}^2$.}

\noindent
Let $y_{n,2},y_{n,3}\in S_{\mu_n}$ be the unique critical points such that $y_{n,2}\in (\frac{\pi}{2},\pi)$ and $y_{n,3}\in (2\pi,\frac{5\pi}{2})$ (that is, $y_{n,3} =y_{n,1} + 2\pi$).
Then we have 
\begin{align*}
\int_{y_{n,2}}^{y_{n,3}} (\mu_n - \sin y) \big (\partial_y^2 v_n - l_n^2 v_n +v_n)  \, v_n \, \dd y - \mu_n \int_{y_{n,2}}^{y_{n,3}} v_n^2 \, \dd y = \delta_n \int_{y_{n,2}}^{y_{n,3}} f_n \, v_n\, \dd y \,,
\end{align*}
and then, by the integration by parts, 
\begin{align}\label{proof.prop.Lambda.Kol.rate.2.8}
& \int_{y_{n,2}}^{y_{n,3}} (\mu_n -\sin y) \big (|\partial_y v_n |^2 + (l_n^2-1) v_n^2 \big ) \dd y  +  \int_{y_{n,2}}^{y_{n,3}} (\mu_n - \frac12\sin y)\, v_n^2 \dd y  \nonumber \\
&  =\frac12 \sum_{j=2,3} |v_n (y_{n,j})|^2(-1)^{j+1}  \cos y_{n,j}  -\delta_n \int_{y_{n,2}}^{y_{n,3}} f_n\, v_n \, \dd y \nonumber \\
&  \leq C |1-\mu_n| \big ( \sum_{j=1,2} \frac{|v_n (y_{n,j})|^2}{|1-\mu_n|^\frac12} + \|f_n \|_{L^2} \big)\quad ({\rm by ~\eqref{eq.cos}, \eqref{proof.prop.Lambda.Kol.rate.2.2'}, ~and}~1-\mu_n \geq \kappa \delta_n^2)\,.
\end{align}
Note that $\mu_n - \sin y\geq 0$ in $(y_{n,2},y_{n,3})$.
When $\frac34 < \mu_\infty\leq 1$ we have the bound $\mu_n-\frac12 \sin y\geq \frac14$ for large $n$.
Since the norm over the interval $[y_{n,2},y_{n,3}]$ is the same as the norm over $[-\frac{\pi}{2},\frac{3\pi}{2}]\setminus (y_{n,1},y_{n,2})$ for $2\pi$ periodic functions, we have  when $\frac34 < \mu_\infty \leq 1$, by using \eqref{proof.prop.Lambda.Kol.rate.2.8} with Step 2 and \eqref{proof.prop.Lambda.Kol.rate.2.2},
\begin{align}\label{proof.prop.Lambda.Kol.rate.2.9}
\lim_{n\rightarrow \infty} \frac{1}{1-\mu_n} \| v_n\|_{L^2 ([-\frac{\pi}{2},\frac{3\pi}{2}]\setminus (y_{n,1},y_{n,2}))}^2 =0\,.
\end{align}

\

\noindent {\bf Step 4:} {\it Estimate of \, $\| (-A_{l_n})^\frac12 v_n \|_{L^2}^2$.}

\noindent The integration by parts and \eqref{proof.prop.Lambda.Kol.rate.2.4} yield for any $\varphi \in H^2 (\T)$, 
\begin{align}\label{proof.prop.Lambda.Kol.rate.2.10}
\langle \partial_y v_n, \partial_y \varphi \rangle_{L^2} + (l_n^2-1) \langle v_n, \varphi \rangle_{L^2}  = - \langle B_{2,l_n} u_n, \varphi \rangle_{L^2} = - \langle \frac{\mu_n v_n + \delta_n f_n}{\mu_n-\sin y}, \varphi \rangle_{L^2} \,.
\end{align}
Note that the term $\frac{\mu_n v_n+f_n}{\mu_n-\sin y}$ belongs to $L^2$ in virtue of \eqref{proof.prop.Lambda.Kol.rate.2.4'} and the Hardy inequality.
Let us estimate the right-hand side of \eqref{proof.prop.Lambda.Kol.rate.2.10}.
\begin{align*}
- \langle \frac{\mu_n v_n + f_n}{\mu_n-\sin y}, \varphi \rangle_{L^2} & = - \int_{-\frac{\pi}{2}}^{\frac{\pi}{2}} \frac{\mu_n v_n +\delta_n f_n}{\mu_n-\sin y}\,  \varphi \, \dd y  - \int_{\frac{\pi}{2}}^\frac{3\pi}{2} \frac{\mu_n v_n + \delta_n f_n}{\mu_n-\sin y}\,  \varphi \, \dd y   =: J_1 [\varphi] + J_2 [\varphi]  \,.
\end{align*}
We are  interested in the case $\varphi=v_n$.
Let us estimate $J_1$, for which detailed computation is needed. 
Let $y_{n,1}$ be the unique point of $S_{\mu_n} \cap (0,\frac{\pi}{2})$.
Let $\kappa_3>0$ be fixed but arbitrary small number.
Take $b_n\in (\frac{1}{10},10)$ such that 
$$\sin a_{n,1} = \sin a_{n,1}'\,, \qquad a_{n,1}=y_{n,1}+ \kappa_3 \frac{\delta_n^2}{|1-\mu_n|^\frac12}\,, \quad  a_{n,1}' = y_{n,1} -b_n \kappa_3 \frac{\delta_n^2}{|1-\mu_n|^\frac12}\,.$$
Set $III_{n,1} = [y_{n,1}-b_n \kappa_3\frac{\delta_n^2}{|1-\mu_n|^\frac12}, y_{n,1} + \kappa_3 \frac{\delta_n^2}{|1-\mu_n|^\frac12}]$ and 
$$IV_{n,1} = [y_{n,1}- \kappa_3 |1-\mu_n|^\frac12, \, y_{n,1} - b_n \kappa_3 \frac{\delta_n^2}{|1-\mu_n|^\frac12}] \cup [y_{n,1}+\kappa_3\frac{\delta_n^2}{|1-\mu_n|^\frac12}, \, y_{n,1} + \kappa_3 |1-\mu_n|^\frac12]\,.$$
Then we have for $a_{n,2} = y_{n,1} + \kappa_3 |1-\mu_n|^\frac12$ and $a_{n,2}' = y_{n,1} - \kappa_3 |1-\mu_n|^\frac12$,
\begin{align}\label{proof.prop.Lambda.Kol.rate.2.11}
& \big |\int_{IV_{n,1}} \frac{1}{(\mu_n-\sin y)} \cos y \dd y \big | \nonumber \\
& = \big |\int_{\sin IV_{n,1}} \frac{1}{\mu_n-s} \dd s \big | \nonumber  \\
& = \big | -  \log (\sin a_{n,2}-\mu_n) +  \log (\mu_n-\sin a_{n,1}) - \log (\mu_n -\sin a_{n,1}')  + \log (\mu_n - \sin a_{n,2}') \big | \nonumber \\
& =\big | \log \frac{\mu_n- \sin a_{n,2}'}{\sin a_{n,2}-\mu_n} \big | \leq C\,.
\end{align}
Here we have used the fact that $\sin a_{n,2}-\mu_n$ and $\mu_n - \sin a_{n,2}'$ are approximated by $\cos y_{n,1} \, \kappa_3 |1-\mu_n|^\frac12 \sim \kappa_3 O(|1-\mu_n|)$; recall \eqref{eq.sin} and \eqref{eq.cos}.
Taking this into account, we can decompose $J_1$ as 
\begin{align*}
J_1 = - \Big ( \int_{III_{n,1}} + \int_{IV_{n,1}} + \int_{[-\frac{\pi}{2},\frac{\pi}{2}]\setminus (III_{n,1} \cup IV_{n,1})} \Big ) = J_{1,1} + J_{1,2} + J_{1,3}\,. 
\end{align*}   
Since $y_{n,1}\in III_{n,1}$, we have from $|\mu_n-\sin y|\geq \frac{\cos y_{n,1}}{C}|y-y_{n,1}|$ for $y\in III_{n,1}$ and the Hardy inequality,
\begin{align*}
|J_{1,1}[\varphi]| & \leq \frac{C}{\cos y_{n,1}}\| \partial_y (\mu_n v_n + \delta_n f_n ) \|_{L^2(III_{n,1})} \| \varphi \|_{L^2(III_{n,1})}\\
& \leq \frac{C}{|1-\mu_n|^\frac12} (\| \partial_y v_n \|_{L^2} + \|\delta_n \partial_y f_n \|_{L^2(III_{n,1})} ) |III_{n,1}|^\frac12 \| \varphi \|_{L^\infty}\\
& \leq \frac{C\kappa_3^\frac12\delta_n}{|1-\mu_n|^\frac34} (1 + \|\delta_n \partial_y f_n \|_{L^2(III_{n,1})} ) \| \varphi \|_{L^\infty}\,.
\end{align*}
Here we have used the normalized condition \eqref{proof.prop.Lambda.Kol.rate.2.2'} for $\|\partial_y v_n \|_{L^2}$. The norm $\| \delta_n \partial_y f_n \|_{L^2 (III_{n,1})}$ is bounded as \eqref{proof.prop.Lambda.Kol.rate.2.7} by the definition of $III_{n,1}$, and thus, in virtue of \eqref{proof.prop.Lambda.Kol.rate.2.2}, we conclude that 
\begin{align}\label{proof.prop.Lambda.Kol.rate.2.12}
\displaystyle \sup_n \frac{\delta_n}{|1-\mu_n|^\frac12} \| \delta_n \partial_y f_n \|_{L^2(III_{n,1})} \leq C<\infty\,,
\end{align}
which gives 
\begin{align}\label{proof.prop.Lambda.Kol.rate.2.13}
|J_{1,1}[\varphi]| \leq C \kappa_3^\frac12 \frac{\| \varphi \|_{L^\infty}}{|1-\mu_n|^\frac14}\,,
\end{align}
and in particular, by using \eqref{proof.prop.Lambda.Kol.rate.2.2''},
\begin{align}\label{proof.prop.Lambda.Kol.rate.2.14}
|J_{1,1}[v_n]| \leq C \kappa_3^\frac12\,.
\end{align}
As for the term $J_{1,2}$, we compute as 
\begin{align*}
J_{1,2} [\varphi] & = - \mu_n \int_{IV_{n,1}} \frac{\cos y}{\mu_n-\sin y} \frac{v_n \varphi}{\cos y} \, \dd y - \int_{IV_{n,1}} \frac{\delta_n f_n}{\mu_n-\sin y} \varphi \, \dd y \\
& = - \mu_n \int_{IV_{n,1}} \frac{\cos y}{\mu_n-\sin y} \big ( \frac{v_n \varphi}{\cos y} - \frac{v_n(y_{n,1}) \varphi (y_{n,1})}{\cos y_{n,1}}\big ) \, \dd y \\
& \quad - \frac{\mu_n v_n(y_{n,1}) \varphi (y_{n,1})}{\cos y_{n,1}} \int_{IV_{n,1}} \frac{\cos y}{\mu_n-\sin y}   \, \dd y- \int_{IV_{n,1}} \frac{\delta_n f_n}{\mu_n-\sin y} \varphi \, \dd y\\
& =: J_{1,2,1} [\varphi] + J_{1,2,2}[\varphi] + J_{1,2,3}[\varphi]\,.
\end{align*}
Then we have from $|\cos y|\geq \frac{1}{C} |1-\mu_n|^\frac12$ for $y\in III_{n,1}\cup IV_{n,1}$,
\begin{align}\label{proof.prop.Lambda.Kol.rate.2.15}
|J_{1,2,1}[\varphi]| & \leq C \kappa_3^\frac12 |1-\mu_n|^\frac14 \| \partial_y ( \frac{v_n \varphi}{\cos y} )\|_{L^2(III_{n,1}\cup IV_{n,1})} \nonumber \\
& \leq C\kappa_3^\frac12 |1-\mu_n|^\frac14 \big ( \frac{\| \partial_y v_n \|_{L^2} \| \varphi \|_{L^\infty} + \| v_n \|_{L^\infty} \| \partial_y \varphi \|_{L^2}}{|1-\mu_n|^\frac12} + \frac{\| v_n \|_{L^\infty} \| \varphi\|_{L^\infty}}{|1-\mu_n|^\frac34} \big ) \nonumber \\
& \leq C\kappa_3^\frac12 \big ( \frac{\|\varphi\|_{L^\infty}}{|1-\mu_n|^\frac14} + \| \partial_y \varphi \|_{L^2} \big ) \quad ({\rm by ~\eqref{proof.prop.Lambda.Kol.rate.2.2'} ~and~ \eqref{proof.prop.Lambda.Kol.rate.2.2''}})\,.
\end{align}
Hence we have 
\begin{align}\label{proof.prop.Lambda.Kol.rate.2.16}
|J_{1,2,1}[v_n]| & \leq C \kappa_3^\frac12\,.
\end{align}
The term $J_{1,2,2}$ is estimated from \eqref{proof.prop.Lambda.Kol.rate.2.11},
\begin{align}\label{proof.prop.Lambda.Kol.rate.2.17}
|J_{1,2,2}[\varphi]| \leq C \frac{|v_n (y_{n,1}) \varphi (y_{n,1})|}{|1-\mu_n|^\frac12}\,,
\end{align}
and thus, 
\begin{align}\label{proof.prop.Lambda.Kol.rate.2.18}
|J_{1,2,2}[v_n]| \leq C \frac{|v_n (y_{n,1})|^2}{|1-\mu_n|^\frac12}\,.
\end{align} 
As for $J_{1,2,3}$, we have from the definition of $IV_{n,1}$ and \eqref{est.sin},
\begin{align}\label{proof.prop.Lambda.Kol.rate.2.19}
|J_{1,2,3}[\varphi]| \leq \| \frac{1}{\mu_n - \sin y} \|_{L^2 (IV_{n,1})} \| \delta_n f_n \|_{L^2} \| \varphi \|_{L^\infty} 
& \leq \frac{C\delta_n \|f_n \|_{L^2} \| \varphi \|_{L^\infty} }{|1-\mu_n|^\frac12 \cdot \kappa_3^\frac12 \frac{\delta_n}{|1-\mu_n|^\frac14}}  \leq \frac{C \|f_n\|_{L^2} \| \varphi\|_{L^\infty}}{\kappa_3^\frac12 |1-\mu_n|^\frac14}\,.
\end{align}  
Thus it follows that
\begin{align}\label{proof.prop.Lambda.Kol.rate.2.20}
|J_{1,2,3}[v_n]| \leq \frac{C\|f_n \|_{L^2}}{\kappa_3^\frac12}\,.
\end{align}
Next we estimate $J_{1,3}$:
\begin{align*}
J_{1,3}[\varphi] & = - \int_{[-\frac{\pi}{2},\frac{\pi}{2}]\setminus (III_{n,1} \cup IV_{n,1})} \frac{\mu_n v_n \, \varphi}{\mu_n-\sin y}\, \dd y - \int_{[-\frac{\pi}{2},\frac{\pi}{2}]\setminus (III_{n,1} \cup IV_{n,1})} \frac{\delta_n f_n \, \varphi}{\mu_n-\sin y} \, \dd y\\
& =: J_{1,3,1}[\varphi] + J_{1,3,2}[\varphi]\,.
\end{align*}
Since $\mu_n-\sin y$ is positive for $y\in [-\frac{\pi}{2},y_{n,1})$, we have from the definition of $III_{n,1}$ and $IV_{n,1}$,
\begin{align}\label{proof.prop.Lambda.Kol.rate.2.21}
J_{1,3,1} [v_n] & = - \int_{-\frac{\pi}{2}}^{y_{n,1}-\kappa_3 |1-\mu_n|^\frac12} \frac{\mu_n v_n^2}{\mu_n-\sin y} \, \dd y - \int_{y_{n,1}+ \kappa_3 |1-\mu_n|^\frac12}^\frac{\pi}{2} \frac{\mu_n v_n^2}{\mu_n-\sin y} \, \dd y \nonumber \\
& \leq - \int_{y_{n,1}+ \kappa_3 |1-\mu_n|^\frac12}^\frac{\pi}{2} \frac{\mu_n v_n^2}{\mu_n-\sin y} \, \dd y  \nonumber \\
& \leq \frac{C\|v_n\|_{L^2([y_{n,1}, y_{n,2}])}^2}{\kappa_3(1-\mu_n)}\,.
\end{align}
On the other hand, for general $\varphi\in H^2(\T)$ we have 
\begin{align}\label{proof.prop.Lambda.Kol.rate.2.22}
|J_{1,3,1}[\varphi]| & \leq |\int_{[-\frac{\pi}{2},\frac{\pi}{2}]\setminus (III_{n,1} \cup IV_{n,1})} \frac{\mu_n (v_n -  v_n (y_{n,1})) \, \varphi}{\mu_n-\sin y}\, \dd y | + |\mu_n v_n (y_{n,1}) \int_{[-\frac{\pi}{2},\frac{\pi}{2}]\setminus (III_{n,1} \cup IV_{n,1})} \frac{\varphi}{\mu_n-\sin y}\, \dd y | \nonumber \\
& \leq \frac{C\| \partial_y v_n \|_{L^2} \| \varphi\|_{L^2}}{|1-\mu_n|^\frac12} + \frac{C |\mu_n v_n (y_{n,1})| \, \| \varphi \|_{L^2}}{\kappa_3^\frac12 |1-\mu_n|^\frac34}\,.
\end{align}
The estimate of $J_{1,3,2}$ is similar to $J_{1,2,2}$, and we have 
\begin{align}\label{proof.prop.Lambda.Kol.rate.2.23}
|J_{1,3,2}[\varphi]|  \leq \| \frac{1}{\mu_n-\sin y} \|_{L^2 ([-\frac{\pi}{2},\frac{\pi}{2}]\setminus (III_{n,1}\cup IV_{n,1}))} \| \delta_n f_n \|_{L^2} \| \varphi \|_{L^\infty} & \leq \frac{C\delta_n \| f_n \|_{L^2} \|\varphi \|_{L^\infty} }{|1-\mu_n|^\frac12 \cdot \kappa_3^\frac12 |1-\mu_n|^\frac14} \nonumber \\
& \leq \frac{C\|f_n\|_{L^2}\| \varphi\|_{L^\infty}}{\kappa_3^\frac12 |1-\mu_n|^\frac14}\,.
\end{align}
Hence, 
\begin{align}\label{proof.prop.Lambda.Kol.rate.2.24}
|J_{1,3,2}[v_n]|\leq \frac{C\|f_n\|_{L^2}}{\kappa_3^\frac14}\,.
\end{align}
Collecting \eqref{proof.prop.Lambda.Kol.rate.2.14}, \eqref{proof.prop.Lambda.Kol.rate.2.16}, \eqref{proof.prop.Lambda.Kol.rate.2.18}, \eqref{proof.prop.Lambda.Kol.rate.2.20}, \eqref{proof.prop.Lambda.Kol.rate.2.21}, \eqref{proof.prop.Lambda.Kol.rate.2.24}, we have 
\begin{align}\label{proof.prop.Lambda.Kol.rate.2.25}
J_1[v_n]\leq \frac{C \|v_n \|_{L^2 ([y_{n,1},y_{n,2}])}^2}{\kappa_3 (1-\mu_n)} + C\frac{|v_n (y_{n,1})|^2}{|1-\mu_n|^\frac12} + \frac{C\|f_n\|_{L^2}}{\kappa_3^\frac12} + C \kappa_3^\frac12\,.
\end{align}
The estimate of $J_2$ is exactly the same as $J_1$. Hence we have 
\begin{align}\label{proof.prop.Lambda.Kol.rate.2.26}
\| \partial_y v_n\|_{L^2}^2 + (l_n^2-1) \| v_n \|_{L^2}^2  & = J_1[v_n] +J_2[v_n]  \nonumber \\
& \leq  \frac{C \|v_n \|_{L^2 ([y_{n,1},y_{n,2}])}^2}{\kappa_3 (1-\mu_n)}  + C \sum_{j=1,2} \frac{|v_n (y_{n,j})|^2}{|1-\mu_n|^\frac12} + \frac{C\|f_n\|_{L^2}}{\kappa_3^\frac12} + C \kappa_3^\frac12\,.
\end{align}
By Step 2 and $\displaystyle \lim_{n\rightarrow \infty} \| f_n \|_{L^2}=0$, we conclude that 
\begin{align}\label{proof.prop.Lambda.Kol.rate.2.27}
\limsup_{n\rightarrow \infty} \, \big ( \| \partial_y v_n \|_{L^2}^2 + (l_n^2-1) \| v_n \|_{L^2}^2 \big ) \leq \frac{C}{\kappa_3} \limsup_{n\rightarrow \infty} \frac{\|v_n \|_{L^2 ([y_{n,1},y_{n,2}])}^2}{1-\mu_n} + C\kappa_3^\frac12  \,,
\end{align}
for any small $\kappa_3>0$. 
Suppose that $\frac12 \leq \mu_\infty<1$ and $l^2_\infty\in [1,\infty)$. 
In this case we have from \eqref{proof.prop.Lambda.Kol.rate.2.13}, \eqref{proof.prop.Lambda.Kol.rate.2.15}, \eqref{proof.prop.Lambda.Kol.rate.2.17}, \eqref{proof.prop.Lambda.Kol.rate.2.19}, \eqref{proof.prop.Lambda.Kol.rate.2.22}, \eqref{proof.prop.Lambda.Kol.rate.2.23},
\begin{align}\label{proof.prop.Lambda.Kol.rate.2.28}
\big | \langle \partial_y v_\infty, \partial_y \varphi\rangle_{L^2} \big | = \big | \lim_{n\rightarrow \infty} \langle \partial_y v_n,\partial_y \varphi \rangle_{L^2} \big | 
& \leq (l_\infty^2-1)  \big |\langle v_\infty,\varphi \rangle_{L^2} \big | + \limsup_{n\rightarrow \infty} \big | \langle \frac{\mu_n v_n + \delta_n f_n}{\mu_n-\sin y}, \varphi \rangle_{L^2}\big | \nonumber \\
& \leq (l_\infty^2-1) \| v_\infty\|_{L^2} \| \varphi\|_{L^2} + C \| \varphi \|_{L^2}\,,
\end{align}
for any $\varphi \in H^2(\T)$.

\

\noindent {\bf Case $\frac12\leq \mu_\infty<1$.}

\noindent Let us first consider the case $\mu_\infty<1$. If $l_\infty=\infty$ then \eqref{proof.prop.Lambda.Kol.rate.2.27} implies that $\displaystyle \lim_{n\rightarrow \infty} \|v_n\|_{L^2}^2=0$, and then, again by \eqref{proof.prop.Lambda.Kol.rate.2.27} and $\mu_\infty<1$, we have $\displaystyle \lim_{n\rightarrow \infty} \|\partial_y v_n\|_{L^2}^2=0$, since $\kappa_3>0$ is arbitrary small. By recalling that $\displaystyle \lim_{n\rightarrow \infty} \delta_n^2\|u_n\|_{L^2}^2=0$ by Step 1, we achieve the contradiction with the normalized condition \eqref{proof.prop.Lambda.Kol.rate.2.2'}.
Hence we conclude that $l_\infty<\infty$. In this case, in virtue of \eqref{proof.prop.Lambda.Kol.rate.2.27} and the assumption $\mu_\infty<1$, we may assume that $\inf_n \| v_n \|_{L^2}>0$. Since $v_n$ converges to $v_\infty$ strongly in $L^2(\T)$, the limit $v_\infty$ is a  nontrivial function. Moreover, we see from \eqref{proof.prop.Lambda.Kol.rate.2.28} that $v_\infty$ belongs to $H^2(\T)$.  The direct computation using the weak formulation implies that $v_\infty$ satisfies $(\mu_\infty- \sin y) (A_{l_\infty}+1) v_\infty -\mu_\infty v_\infty=0$ for $y\in \T\setminus  S_{\mu_\infty}$. Since we have already shown that $v_\infty\in H^2(\T)$, we conclude that $u_\infty = A_{l_\infty} v_\infty$ is an eigenfunction of $\hat{\Lambda}_{l_\infty}$ in $L^2(\T)$. By Proposition \ref{prop.kernel.Kol} this is a contradiction since $\mu_\infty\in [\frac12,1)$.  The case $\frac12 \leq \mu_\infty <1$ is settled.

\

\noindent {\bf Case $\mu_\infty=1$.}

\noindent  In this case we need additional steps to achieve the contradiction.

\noindent {\bf Step 5: $\displaystyle \limsup_{n\rightarrow \infty} \frac{\|v_n \|_{L^2 ([y_{n,1},y_{n,2}])}^2}{1-\mu_n}>0$ when $\mu_\infty=1$.}

\noindent Suppose that $\mu_\infty=1$ and $\displaystyle \lim_{n\rightarrow \infty} \frac{\|v_n \|_{L^2 ([y_{n,1},y_{n,2}])}^2}{1-\mu_n}=0$.  Then \eqref{proof.prop.Lambda.Kol.rate.2.9} implies $\displaystyle \lim_{n\rightarrow \infty} \frac{\|v_n\|_{L^2}^2}{1-\mu_n}=0$, and \eqref{proof.prop.Lambda.Kol.rate.2.27} implies $\displaystyle \lim_{n\rightarrow \infty} \| \partial_y v_n \|_{L^2}^2=0$ as well. Since $\displaystyle \lim_{n\rightarrow \infty} \delta_n^2 \| u_n\|_{L^2}^2=0$ by Step 1, we achieve the contradiction due to the normalized condition \eqref{proof.prop.Lambda.Kol.rate.2.2'}.

\

\noindent {\bf Step 6: Rescaling and limiting process for the case $\mu_\infty=1$.}

\noindent By Step 5 we may assume that $\displaystyle \inf_n \frac{\| v_n\|_{L^2([y_{n,1},y_{n,2}])}^2}{1-\mu_n}>0$ (by taking suitable subsequence if necessary). Set 
\begin{align}\label{proof.prop.Lambda.Kol.rate.2.29}
w_n (\xi) = \frac{1}{|1-\mu_n|^{\frac14}} v_n (\frac{\pi}{2} + |1-\mu_n|^\frac12 \xi)\,, \qquad \xi \in [-2, 2]\,.
\end{align}
Let $c_n>0$ be the number such that 
\begin{align}\label{proof.prop.Lambda.Kol.rate.2.30}
y_{n,1} = \frac{\pi}{2} - |1-\mu_n|^\frac12 c_n\,, \qquad y_{n,2} = \frac{\pi}{2} + |1-\mu_n|^\frac12 c_n\,.
\end{align}
Note that
\begin{align*}%\label{proof.prop.Lambda.Kol.rate.2.39}
\mu_n = \sin y_{n,j} = \sin \frac{\pi}{2} -\frac12 \,  (\frac{\pi}{2}-y_{n,j})^2 + \frac{1}{4!}  (\frac{\pi}{2}-y_{n,j})^4 + \cdots\,,
\end{align*}
which gives
\begin{align*}%\label{proof.prop.Lambda.Kol.rate.2.40}
(y_{n,j}-\frac{\pi}{2})^2 = 2 (1-\mu_n) + O (|1-\mu_n|^2)\,.
\end{align*}
Hence we see from $(y_{n,j}-\frac{\pi}{2})^2 = (1-\mu_n ) c_n^2$ by its definition,
\begin{align}\label{proof.prop.Lambda.Kol.rate.2.31}
c_n = \sqrt{2} + O(|1-\mu_n|^\frac12)\,.
\end{align} 
In particular, we have 
\begin{align}\label{proof.prop.Lambda.Kol.rate.2.32}
1\leq c_n \leq 2 \qquad \text{ for all large }~n\,, \qquad \lim_{n\rightarrow \infty} c_n = \sqrt{2}\,.
\end{align}
In virtue of the normalized condition \eqref{proof.prop.Lambda.Kol.rate.2.2'} and Step 5, we have 
\begin{align}\label{proof.prop.Lambda.Kol.rate.2.33}
\begin{split}
\|\partial_\xi w_n \|_{L^2(-2,2))}^2 + \| w_n \|_{L^2((-2,2))}^2 & \leq  \|\partial_y v_n \|_{L^2}^2 + \frac{\| v_n \|_{L^2}^2}{|1-\mu_n|} \leq 1\,, \\
\inf_n \|w_n \|_{L^2((-c_n,c_n))}^2 & = \inf_n \frac{\| v_n \|_{L^2((y_{n,1},y_{n,2}))}^2}{1-\mu_n} >0\,.
\end{split}
\end{align} 
That is, the sequence $\{w_n\}$ is uniformly bounded in $H^1((-2, 2))$, and thus, we may assume that $w_n$ weakly converges in $H^1((-2,2))$ to some $w_\infty\in H^1 ((-2,2))$ and strongly converges in $L^2((-2,2))$ as well as in $C^\eta ([-2,2])$ for some $\eta>0$. Moreover, by the uniform lower bound in \eqref{proof.prop.Lambda.Kol.rate.2.33}, the limit $w_\infty$ is nontrivial.
The direct computation shows that $w_n$ satisfies 
\begin{align}\label{proof.prop.Lambda.Kol.rate.2.34}
\begin{split}
(\mu_n - \sin_n \xi ) \partial_{\xi}^2 w_n & = |1-\mu_n| \, \bigg ( (\mu_n - \sin_n \xi ) (l_n^2-1) w_n + \mu_n w_n  + \delta_n g_n\bigg )\,, \\
&\qquad \xi\in (-2,2)\,.
\end{split}
\end{align}
Here we have set 
\begin{align*}
\sin_n \xi = \sin (\frac{\pi}{2}+ |1-\mu_n|^\frac12 \xi)\,, \qquad g_n (\xi) =|1-\mu_n|^{-\frac14} f_n (\frac{\pi}{2}+|1-\mu_n|^\frac12\xi)\,.
\end{align*}
On the points $\xi=\pm c_n$ we have 
\begin{align}\label{proof.prop.Lambda.Kol.rate.2.35}
w_n (-c_n) = |1-\mu_n|^{-\frac14} v_n (y_{n,1})\,, \qquad w_n (c_n) = |1-\mu_n|^{-\frac14} v_n (y_{n,2})\,.
\end{align}
Thus  Step 2 gives
\begin{align}\label{proof.prop.Lambda.Kol.rate.2.36}
w_\infty (\pm \sqrt{2}) =0\,.
\end{align}
Next we see 
\begin{align*}
\sin_n \xi = \sin \frac{\pi}{2} - \frac12 |1-\mu_n| \xi^2  + \frac{1}{4!} |1-\mu_n|^2 \xi^4 \cdots
\end{align*}
and 
\begin{align*}
\mu_n = \sin y_{n,2} = \sin (\frac{\pi}{2} + |1-\mu_n|^\frac12 c_n) = 1  -\frac12 |1-\mu_n| c_n^2 + \frac{1}{4!} |1-\mu_n|^2 c_n^4 + \cdots
\end{align*}  
which gives
\begin{align*}
\mu_n  -\sin_n\xi = \frac12 |1-\mu_n| (\xi^2 -c_n^2) \, \bigg ( 1 + |1-\mu_n| \,  q_n (\xi) \bigg )\,.
\end{align*}
Here $q_n$ is a  smooth function on $[-2,2]$ satisfying the uniform bound 
$$\sup_n \| \frac{\dd^k q_n}{\dd \xi^k}  \|_{L^\infty ((-2,2))} <\infty\quad  \text{for}~~ k=0,1,2\,. $$
Thus \eqref{proof.prop.Lambda.Kol.rate.2.34} is written as 
\begin{align}\label{proof.prop.Lambda.Kol.rate.2.37}
\begin{split}
& \frac12  (\xi^2-c_n^2  )\bigg (1+|1-\mu_n|\, q_n \bigg) \partial_\xi^2w_n \\
&\qquad = \frac12|1-\mu_n| (l_n^2-1) (\xi^2-c_n^2)\bigg (1+|1-\mu_n| \, q_n \bigg)  w_n + \mu_n w_n + \delta_n g_n\,.
\end{split}
\end{align}
Note that, from $1-\mu_n\geq \kappa \delta_n^2$ and $\|f_n \|_{L^2}\rightarrow 0$ as $n\rightarrow \infty$, 
\begin{align}\label{proof.prop.Lambda.Kol.rate.2.38}
\delta_n \| g_n\|_{L^2((-2, 2))} \leq \frac{C \delta_n}{|1-\mu_n|^\frac12} \| f_n \|_{L^2} \rightarrow 0\qquad n \rightarrow \infty\,.
\end{align}
We have from \eqref{proof.prop.Lambda.Kol.rate.2.37} that for any test function $\varphi\in C_0^\infty ((-2,2))$,
\begin{align}\label{proof.prop.Lambda.Kol.rate.2.39}
\begin{split}
& -\langle \partial_\xi w_n, \partial_\xi \varphi \rangle_{L^2((-2,2))} \\
&\quad = |1-\mu_n|(l_n^2-1) \, \langle w_n, \varphi\rangle_{L^2((-2,2))} + 2 \langle \frac{\mu_n w_n + \delta_n g_n}{\xi^2-c_n^2}, \frac{\varphi}{1+|1-\mu_n|\, q_n} \rangle_{L^2((-2,2))} \,.
\end{split}
\end{align}
Let us focus on the second term of the right-hand side of \eqref{proof.prop.Lambda.Kol.rate.2.39}.
\begin{align*}
\langle \frac{\mu_n w_n + \delta_n g_n}{\xi^2-c_n^2}, \frac{\varphi}{1+|1-\mu_n|\, q_n} \rangle_{L^2((-2,2))} & = \int_{[0,2]} \frac{1}{\xi^2-c_n^2} \Big[h_n \Big]_{(+)}  \dd y  + \int_{[-2,0]} \frac{1}{\xi^2-c_n^2} \Big[h_n \Big]_{(-)}  \dd y   \,.
\end{align*}
Here $$h_n = \frac{\Big (\mu_n w_n + \delta_n g_n\Big ) \varphi}{1+|1-\mu_n| q_n}\,, \qquad [q]_{(\pm)} (\xi) := q (\xi) - q (\pm c_n)\,,$$
where we have used $(\mu_n w_n+\delta_n g_n)(\pm c_n)=0$.
Let us decompose $h_n$ as 
\begin{align*}
& h_n = h_{n,1} + h_{n,2}\,, \qquad h_{n,1} = \frac{\mu_n w_n \varphi}{1+|1-\mu_n| q_n}\,, \qquad  h_{n,2} = \frac{\delta_n g_n \varphi}{1+|1-\mu_n| q_n}\,.
\end{align*}
Then, by the weak convergence in $H^1$ for $w_n$ and strong convergence to zero in $H^2$ for $|1-\mu_n|q_n$, it is not difficult to show the convergence 
\begin{align}\label{proof.prop.Lambda.Kol.rate.2.40}
& \int_{[0,2]} \frac{1}{\xi^2-c_n^2} \Big[h_{n,1} \Big]_{(+)}  \dd y + \int_{[-2,0]} \frac{1}{\xi^2-c_n^2} \Big[h_{n,1} \Big]_{(-)}  \dd \xi \nonumber \\
& \rightarrow \int_{[0,2]}\frac{1}{\xi^2-2} \Big[w_\infty \varphi \Big]_{(+)}  \dd y + \int_{[-2,0]} \frac{1}{\xi^2-2} \Big[w_\infty \varphi \Big]_{(-)}  \dd \xi   = \int_{-2}^2 \frac{w_\infty \varphi}{\xi^2-2} \, \dd \xi\,.
\end{align}
Here $[w_\infty \varphi]_{(\pm)} (\xi) := (w_\infty \varphi) (\xi) - (w_\infty \varphi) (\pm \sqrt{2}) = (w_\infty \varphi) (\xi)$ is used in the last line since $w_\infty (\pm\sqrt{2})=0$. 
Next we show 
\begin{align}\label{proof.prop.Lambda.Kol.rate.2.41}
\big | \int_{[0,2]} \frac{1}{\xi^2-c_n^2} \Big[h_{n,2} \Big]_{(+)}  \dd y \big |+ \big | \int_{[-2,0]} \frac{1}{\xi^2-c_n^2} \Big[h_{n,2} \Big]_{(-)}  \dd y \big | \rightarrow 0\,.
\end{align}
It suffices to consider the integral over $[-2,0]$. Recall from $\displaystyle \lim_{n\rightarrow \infty} w_n (-c_n)=0$ and $(\mu_n w_n + \delta_n g_n )(\pm c_n)=0$ that we have $\displaystyle \lim_{n\rightarrow \infty} \delta_n g_n (-c_n) =0$. Then we have 
\begin{align*}
\int_{[-2,0]} \frac{1}{\xi^2-c_n^2} \Big[h_{n,2} \Big]_{(-)}  \dd \xi & =  \int_{[-2,0]}  \frac{1}{\xi^2-c_n^2} \Big[\delta_n g_n  \Big]_{(-)} \,  \varphi_n \, \dd \xi + o(1)\,, \qquad \varphi_n : = \frac{\varphi}{1+|1-\mu_n| q_n} \,.
\end{align*}
We take arbitrary small $\kappa''>0$ and set $\tilde \Pi_n =[-c_n- \kappa'' \frac{\delta_n^2}{1-\mu_n}, -c_n + \kappa'' \frac{\delta_n^2}{1-\mu_n}]$. Then we have 
\begin{align*}
\big | \int_{\tilde \Pi_n} \frac{1}{\xi^2-c_n^2} \Big[\delta_n g_n  \Big]_{(-)}  \varphi_n \,  \dd \xi \big | & \leq C (\kappa'')^\frac12 \frac{\delta_n}{|1-\mu_n|^\frac12} \| \delta_n \partial_\xi g_n \|_{L^2 (\tilde \Pi_n)} \| \varphi \|_{L^\infty} \\
& \leq C  (\kappa'')^\frac12 \frac{\delta_n}{|1-\mu_n|^\frac12} \| \delta_n \partial_y f_n \|_{L^2 (\Pi_n)} \| \varphi \|_{L^\infty} \,.
\end{align*}
Here $\Pi_n = [y_{n,1}-\kappa'' \frac{\delta_n^2}{|1-\mu_n|^\frac12}, y_{n,1} + \kappa'' \frac{\delta_n^2}{|1-\mu_n|^\frac12}]$, and the norm $\| \delta_n \partial_y f_n \|_{L^2 (\Pi_n)}$ is already estimated as in \eqref{proof.prop.Lambda.Kol.rate.2.12}. Thus we have 
\begin{align*}
\limsup_{n\rightarrow \infty} \big | \int_{\tilde \Pi_n} \frac{1}{\xi^2-c_n^2} \Big[\delta_n g_n  \Big]_{(-)} \frac{\varphi}{1+|1-\mu_n| q_n} \,  \dd \xi \big | \leq C (\kappa'')^\frac12 \,.
\end{align*}
On the other hand, we have 
\begin{align*}
\big | \int_{[-2,0]\setminus \tilde \Pi_n} \frac{1}{\xi^2-c_n^2} \Big[\delta_n g_n  \Big]_{(-)} \varphi_n  \dd \xi \big | & \leq \big | \int_{[-2,0]\setminus \tilde \Pi_n} \frac{1}{\xi^2-c_n^2} \delta_n g_n  \varphi_n \dd \xi \big | \nonumber \\
& \quad + | \delta_n g_n (-c_n) \int_{[-2,0]\setminus \tilde \Pi_n} \frac{\varphi_n-\varphi_n (-c_n)}{\xi^2-c_n^2}  \dd \xi \big | \\
& \quad + | \delta_n g_n (-c_n) \varphi_n (-c_ n) \int_{[-2,0]\setminus \tilde \Pi_n} \frac{1}{\xi^2-c_n^2}  \dd \xi \big |  \nonumber \\
& \leq \frac{C|1-\mu_n|^\frac12}{(\kappa'')^\frac12 \delta_n} \| \delta_n g_n \|_{L^2 ((-2,2))} \| \varphi_n \|_{L^\infty} +  C \| \varphi_n \|_{H^1} |\delta_n g_n (-c_n)| \nonumber\\
& \leq \frac{C}{(\kappa'')^\frac12} \| f_n \|_{L^2} + C |\delta_n g_n (-c_n)|  \quad \rightarrow 0 \quad (n\rightarrow \infty)\,.
\end{align*}
Here \eqref{proof.prop.Lambda.Kol.rate.2.38} and $\displaystyle \sup_n \big |\int_{[-2,0]\setminus \tilde \Pi_n} \frac{1}{\xi^2-c_n^2} \, \dd \xi \big | <\infty$ are used. Hence \eqref{proof.prop.Lambda.Kol.rate.2.41} holds, by taking $\kappa''\rightarrow 0$ after $n\rightarrow \infty$.
Collecting these, we conclude that 
\begin{align}\label{proof.prop.Lambda.Kol.rate.2.42}
\begin{split}
2 \lim_{n\rightarrow \infty} \langle \frac{\mu_n w_n + \delta_n g_n}{\xi^2-c_n^2}, \frac{\varphi}{1+|1-\mu_n|\, q_n} \rangle_{L^2((-2,2))}  = 2 \int_{-2}^2 \frac{w_\infty \varphi }{\xi^2-2} \dd \xi \,, \qquad w_\infty (\pm \sqrt{2})=0 \,.
\end{split}
\end{align}

\

\noindent {(i) When $\displaystyle \limsup_{n\rightarrow \infty}|1-\mu_n| (l_n^2-1)=\infty$:}

\noindent In this case we may assume that $\displaystyle \lim_{n\rightarrow \infty}|1-\mu_n| (l_n^2-1)=\infty$ by taking a suitable subsequence. Then we divide both sides of \eqref{proof.prop.Lambda.Kol.rate.2.39} by $|1-\mu_n| (l_n^2-1)$ and consider the weak formulation with arbitrary test function $\varphi\in C_0^\infty ((-2,2))$:
\begin{align*}
& \frac{-1}{|1-\mu_n| (l_n^2-1)}\langle \partial_\xi w_n, \partial_\xi \varphi \rangle_{L^2 ((-2,2))} \\
& = \langle  w_n , \varphi \rangle_{L^2 ((-2,2))} + \frac{2}{|1-\mu_n| (l_n^2-1)}  \langle \frac{\mu_n w_n + \delta_n g_n}{\xi^2-c_n^2}, \frac{\varphi}{1+|1-\mu_n|\, q_n} \rangle_{L^2((-2,2))}\,. 
\end{align*}  
Then, by taking the limit $n\rightarrow \infty$ in the above weak formulation and using \eqref{proof.prop.Lambda.Kol.rate.2.42},
we obtain $w_\infty=0$, which is a contradiction.   

\

\noindent  {(ii) When $\displaystyle \limsup_{n\rightarrow \infty}|1-\mu_n| (l_n^2-1)<\infty$:}

\noindent In this case we may assume that $\displaystyle \lim_{n\rightarrow \infty}|1-\mu_n| (l_n^2-1)=d_\infty\in [0,\infty)$ by taking a suitable subsequence. 
Then, by considering the weak formulation for \eqref{proof.prop.Lambda.Kol.rate.2.39} with arbitrary test function $\phi\in C_0^\infty ((-2, 2))$ as above, we verify from \eqref{proof.prop.Lambda.Kol.rate.2.42} that the limit  $w_\infty\in H^1 ((-2,2))$ satisfies $w_\infty (\pm \sqrt{2})=0$ and
\begin{align}\label{proof.prop.Lambda.Kol.rate.2.493}
- \langle \partial_\xi w_\infty, \partial_\xi \varphi \rangle_{L^2 ((-2,2))} = d_\infty \langle w_\infty, \varphi\rangle_{L^2 ((-2,2))} + 2 \int_{-2}^2 \frac{w_\infty \varphi }{\xi^2-2} \dd \xi \,.
\end{align}
By the Hardy inequality the second term in the right-hand side of \eqref{proof.prop.Lambda.Kol.rate.2.493} is bounded from above by $C \|w_\infty \|_{H^1} \| \varphi \|_{L^2}$, which implies $w_\infty \in H^2 ((-2,2))$. 
In particular, $w_\infty\in C^{1+\eta} ((-2,2))$ for some $\eta>0$.
By considering the test function of the form $(\xi^2-2)\phi$  with $\phi\in C_0^\infty ((-2,2))$, we also have from \eqref{proof.prop.Lambda.Kol.rate.2.493} that 
\begin{align}\label{proof.prop.Lambda.Kol.rate.2.49}
(\xi^2-2) \partial_\xi^2 w_\infty =d_\infty  (\xi^2-2)  w_\infty + 2 w_\infty\,,  \qquad \xi \in (-2,2)\setminus \{\pm \sqrt{2}\}\,, \qquad w_\infty (\pm \sqrt{2})=0\,.
\end{align}
We can show from \eqref{proof.prop.Lambda.Kol.rate.2.49} that $w_\infty\in H^2 ((-2,2))$ is smooth except for the points $\xi=\pm \sqrt{2}$,
and thus \eqref{proof.prop.Lambda.Kol.rate.2.49} is satisfied pointwise in $(-2,2)\setminus\{\pm\sqrt{2}\}$. Our aim is to show that $w_\infty=0$ in $(-2,2)$, for it leads to the contradiction.
The difficulty is that the polynomial $\xi^2-2$ satisfies \eqref{proof.prop.Lambda.Kol.rate.2.49} at least when $d_\infty=0$, therefore we need to derive some additional estimate for $w_\infty$.
The key is the estimate of $v_n$ outside the interval $(y_{n,1},y_{n,2})$, which is obtained in Step 3.
Indeed, \eqref{proof.prop.Lambda.Kol.rate.2.29} implies that
\begin{align*}
\int_{-2}^{-c_n} |w_n |^2 \dd \xi  & = |1-\mu_n|^{-\frac12} \int_{-2}^{-c_n} |v_n (\frac{\pi}{2} + |1-\mu_n|^\frac12 \xi )|^2 \dd \xi  = \frac{1}{1-\mu_n}\int_{\frac{\pi}{2}-2|1-\mu_n|^\frac12}^{y_{n,1}} v_n^2 \dd y \,,
\end{align*}
and hence, since $0<\frac{\pi}{2}-2|1-\mu_n|^\frac12<y_{n,1}$, estimate \eqref{proof.prop.Lambda.Kol.rate.2.9} leads to 
\begin{align*}
\lim_{n\rightarrow \infty} \int_{-2}^{-c_n} | w_n |^2 \dd \xi  =0 \,.
\end{align*}
That is, $w_\infty=0$ for $\xi\in (-2,-\sqrt{2})$. Since $w_\infty\in C^{1+\eta} ((-2,2))$ we conclude that 
\begin{align}\label{proof.prop.Lambda.Kol.rate.2.51}
\partial_\xi w_\infty (-\sqrt{2})=0\,.
\end{align}
Note that the singularity $(\xi^2-2)^{-1}$ at $\xi=-\sqrt{2}$ is first order.
Then it is easy to show that the solution $w_\infty\in H^2((-2,2))$ to \eqref{proof.prop.Lambda.Kol.rate.2.49} satisfying the {\it initial} condition $w_\infty (-\sqrt{2})=\partial_\xi w_\infty (-\sqrt{2})=0$ must be trivial, i.e., $w_\infty=0$ in $(-2,2)$.
This contradicts with $\|w_\infty\|_{L^2((-2,2))}>0$. The proof of \eqref{est.lem.Lambda.Kol.rate.2} is complete.

%\begin{lem}\label{lem.Lambda.Kol.rate.2-2}  The inequality \eqref{est.lem.Lambda.Kol.rate.2'} holds.
% \end{lem}

\
 
\noindent {\it Proof of \eqref{est.lem.Lambda.Kol.rate.2'}.} Again we will use the contradiction argument.
Suppose that 
there exist $\{\delta_n, l_n, \mu_n\}_{n\in \N}$, $\delta_n\in (0,1]$, $l_n\in \Z\setminus\{0\}$, $\mu_n\in (-1+\kappa\delta_n^2,-\frac12]\cup [\frac12, 1-\kappa \delta_n^2)$, and $\{u_n\}\subset H^2 (\T; \R)$ such that 
\begin{align*}
& \lim_{n\rightarrow \infty} \delta_n=\delta_\infty\in [0,1]\,, \quad \lim_{n\rightarrow \infty} l_n = l_\infty\in \{\pm \infty\}\cup \Z\setminus\{0\}\,,\\
& \lim_{n\rightarrow \infty}\mu_n=\mu_\infty \in [-1+\kappa\delta^2_\infty,-\frac12]\cup [\frac12, 1-\kappa \delta_\infty^2]\,,
\end{align*}
and 
\begin{align}\label{proof.lem.Lambda.Kol.rate.2-2.-1}
\begin{split}
& \| M_{\cos y}  B_{2,l_n} u_n  \|_{L^2}^2 =1\,, \\
& \lim_{n\rightarrow \infty}  \bigg ( \delta_n^{-2}  \|  (\mu_n - \hat{\Lambda}_{l_n}) u_n \|_{L^2}^2 + \delta_n^2 (1-|\mu_n|) \| (-A_{l_n})^\frac12  u_n\|_{L^2}^2 \bigg ) =0\,.
\end{split}
\end{align}
As in the previous lemma, set
\begin{align}\label{proof.lem.Lambda.Kol.rate.2-2.-2}
f_n = \delta_n^{-1} (\mu_n   - \hat{\Lambda}_{l_n} ) u_n\,, \qquad v_n = A_{l_n}^{-1} u_n\,.
\end{align}
Since $u_n$ is real valued, so is $v_n$, and $v_n$ satisfies 
\begin{align}\label{proof.lem.Lambda.Kol.rate.2-2.-3}
(\mu_n - M_{\sin y}) \big ( A_{l_n} + 1) v_n - \mu_n  v_n = \delta_n f_n\,.
\end{align}
From \eqref{est.lem.Lambda.Kol.rate.2} and the condition $1-|\mu_n|\geq \kappa \delta_n^2$, we have 
\begin{align}\label{proof.lem.Lambda.Kol.rate.2-2.-4}
\lim_{n\rightarrow \infty} \big (\delta_n^2 \| u_n \|_{L^2}^2 + \| \partial_y v_n\|_{L^2}^2 + \frac{1}{1-|\mu_n|} \| v_n \|_{L^2}^2 \big) =0\,,
\end{align}
which is essential in the proof below. It suffices to consider the case $\delta_\infty=0$; otherwise we have $\displaystyle \lim_{n\rightarrow \infty} \|u_n\|_{L^2}=0$ by \eqref{proof.lem.Lambda.Kol.rate.2-2.-4}, which implies $\displaystyle \lim_{n\rightarrow \infty} \|M_{\cos y} B_{2,l_n} u_n\|_{L^2}=0$ and we achieve the contradiction. 
We also note that we may assume $\mu_n\geq \frac12$ for all $n$, for the case $\mu_n\leq -\frac12$ is handle in the same manner.
Let us recall that $y_{n,j}\in S_{\mu_n}$ are the critical points, $\sin y_{n,j}=\mu_n$, such that $y_{n,1}\in (0,\frac{\pi}{2})$ and $y_{n,2}\in (\frac{\pi}{2}, \pi)$.
Then \eqref{proof.lem.Lambda.Kol.rate.2-2.-3} gives the identity 
\begin{align}\label{proof.lem.Lambda.Kol.rate.2-2.-5}
\mu_n v_n (y_{n,j}) + \delta_n f_n (y_{n,j})=0 \,.
\end{align}

Let $\kappa_5>0$ be fixed and sufficiently small number.
We decompose the interval $[-\frac{\pi}{2},\frac{3\pi}{2}]$ into $\tilde I_n$ and $\tilde I_n^c=[-\frac{\pi}{2},\frac{3\pi}{2}]\setminus \tilde I_n$, where 
\begin{align*}
\tilde I_n & = \big \{y\in [-\frac{\pi}{2},\frac{3\pi}{2}]~ \big |~{\rm dist}\, (y, S_{\mu_n}) \leq \kappa_5 \delta_n \big \}\,.
\end{align*}
We also set 
\begin{align*}
\widetilde{II}_n = \widetilde{II}_{n,1} \cup \widetilde{II}_{n,2}\,, \qquad \widetilde{II}_{n,1} = [-\frac{\pi}{2}, y_{n,1}+\kappa_5 \delta_n] \qquad \widetilde{II}_{n,2} = [y_{n,2}-\kappa_5 \delta_n, \frac{3\pi}{2}]\,.
\end{align*}
Note that $\tilde I_n \subset \widetilde{II}_{n}$ and $|\widetilde{II}_{n,j}|\geq \frac{1}{C}$ hold. 
Then we have from $|\cos y|\leq C |1-\mu_n|^\frac12$ for $y\in \tilde I_n$ and $|\cos y|\geq \frac{1}{C} |1-\mu_n|^\frac12$ for $y\in \widetilde{II}_{n,1} \cup \widetilde{II}_{n,2}$,
\begin{align}
\| M_{\cos y} B_{2,l_n} u_n \|_{L^2(\tilde I_n)} \leq  C |\tilde I_n|^\frac12 |1-\mu_n|^\frac12 \| B_{2,l_n} u_n \|_{L^\infty(\tilde{I}_n)}  
& \leq C \delta_n^\frac12 |1-\mu_n|^\frac12 \|  B_{2,l_n} u_n \|_{H^1(\widetilde{II}_n)}^\frac12 \|B_{2,l_n} u_n \|_{L^2(\widetilde{II}_n)}^\frac12 \nonumber \\
& \leq C \delta_n^\frac12 |1-\mu_n|^\frac14  \| (-A_{l_n})^\frac12 u_n \|_{L^2}^\frac12 \| M_{\cos y} B_{2,l_n} u_n \|_{L^2(\widetilde{II}_n)}^\frac12 \nonumber \\
& \leq C \delta_n^\frac12 |1-\mu_n|^\frac14  \| (-A_{l_n})^\frac12 u_n \|_{L^2}^\frac12 \nonumber \\
& \rightarrow 0 \quad (n\rightarrow \infty)\,. \nonumber 
\end{align}
Here we have used  \eqref{proof.lem.Lambda.Kol.rate.2-2.-1} in the last line.
Next, \eqref{proof.lem.Lambda.Kol.rate.2-2.-3} gives 
\begin{align}\label{proof.lem.Lambda.Kol.rate.2-2.-6}
B_{2,l_n} u_n & = \frac{\mu_n v_n+\delta_n f_n}{\mu_n - \sin y} \,.
\end{align}
We decompose $\tilde I_n^c$ as $\tilde I_n^c = (\tilde I_n^c \cap [0,\frac{\pi}{2}] )\cup (\tilde I_n^c \cap [\frac{\pi}{2}, \pi]) \cup (\tilde I_n^c \cap [\pi,2\pi])=: \tilde I_{n,1}^c \cup \tilde I_{n,2}^c \cup \tilde I_{n,3}^c$. Then, since $\frac12\leq \mu_n<1$ we find that there exists $C>0$ such that 
\begin{align}\label{est.sin'}
\begin{split}
& |\mu_n-\sin y|\geq \frac{|\cos y|}{C}\,  |y-y_{n,1}|\,, \qquad y\in \tilde I_{n,1}^c\,, \\
& |\mu_n-\sin y|\geq \frac{|\cos y|}{C} \, |y-y_{n,2}|\,, \qquad y\in \tilde I_{n,2}^c\,, \\ 
& |\mu_n - \sin y| \geq \frac{1}{2}\,, \qquad y\in \tilde I_{n,3}^c\,.
\end{split}
\end{align}
Here $y_{n,j}$ is a unique point of $S_{\mu_n}$ such that $y_{n,1}\in (0,\frac{\pi}{2})$ and $y_{n,2}\in (\frac{\pi}{2},\pi)$, respectively.
Then we have from \eqref{proof.lem.Lambda.Kol.rate.2-2.-6} and the definition of $\tilde I_{n,1}^c$,
\begin{align*}
 \| M_{\cos y} B_{2,l_n} u_n \|_{L^2(\tilde I_{n,1}^c)} & \leq C \| \frac{v_n - v_n (y_{n,1})}{y-y_{n,1}} \|_{L^2(\tilde I_{n,1}^c)} + \mu_n |v_n (y_{n,1})| \, \| \frac{\cos y}{\mu_n-\sin y}\|_{L^2(\tilde I_{n,1}^c)} + \frac{C\delta_n \| f_n \|_{L^2}}{\kappa_5 \delta_n}  \\
& \leq C \| \partial_y v_n \|_{L^2} + \frac{C|v_n (y_{n,1})|}{\kappa_5^\frac12 \delta_n^\frac12} + \frac{C\| f_n \|_{L^2}}{\kappa_5}\,.
\end{align*}
The same estimate holds for $\| M_{\cos y} B_{2,l_n} u_n \|_{L^2(\tilde I_{n,2}^c)}$, while we have from \eqref{est.sin'},
\begin{align*}
\| M_{\cos y} B_{2,l_n} u_n \|_{L^2 (\tilde I_{n,3}^c)} = \| \frac{\mu_n v_n + \delta_n f_n}{\mu_n-\sin y} \|_{L^2 (\tilde I_{n,3}^c)}\leq C  \big ( \| v_n \|_{L^2} + \delta_n \| f_n \|_{L^2} \big ) \,.
\end{align*}
Collecting these, we obtain from \eqref{proof.lem.Lambda.Kol.rate.2-2.-4},
\begin{align}\label{proof.lem.Lambda.Kol.rate.2-2.-7}
\limsup_{n\rightarrow \infty} \| M_{\cos y} B_{2,l_n} u_n \|_{L^2} \leq  \frac{C}{\kappa_5^\frac12} \sum_{j=1,2} \limsup_{n\rightarrow \infty} \frac{|v_n (y_{n,j})|}{\delta_n^\frac12}\,.
\end{align}
Let us estimate $\frac{|v_n (y_{n,j})|}{\delta_n^\frac12}$. It suffices to consider the case $j=1$, for the case $j=2$ is handled in the same manner.
The argument is similar to Step 2 in the proof of  \eqref{est.lem.Lambda.Kol.rate.2}.
For sufficiently small $\kappa'>0$ as above, we set $\widetilde T_n$ as
\begin{align*}
\widetilde T_n = [y_{n,1}, y_{n,1}+ {\kappa'}^2 \delta_n] \subset (0,\frac{\pi}{2}) \,. 
\end{align*}
We take $z_{n,1}\in \widetilde T_n$ so that 
\begin{align*}
|f(z_{n,1})|^2 \leq \frac{C}{{\kappa'}^2\delta_n} \| f_n \|_{L^2(\widetilde T_n)}^2 \leq \frac{C}{{\kappa'}^2 \delta_n} \| f_n \|_{L^2}^2\,.
\end{align*}
Then we have from \eqref{proof.lem.Lambda.Kol.rate.2-2.-3},
\begin{align}
\mu_n | v_n (y_{n,1})| = \delta_n |f_n (y_{n,1})| & \leq \delta_n  |f_n (y_{n,1})-f_n (z_{n,1})| +\delta_n  |f_n (z_{n,1})|  \nonumber \\
& \leq  {C\kappa'} \delta_n^\frac32 \| \partial_y f_n \|_{L^2(\widetilde T_n)} + \frac{C\delta_n^\frac12}{{\kappa'}} \| f_n \|_{L^2} \,.  \nonumber 
\end{align}
Next by using the identity \eqref{proof.prop.Lambda.Kol.rate.2.6} for $\delta_n \partial_y f_n$ and also \eqref{proof.lem.Lambda.Kol.rate.2-2.-1},
\begin{align*}
\delta_n \| \partial_y f_n \|_{L^2 (\widetilde T_n)} \leq \| M_{\cos y} B_{2,l_n} u_n \|_{L^2(\widetilde T_n)} + \| \partial_y v_n \|_{L^2}  + C \delta_n |1-\mu_n|^\frac12 \| \partial_y u_n \|_{L^2} \leq C\,,
\end{align*}
and thus, 
\begin{align*}
\limsup_{n\rightarrow \infty} \frac{|v_n (y_{n,1})|}{\delta_n^\frac12} \leq C\kappa' \rightarrow 0 \quad (\kappa'\rightarrow 0)\,.
\end{align*}
Hence \eqref{proof.lem.Lambda.Kol.rate.2-2.-7} implies $\displaystyle \limsup_{n\rightarrow \infty} \| M_{\cos y} B_{2,l_n} u_n \|_{L^2}=0$, which contradicts with \eqref{proof.lem.Lambda.Kol.rate.2-2.-1}. The proof of \eqref{est.lem.Lambda.Kol.rate.2'} is complete.
This completes the proof of Lemma \ref{lem.Lambda.Kol.rate.2}.

\

Finally we consider the case $|\mu|<\frac12$.
The proof is similar to Case 2 in the proof of Lemma \ref{lem.Lambda.Kol.rate.2}.
The only difference is the influence of the projection $\mathbb{Q}_l$ when $l=\pm 1$, which yields an additional nonlocal term in the limit equation when we perform the contradiction argument.
\begin{lem}\label{lem.Lambda.Kol.rate.3}
Let $\kappa\in (0,1)$ be the number in Lemma \ref{lem.Lambda.Kol.rate.1}.
There exists $C>0$  such that if $\delta\in (0,1]$, $l\in \Z\setminus\{0\}$, and $\mu\in \R$ with $|\mu| < \frac12$, then 
\begin{align}\label{est.lem.Lambda.Kol.rate.3}
\begin{split}
\delta^2 \| u\|^2 + \| (-A_l)^{-\frac12} u \|_{L^2}^2 +  \| A_l^{-1} u\|_{L^2}^2 & \leq C \bigg ( \delta^{-2} \| \mathbb{Q}_l (\mu-\hat{\Lambda}_l) u \|_{L^2}^2 + \delta^6 \| (-A_l)^\frac12 u \|_{L^2}^2  \bigg ) \,, 
\end{split}
\end{align}
and 
\begin{align}\label{est.lem.Lambda.Kol.rate.4}
\begin{split}
\| M_{\cos y} B_{2,l} u\|_{L^2}^2 & \leq C \bigg ( \delta^{-2} \| \mathbb{Q}_l  (\mu-\hat{\Lambda}_l) u \|_{L^2}^2 + \delta^2 \| (-A_l)^\frac12 u \|_{L^2}^2  \bigg ) \,, 
\end{split}
\end{align}
for all $u\in H^1 (\T) \cap Y_l$. Here $\mathbb{Q}_l : L^2 (\T) \rightarrow Y_l$ is the orthogonal projection on $Y_l$.
\end{lem}

%\begin{rem}{\rm In \eqref{est.lem.Lambda.Kol.rate.4} the projection $\mathbb{Q}_l$ is not acted on the term $(\mu-\hat{\Lambda}_l)u$, which crucially relied on the fact that the eigenspace of $\hat{\Lambda}_l$ is contained in ${\rm Ker}\, B_{2,l}$. } \end{rem}

Again the proof consists of several steps. We first consider \eqref{est.lem.Lambda.Kol.rate.3}.
% \begin{lem}\label{lem.Lambda.Kol.rate.3-1} The inequality \eqref{est.lem.Lambda.Kol.rate.3} holds.

% \end{lem}

\noindent {\it Proof of  \eqref{est.lem.Lambda.Kol.rate.3}.} The proof is very similar to the 
proof of  \eqref{est.lem.Lambda.Kol.rate.2} and is based on the contradiction argument. Again it suffices to consider  real valued functions and to show the claim for $u\in H^2(\T)\cap Y_l$.
Suppose that the estimate 
\begin{align}\label{proof.lem.Lambda.Kol.rate.3.1}
\begin{split}
\delta^2 \| u \|_{L^2}^2 + \| (-A_l)^{-\frac12} u \|_{L^2}^2 + &  \| A_l^{-1} u\|_{L^2}^2 \leq C \bigg ( \delta^{-2} \| \mathbb{Q}_l (\mu  - \hat{\Lambda}_l) u \|_{L^2}^2 + \delta^6 \| (-A_l)^\frac12  u\|_{L^2}^2 \bigg )\,, \\
& \delta\in (0,1]\,, \quad l\in \Z\setminus\{0\}\,, \quad  |\mu|<\frac12\,, \quad u\in H^2 (\T; \R) \cap Y_l
\end{split}
\end{align}
does not hold.
Then there exist $\{\delta_n, l_n, \mu_n\}_{n\in \N}$, $\delta_n \in (0,1]$, $l_n\in \Z\setminus\{0\}$, $\mu_n \in (-\frac12,\frac12)$, and $\{u_n\}\subset H^2 (\T;\R)\cap Y_{l_n}$ such that 
\begin{align*}
& \lim_{n\rightarrow \infty} \delta_n=\delta_\infty\in [0,1]\,, \quad \lim_{n\rightarrow \infty} l_n = l_\infty\in \{\pm \infty\}\cup \Z\setminus\{0\}\,,\\
& \lim_{n\rightarrow \infty}\mu_n=\mu_\infty \in [-\frac12, \frac12]\,,
\end{align*}
and 
\begin{align}\label{proof.lem.Lambda.Kol.rate.3.2}
\begin{split}
& \delta_n^2 \| u_n \|_{L^2}^2 + \| (-A_{l_n})^{-\frac12} u_n \|_{L^2}^2 + \| A_{l_n}^{-1} u_n \|_{L^2}^2 =1\,, \\
& \lim_{n\rightarrow \infty}  \bigg ( \delta_n^{-2}  \| \mathbb{Q}_{l_n}  (\mu_n - \hat{\Lambda}_{l_n}) u_n \|_{L^2}^2 + \delta_n^6 \| (-A_{l_n} )^\frac12 u_n \|_{L^2}^2 \bigg ) =0\,.
\end{split}
\end{align}
We first observe that $\delta_\infty=0$, otherwise we have $\|(-A_{l_n})^\frac12 u_n \|_{L^2}\rightarrow 0  \, (n\rightarrow \infty)$ due to the second condition in \eqref{proof.lem.Lambda.Kol.rate.3.2}, from which we easily reach a contradiction to the normalized condition in \eqref{proof.lem.Lambda.Kol.rate.3.2}. 
Moreover, if $|l_\infty|\ne 1$ then the situation is exactly the same as the case $\frac12 \leq \mu_\infty<1$ in the proof of Lemma \ref{lem.Lambda.Kol.rate.2}, for $\mathbb{Q}_l=I$ when $|l|\ne 1$. Therefore, it remains to consider the case $|l_\infty|=1$ and $|\mu_\infty|\leq \frac12$. Let us focus on the case $l_\infty=1$ and $0\leq \mu_\infty\leq \frac12$; the other cases are handled in the same manner. Then we may assume that $l_n=1$ for all $n$ by taking a subsequence if necessary,
though we often keep the notation $l_n$ for convenience.
 
Set
\begin{align}\label{proof.lem.Lambda.Kol.rate.3.3}
f_n = \delta_n^{-1} (\mu_n  - \hat{\Lambda}_{l_n} ) u_n\,, \qquad v_n = A_{l_n}^{-1} u_n\,.
\end{align}
We may assume that $v_n$ converges to a function $v_\infty$ strongly in $L^2(\T)$ and weakly in $H^1(\T)$.
Since $u_n$ is real valued, so is $v_n$, and $v_n$ satisfies 
$(\mu_n  - M_{\sin y}) \big ( A_{l_n} + 1) v_n  = \mu_n v_n + \delta_n f_n$, and then by $l_n^2=1$,
\begin{align}\label{proof.lem.Lambda.Kol.rate.3.4}
(\mu_n  - M_{\sin y}) \partial_y^2 v_n  = \mu_n v_n + \delta_n f_n\,.
\end{align}
It is convenient to introduce the value
\begin{align*}
\vartheta_n =
& \displaystyle \frac{1}{2\pi} \int_0^{2\pi} \delta_n f_n \dd y\,,
\end{align*}
which gives $\delta_n f_n = \delta_n \mathbb{Q}_{l_n} f_n + \vartheta_n$. 
The value $\vartheta_n$ is computed from \eqref{proof.lem.Lambda.Kol.rate.3.4} and the condition $v_n \in Y_{l_n}$ as 
\begin{align}\label{proof.lem.Lambda.Kol.rate.3.5}
\vartheta_n = -\frac{1}{2\pi} \int_0^{2\pi} \sin y \, \partial_y^2 v_n \dd y = \frac{1}{2\pi} \int_0^{2\pi} \sin y\, v_n \dd y\,.
\end{align}
Then \eqref{proof.lem.Lambda.Kol.rate.3.4} is written as 
\begin{align}\label{proof.lem.Lambda.Kol.rate.3.6}
(\mu_n  - M_{\sin y}) \partial_y^2  v_n  = \mu_n  v_n + \vartheta_n + \delta_n \mathbb{Q}_{l_n} f_n\,.
\end{align}
The trace relation of \eqref{proof.prop.Lambda.Kol.rate.2.4'} in the present case is 
\begin{align}\label{proof.lem.Lambda.Kol.rate.3.7}
\mu_n v_n (y_{\mu_n}) + \vartheta_n +  \delta_n (\mathbb{Q}_{l_n} f_n ) (y_{\mu_n}) =0\,, \qquad y_{\mu_n}\in S_{\mu_n} \,.
\end{align}
Here $S_{\mu_n}$ is the set of critical points.
A  key difference from the proof of \eqref{est.lem.Lambda.Kol.rate.2} is that the role of $\mu_n v_n (y_{n,j})$, where $y_{n,j}$ is the critical point, has to be replaced by $\mu_n v_n (y_{n,j}) + \vartheta_n$, and similarly, the role of $f_n$ is replaced by $\mathbb{Q}_{l_n} f_n$. 
The other part of the argument is similar to the proof of \eqref{est.lem.Lambda.Kol.rate.2} for the case $\mu_\infty<1$.
Indeed, the same argument as Step 1 in the proof of \eqref{est.lem.Lambda.Kol.rate.2} leads to 
\begin{align}\label{proof.lem.Lambda.Kol.rate.3.8}
\lim_{n\rightarrow \infty} \delta_n \| u_n \|_{L^2} =0\,,
\end{align}  
while the argument of Step 2 gives 
\begin{align*}
\lim_{n\rightarrow \infty} \, (\mu_n v_n (y_{n,j}) + \vartheta_n ) = \lim_{n\rightarrow \infty} \delta_n (\mathbb{Q}_{l_n}) f_n (y_{n,j}) = 0\,.
\end{align*}
In particular, we have 
\begin{align}\label{proof.lem.Lambda.Kol.rate.3.9}
\mu_\infty v_\infty (y_{\mu_\infty}) + \vartheta_\infty =0\,, \qquad y_{\mu_\infty} \in S_{\mu_\infty}\,,
\end{align}
where  
\begin{align*}
\vartheta_\infty = \lim_{n\rightarrow \infty} \frac{1}{2\pi} \int_0^{2\pi} \sin y \, v_n \, \dd y = \frac{1}{2\pi} \int_0^{2\pi} \sin y\, v_\infty \, \dd y\,.
\end{align*}
To estimate $\|\partial_y v_n\|_{L^2}^2$ we use 
\begin{align}\label{proof.lem.Lambda.Kol.rate.3.10}
\partial_y^2 v_n = \frac{\mu_n v_n +\vartheta_n + \delta_n \mathbb{Q}_{l_n} f_n}{\mu_n  - \sin y}\,.
\end{align}
By taking the inner product with $v_n$, we obtain 
\begin{align}\label{proof.lem.Lambda.Kol.rate.3.11}
\|\partial_y v_n \|_{L^2}^2 =  - \int_{-\frac{\pi}{2}}^{\frac{3\pi}{2}} \frac{\big ( \mu_n v_n + \vartheta_n  + \delta_n \mathbb{Q}_{l_n} f_n\big ) v_n }{\mu_n -\sin y}  \dd y  \,.
\end{align}  
Then we apply the argument of Step 4 by replacing $\mu_n v_n$ and $f_n$ there by $\mu_n v_n + \vartheta_n$ and $\mathbb{Q}_{l_n} f_n$, respectively.
We obtain the estimate of the form
\begin{align}\label{proof.lem.Lambda.Kol.rate.3.12}
\limsup_{n\rightarrow \infty}  \| \partial_y v_n \|_{L^2}^2 \leq \frac{C}{\kappa_3} \limsup_{n\rightarrow \infty} \| v_n \|_{L^2}^2 + C \kappa_3^\frac12
\end{align}
for any sufficiently small $\kappa_3>0$, and also obtain the estimate 
\begin{align}\label{proof.lem.Lambda.Kol.rate.3.13}
\big |\langle \partial_y v_\infty, \partial_y \varphi\rangle_{L^2}\big | \leq C \| \varphi \|_{L^2}
\end{align}
for any $\varphi\in H^2 (\T)$. This ensure the regularity $v_\infty \in H^2 (\T)$.
Estimate \eqref{proof.lem.Lambda.Kol.rate.3.12} together with \eqref{proof.lem.Lambda.Kol.rate.3.8} and the normalized condition \eqref{proof.lem.Lambda.Kol.rate.3.2} implies that we may assume $\inf_n \|v_n\|_{L^2}>0$, and thus, the limit $v_\infty$ must be nontrivial.
We can also show that the limit $v_\infty\in H^2(\T) \cap Y_1$ satisfies 
\begin{align}\label{proof.lem.Lambda.Kol.rate.3.14}
(\mu_\infty - \sin y) \partial_y^2 v_\infty = \mu_\infty v_\infty + \vartheta_\infty\,.
\end{align}
If $\mu_\infty=0$ then $\vartheta_\infty=0$ by \eqref{proof.lem.Lambda.Kol.rate.3.9}, and thus, $\partial_y^2 v_\infty=0$ by \eqref{proof.lem.Lambda.Kol.rate.3.14}. Hence $v_\infty$ is a constant. Since $v_\infty\in Y_1$ we conclude that the constant $v_\infty$ must be zero. This is a contradiction.
 The proof is complete for the case $\mu_\infty=0$.
When $0<\mu_\infty\leq \frac12$ 
let $y_{\infty,j}\in S_{\mu_\infty}$, $j=2,3$, be such that $y_{\infty,2}\in (\frac{\pi}{2},\pi)$ and $y_{n,3}\in (2\pi, \frac{5\pi}{2})$. Then $\mu_\infty -\sin y\geq 0$ for $y\in (y_{\infty,2},y_{\infty,3})$.
Thus we see 
\begin{align*}
\int_{y_{\infty,2}}^{y_{\infty,3}} \partial_y^2 v_\infty \,  (\mu_\infty v_\infty + \vartheta_\infty ) \dd y = \int_{y_{\infty,2}}^{y_{\infty,3}} \frac{(\mu_\infty v_\infty+\vartheta_\infty)^2}{\mu_\infty -\sin y}  \dd y\,,
\end{align*}
which makes sense by the $H^2$ regularity of $v_\infty$ and the condition \eqref{proof.lem.Lambda.Kol.rate.3.9}. 
The integration by parts and \eqref{proof.lem.Lambda.Kol.rate.3.14} imply
\begin{align}\label{proof.lem.Lambda.Kol.rate.3.15}
\mu_\infty \int_{y_{\infty,2}}^{y_{\infty,3}} |\partial_y v_\infty|^2 \dd y + \int_{y_{\infty,2}}^{y_{\infty,3}} \frac{(\mu_\infty v_\infty + \vartheta_\infty)^2}{\mu_\infty - \sin y} \dd y =0\,.
\end{align}
Therefore, we conclude that $\mu_\infty v_\infty + \vartheta_\infty=0$ on $[y_{\infty,2},y_{\infty,3}]$.
Set $w_\infty = \partial_y v_\infty\in H^1 (\T)$, which then satisfies $w_\infty =0$ on  $[y_{\infty,2},y_{\infty,3}]$.
Moreover, from \eqref{proof.lem.Lambda.Kol.rate.3.15} and also from $\vartheta_\infty=-\mu_\infty v(y_{\infty,j})$, we see
\begin{align}\label{proof.lem.Lambda.Kol.rate.3.16}
\partial_y w_\infty = \frac{\mu_\infty v_\infty +\vartheta_\infty}{\mu_\infty-\sin y} =  \frac{\mu_\infty}{\mu_\infty-\sin y} \int_{y_{\infty,2}}^y w_\infty \dd z\,, \quad y>y_{\infty,2}\,, \qquad w_\infty (y_{\infty,2})=0\,.
\end{align}
Then it is easy to see that $w_\infty =0$ for $y\in (y_{\infty,2}, y_{\infty,2} + \tau)$ for some $\tau>0$, and thus, $w_\infty =0$ for all $y>y_{\infty,2}$. Hence $w_\infty=0$, i.e., $v_\infty$ is a constant. Since $v_\infty \in Y_1$ we must have $v_\infty=0$, which is a contradiction. 
The proof of \eqref{est.lem.Lambda.Kol.rate.3} is complete.

\

% \begin{lem}\label{lem.Lambda.Kol.rate.3-2} The inequality \eqref{est.lem.Lambda.Kol.rate.4} also holds.
% \end{lem}

\noindent {\it Proof of  \eqref{est.lem.Lambda.Kol.rate.4}.} 
The proof is again very similar with the proof of \eqref{est.lem.Lambda.Kol.rate.2'} and is based on a contradiction argument. 
Then the problem is reduced to the analysis of the sequence $\{\tilde u_n\}$, where $\tilde u_n = M_{\cos y} B_{2,l_n} u_n$, $u_n\in H^2 (\T) \cap Y_{l_n}$.
As in proof of  \eqref{est.lem.Lambda.Kol.rate.3},  it suffices to consider the case $|l_n|=1$ and $\delta_\infty=0$,
and without loss of generality we may assume that $l_n=1$ for all $n$. 
From the hypothesis of the contradiction argument, we have the convergence
\begin{align}\label{proof.lem.Lambda.Kol.rate.3-2.1}
\lim_{n\rightarrow \infty} \bigg ( \delta_n^{-2} \| (\mu_n - \hat{\Lambda}_{l_n}) u_n \|_{L^2}^2 + \delta_n^2 \| (-A_{l_n})^\frac12 u_n \|_{L^2}^2 \bigg ) =0\,,
\end{align}
and therefore, from \eqref{est.lem.Lambda.Kol.rate.3} we have 
\begin{align}\label{proof.lem.Lambda.Kol.rate.3-2.2}
\lim_{n\rightarrow \infty} (\delta_n^2 \| u_n \|_{L^2}^2 + \|  v_n \|_{H^1}^2) =0\,.
\end{align}
Then we can apply the same argument as in the proof of  
% Lemma \ref{lem.Lambda.Kol.rate.2-2} 
\eqref{est.lem.Lambda.Kol.rate.2'}, for the argument there relies only on \eqref{proof.lem.Lambda.Kol.rate.3-2.1} and \eqref{proof.lem.Lambda.Kol.rate.3-2.2}. 
The only difference from  the proof of  \eqref{est.lem.Lambda.Kol.rate.2'} is that the role of $\mu_n v_n$ and $f_n$ is again replaced by $\mu_n v_n + \vartheta_n$ and $\mathbb{Q}_{l_n} f_n$, respectively.
We omit the details. The proof of Lemma \ref{lem.Lambda.Kol.rate.3} is complete.

\

%\noindent {\it End of the proof of Lemma \ref{lem.Lambda.Kol.rate.3}.}
%Lemma \ref{lem.Lambda.Kol.rate.3} directly follows from Lemmas \ref{lem.Lambda.Kol.rate.3-1} and \ref{lem.Lambda.Kol.rate.3-2}. The proof is complete.

\

Let us set for $m\geq 1$ and $\mu \in \R$,
\begin{align}\label{def.h_1}
h_{1} (m, \mu) =
\begin{cases}
& \displaystyle 0 \qquad \qquad \qquad \qquad  \quad \text{if }~~~ |\mu| >1+\frac{\kappa}{m}\,, \\
& \displaystyle m^{-\frac12} \qquad \qquad \qquad \quad  ~ ~~\text{if }~~~ 1-\frac{\kappa}{m} < |\mu| \leq 1 + \frac{\kappa}{m}\,,\\
& \displaystyle  m^{-1} ( 1- |\mu|)^{-\frac12} \qquad ~ \quad ~ ~\text{if }~~~  |\mu| \leq 1- \frac{\kappa}{m}\,,
\end{cases}
\end{align}
\begin{align}\label{def.h_2}
h_{2} (m, \mu) = 
\begin{cases}
& \displaystyle 0 \qquad \qquad \qquad   \qquad ~~ \text{if }~~~ |\mu| >1+\frac{\kappa}{m^2}\,,\\
& \displaystyle m^{-2}  \qquad \qquad  \qquad \quad \text{if }~~ 1 -\frac{\kappa}{m^2} < |\mu|  \leq 1+ \frac{\kappa}{m^2}\,,\\
& \displaystyle m^{-1} (1- |\mu|)^\frac12 \qquad ~~~~~ \text{if }~~|\mu| \leq 1- \frac{\kappa}{m^2}\,,
\end{cases}
\end{align}
and 
\begin{align}\label{def.h_3}
h_{3} (m, \mu) = 
\begin{cases}
& \displaystyle 0 \qquad \qquad \qquad  \qquad ~ ~~ \text{if }~~~ |\mu| >1+\frac{\kappa}{m^2}\,,\\
& \displaystyle m^{-2}  \qquad \qquad \qquad  \qquad \text{if }~~ 1 -\frac{\kappa}{m^2} < |\mu|  \leq 1 + \frac{\kappa}{m^2}\,,\\
& \displaystyle m^{-3} (1- |\mu|)^{-\frac12} \qquad \quad ~~~ \text{if }~~|\mu| \leq 1- \frac{\kappa}{m^2}\,.
\end{cases}
\end{align}
Note that each $h_j$ satisfies $h_{j}(m,-\mu)=h_{j} (m, \mu)$ and $\displaystyle \lim_{m\rightarrow \infty} \sup_{\mu\in \R} h_j (m,\mu)=0$. Moreover, we have $h_3(m,\mu)\leq C h_2(m,\mu)$. 
By Lemmas \ref{lem.Lambda.Kol.rate.1}, \ref{lem.Lambda.Kol.rate.2}, and \ref{lem.Lambda.Kol.rate.3}, we obtain
\begin{prop}\label{thm.Lambda.Kol.rate.1} There exist $C, \, \kappa>0$ such that for all $m\geq 1$, $\mu\in \R$, and $l\in \Z\setminus\{0\}$,
\begin{align}
\begin{split}
\| u\|_{L^2}^2 & \leq C \bigg ( m^2 \| (\mu -  \hat{\Lambda}_l) u \|_{L^2}^2 +  h_{1}^2 (m,\mu) \, \| (-A_l)^\frac12 u\|_{L^2}^2 \bigg )\,, \\
& \qquad \qquad \qquad \qquad \qquad  u\in H^1(\T)\,,\quad \frac12 \leq | \mu |\leq 1+ \frac{\kappa}{m}\,, \label{est.thm.Lambda.Kol.rate.1.1}
\end{split}\\
\begin{split}
\| u\|_{L^2}^2 & \leq C \bigg ( m^2 \| \mathbb{Q}_l (\mu - \hat{\Lambda}_l) u \|_{L^2}^2 +  h_{1}^2 (m,\mu) \, \| (-A_l)^\frac12 u \|_{L^2}^2 \bigg )\,,\\
& \qquad  \quad \qquad \qquad \qquad u\in H^1(\T)\cap Y_l\,, \quad |\mu|<\frac12\,.\label{est.thm.Lambda.Kol.rate.1.2}
\end{split}\\
\begin{split}
\| M_{\cos y} B_{2,l} u \|_{L^2}^2 & \leq C \bigg ( m^2 \| (\mu - \hat{\Lambda}_l) u \|_{L^2}^2 +  h_{2}^2 (m,\mu) \, \| (-A_l)^\frac12 u\|_{L^2}^2 \bigg )\,, \\
& \qquad \qquad \qquad \quad \qquad u\in H^1(\T)\,,\quad \frac12 \leq |\mu|\leq 1+ \frac{\kappa}{
 m^2}\,,\label{est.thm.Lambda.Kol.rate.1.3}
\end{split}\\
\begin{split}
\| M_{\cos y} B_{2,l} u \|_{L^2}^2 & \leq C \bigg ( m^2 \| \mathbb{Q}_l (\mu - \hat{\Lambda}_l) u \|_{L^2}^2 +  h_{2}^2 (m,\mu)  \, \| (-A_l)^\frac12 u \|_{L^2}^2 \bigg )\,,\\
& \qquad  \quad \qquad \qquad \qquad u\in H^1(\T) \cap Y_l\,,\quad | \mu| < \frac12\,,\label{est.thm.Lambda.Kol.rate.1.4}
\end{split}
\end{align}
and 
\begin{align}
\begin{split}
\| (-A_l)^{-\frac12} u \|_{L^2}^2 & \leq C \bigg ( m^2 \| (\mu - \hat{\Lambda}_l) u \|_{L^2}^2 +  h_{3}^2 (m,\mu) \, \| (-A_l)^\frac12 u\|_{L^2}^2 \bigg )\,, \\
& \qquad \qquad \qquad \quad \qquad u\in H^1(\T)\,,\quad \frac12 \leq |\mu|\leq 1+ \frac{\kappa}{
 m^2}\,,\label{est.thm.Lambda.Kol.rate.1.1'}
\end{split}\\
\begin{split}
\| (-A_l)^{-\frac12} u \|_{L^2}^2 & \leq C \bigg ( m^2 \| \mathbb{Q}_l (\mu -  \hat{\Lambda}_l) u \|_{L^2}^2 +  h_{3}^2 (m,\mu)  \, \| (-A_l)^\frac12 u \|_{L^2}^2 \bigg )\,,\\
& \qquad  \quad \qquad \qquad \qquad u\in H^1(\T) \cap Y_l\,,\quad | \mu| < \frac12\,.\label{est.thm.Lambda.Kol.rate.1.2'}
\end{split}
\end{align}
On the other hand, if $|\mu|>1$ then 
\begin{align}
\| u\|_{L^2}^2 & \leq \frac{C}{(|\mu|-1)^2} \| (\mu -  \hat{\Lambda}_l) u\|_{L^2}^2\,, \qquad u\in H^1 (\T)\,,\label{est.thm.Lambda.Kol.rate.1.5}\\
\| M_{\cos y} B_{2,l} u\|_{L^2}^2 & \leq \frac{C}{(|\mu|-1)} \| (\mu - \hat{\Lambda}_l) u\|_{L^2}^2\,, \qquad u\in H^1 (\T)\,,\label{est.thm.Lambda.Kol.rate.1.6}
\end{align}
and 
\begin{align}
\| (-A_l)^{-\frac12} u\|_{L^2}^2 & \leq \frac{C}{|\mu| \, (|\mu|-1)} \| (\mu - \hat{\Lambda}_l) u\|_{L^2}^2\,, \qquad u\in H^1 (\T)\,.\label{est.thm.Lambda.Kol.rate.1.7}
\end{align}
\end{prop}

\

Note that the constants $C$ and $m_0$ in Proposition \ref{thm.Lambda.Kol.rate.1} are independent of $l\in \Z\setminus\{0\}$. Let us recall that $L_{\alpha,l}$ is defined as $L_{\alpha,l} = A_l - i \alpha l \hat{\Lambda}_l$, and thus, it is convenient to introduce 
\begin{align}\label{def.tilde_alpha_l}
\tilde \alpha = \tilde \alpha (l) = \alpha l\,.
\end{align}
We are interested in the estimate of $\|(i\lambda + \mathbb{Q}L_{\alpha,l})^{-1} \|_{Y_l\rightarrow Y_l}$ by applying Theorem \ref{thm.abstract.2}. In particular, the dependence of the estimate on $\tilde \alpha$ is important. Let us recall that $B_3 = B_1 T_l A_l^{-1} + M_{\cos y} B_{2,l}=M_{\sin y} (\partial_y - l) A_l^{-1} + M_{\cos y} B_{2,l}$, and hence, 
\begin{align*}
\| B_3 u \|_{L^2} \leq C \| (-A_l)^{-\frac12} u\|_{L^2} + \| M_{\cos y} B_{2,l} u \|_{L^2}\,.
\end{align*}
Since $h_3(m, \mu)\leq C h_2 (m,\mu)$ by their definitions, it suffices to consider the function $F(\tilde \alpha,\mu)$ defined as
\begin{align}\label{def.F_l}
F (\tilde \alpha,\mu) = \inf_{m_1,m_2\geq m_0} \bigg ( \frac{m_1}{|\tilde \alpha|}  + \frac{m_1^2m_2^2}{\tilde \alpha^2} + \frac{m_1^2h_{2} (m_2,\mu)}{|\tilde \alpha|}  + h_{1} (m_1,\mu)^2 \bigg ) \,.
\end{align}
Here $h_{j}$ are defined by \eqref{def.h_1} and \eqref{def.h_2}.
Our aim is to obtain the upper bound for  $F(\tilde \alpha,\mu)$.

\noindent {\bf Case 1: $|\mu| > 1 + \frac{\kappa^\frac34}{|\tilde \alpha|^\frac12}$.}
In this case let us take 
\begin{align}\label{proof.F_l.6}
m_1 = m_2^2 = \frac{2\kappa}{|\mu|-1}\,.
\end{align}
Then we have 
\begin{align*}
1+\frac{\kappa}{m_1} = 1+\frac{\kappa}{m_2^2} = 1 +  \frac{|\mu|-1}{2} <|\mu|\,.
\end{align*}
Thus, $h_{j} (m_j,\mu)=0$ by the definition, and 
\begin{align*}
& \frac{m_1}{|\tilde \alpha|} =\frac{2\kappa}{|\tilde \alpha| (|\mu|-1)} \,,\\
& \frac{m_1^2m_2^2}{\tilde \alpha^2} \leq \frac{C}{\tilde \alpha^2 (|\mu|-1)^3} \leq \frac{C}{|\tilde \alpha| (|\mu|-1)} \,.
\end{align*}
Hence we obtain
\begin{align}\label{proof.F_l.7}
\frac{m_1}{|\tilde \alpha|}  + \frac{m_1^2m_2^2}{\tilde \alpha^2} + \frac{m_1^2h_{2} (m_2,\mu)}{|\tilde \alpha|}  + h_{1} (m_1,\mu)^2  \leq \frac{C}{|\tilde \alpha| (|\mu|-1)}\,,
\end{align}
as desired.

\noindent {\bf Case 2: $\displaystyle  1- \frac{\kappa^\frac34}{|\tilde \alpha|^\frac12} <  |\mu| \leq 1 + \frac{\kappa^\frac34}{|\tilde \alpha|^\frac12}$.} In this case we take $m_1$ and $m_2$ as 
\begin{align}\label{proof.F_l.4}
m_1 = m_2^2 = \kappa^{\frac12} |\tilde \alpha|^\frac12\,.
\end{align}
Then we see from $\kappa\in (0,1)$,
\begin{align*}
\frac{\kappa}{m_1} = \frac{\kappa}{m_2^2} = \frac{\kappa^{1-\frac12}}{|\tilde \alpha|^\frac12} >\frac{\kappa^\frac34}{|\tilde \alpha|^\frac12} \geq  \big | 1-|\mu|\big |\,.
\end{align*} 
Hence $h_{1} (m_1,\mu) = m_1^{-\frac12}$ and $h_{2} (m_2,\mu) = m_2^{-2}$ for this choice of $m_1$ and $m_2$.
We can also check that there exist $C, C'>0$ depending only on $\kappa$ such that 
\begin{align*}
\frac{1}{C} \frac{m_1^2 m_2^2}{\tilde \alpha^2} \leq  \frac{m_1^2 h_2 (m_2,\mu)}{|\tilde \alpha|} \leq  C' h_{1} (m_1,\mu)^2 \leq C \frac{m_1^2 m_2^2}{\tilde \alpha^2} \,.
\end{align*}
Let us compute the size of each term as in Case $1$:
\begin{align*}
& \frac{m_1}{|\tilde \alpha|} = \frac{\kappa^\frac12}{|\tilde \alpha|^\frac12} \,,\qquad  \frac{m_1^2m_2^2}{\tilde \alpha^2}  = \frac{\kappa^\frac32}{|\tilde \alpha|^\frac12}\,.
\end{align*}
Thus we have 
\begin{align}\label{proof.F_l.5}
\frac{m_1}{|\tilde \alpha|}  + \frac{m_1^2m_2^2}{\tilde \alpha^2} + \frac{m_1^2h_{2} (m_2,\mu)}{|\tilde \alpha|}  + h_{1} (m_1,\mu)^2 \leq \frac{C}{|\tilde \alpha|^{\frac12}} \,.
\end{align}

\noindent {\bf Case 3: $\displaystyle |\mu|\leq 1- \frac{\kappa^\frac34}{|\tilde \alpha|^\frac12}$.}
In this case we take $m_1$ and $m_2$ as
\begin{align}\label{proof.F_l.1}
m_1 = (\frac{|\tilde \alpha|}{1-|\mu|})^\frac13\,, \qquad m_2 = \bigg (|\tilde \alpha| (1-|\mu|)^\frac12 \bigg )^\frac13 \,.
\end{align}
Then we have 
\begin{align*}
\frac{\kappa}{m_1 (1-|\mu|)} =\frac{\kappa}{|\tilde \alpha|^\frac13 (1-|\mu|)^{\frac23}} \leq  \frac{\kappa}{|\tilde \alpha|^\frac13} (\frac{|\tilde \alpha|^\frac12}{\kappa^\frac34})^\frac23 = \kappa^\frac12 <1\,,
\end{align*}
and 
\begin{align*}
\frac{\kappa}{m_2^2 (1-|\mu|)} = \frac{\kappa}{|\tilde \alpha|^\frac23 (1-|\mu|)^\frac43}\leq \frac{\kappa}{|\tilde \alpha|^\frac23} (\frac{|\tilde \alpha|^\frac12}{\kappa^\frac34})^\frac43 =1\,.
\end{align*}
Hence $h_{1} (m_1,\mu) =  m_1^{-1} (1-|\mu|)^{-\frac12}$ and $h_{2} (m_2,\mu) = m_2^{-1} (1-|\mu|)^\frac12$ by their definitions. Thus it follows that
\begin{align*}
\frac{m_1^2m_2^2}{\tilde \alpha^2} = \frac{m_1^2h_{2} (m_2,\mu)}{|\tilde \alpha|} \qquad \text{for any} ~~m_1>0\,,
\end{align*}
and we can also check the balance 
\begin{align*}
\frac{m_1^2 h_{2} (m_2,\mu)}{|\tilde \alpha|} =h_{1} (m_1,\mu)^2 \,.
\end{align*}
Let us now compute the size of each term in the right-hand side of \eqref{def.F_l} for $m_1$ and $m_2$ defined as \eqref{proof.F_l.1}:
\begin{align*}
& \frac{m_1}{|\tilde \alpha|}  = \frac{1}{|\tilde \alpha|^\frac23 (1-|\mu|)^\frac13}\,,\\
& \frac{m_1^2 m_2^2}{\tilde \alpha^2}  = \frac{|\tilde \alpha|^\frac23}{(1-|\mu|)^\frac23} |\tilde \alpha|^\frac23 (1-|\mu|)^\frac13 \tilde \alpha^{-2} = \frac{1}{|\tilde \alpha|^\frac23 (1-|\mu|)^\frac13}\,.
\end{align*}
Thus we have 
\begin{align}\label{proof.F_l.3}
 \frac{m_1}{|\tilde \alpha|}  + \frac{m_1^2m_2^2}{\tilde \alpha^2} + \frac{m_1^2h_{2} (m_2,\mu)}{|\tilde \alpha^2|}  + h_{1} (m_1,\mu)^2 \leq \frac{C}{|\tilde \alpha|^\frac23 (1-|\mu|)^\frac13}\,.
\end{align}
Here $C$ is also  independent  of $l$.

As a summary, we have arrived at the upper bound of $F (\tilde \alpha,\mu)$ such that 
\begin{align}\label{proof.F_l.8}
F (\tilde \alpha,\mu) \leq C |\tilde \alpha|^{-\frac12} = C|\alpha l |^{-\frac12}\,,
\end{align}
where $C$ is independent of $\alpha$, $\mu$, and $l$.
Then Theorem \ref{thm.abstract.2} implies the following result.
\begin{thm}\label{thm.Lambda.Kol.rate.2} There exist positive numbers $c$ and $C$ such that the following statement holds for all sufficiently large $|\alpha|$. 
Let $\lambda\in \R$ and $l\in \Z\setminus \{0\}$. Then  $\mathbb{Q}_l L_{\alpha,l}$  in $Y_l = \mathbb{Q}_l L^2(\T)$ satisfies  
\begin{align}
\sup_{\zeta \in \sigma (\mathbb{Q}_l L_{\alpha,l})} \Re \zeta  \leq - c \, |\alpha l|^\frac12\,, \label{est.thm.Lambda.Kol.2.3}
\end{align}
and 
\begin{align}\label{est.thm.Lambda.Kol.2.4}
\| (i\lambda + \mathbb{Q}_l L_{\alpha,l})^{-1}\|_{Y_l \rightarrow Y_l}   \leq 
\begin{cases}
& \displaystyle \frac{C}{|\alpha l| \, (|\frac{\lambda}{\alpha l}|-1)} ~~\qquad {\rm if}~~|\frac{\lambda}{\alpha l}|>1+\frac{1}{|\alpha l|^\frac12}\,,\\
& \displaystyle \frac{C}{|\alpha l |^\frac12} \qquad \qquad \qquad {\rm if}~~1-\frac{1}{|\alpha l|^\frac12} < |\frac{\lambda}{\alpha l}| \leq 1 + \frac{1}{|\alpha l|^\frac12}\,,\\
& \displaystyle \frac{C}{|\alpha l|^\frac23 (1-|\frac{\lambda}{\alpha l}|)^\frac13} \qquad {\rm if}~~ |\frac{\lambda}{\alpha l}| \leq 1 - \frac{1}{|\alpha l|^\frac12}\,.
\end{cases}
\end{align}
Here $c$ and $C$ are independent of $\alpha$ and $l$.
\end{thm}

Recall that $L_{\alpha,l}=A_l-i\alpha l \hat{\Lambda}_l$. 
Theorem \ref{thm.Lambda.Kol.rate.2} follows from Theorem \ref{thm.abstract.2} and the estimates \eqref{proof.F_l.7}, \eqref{proof.F_l.5}, \eqref{proof.F_l.3}, and \eqref{proof.F_l.8} for $F_l(\tilde \alpha,\mu)$ with $\tilde \alpha = \alpha l$. Note that $\mathbb{Q}_l=I$ when $|l|\geq 2$.

\subsection{Estimate for semigroup}\label{subsec.Kol.semigroup}

The resolvent estimates in Theorem \ref{thm.Lambda.Kol.rate.2} provide a crucial information on the solution to the following nonstationary problem both in qualitative and quantitative point of views.
\begin{align}\label{proof.thm.semigroup.1}
\left\{
\begin{aligned}
& \frac{\dd  w}{\dd t} - L_\alpha w  =0\,, \qquad t>0\,, \\
& w|_{t=0}   = f\in L^2_0 (\T^2)\,.
\end{aligned}\right.
\end{align}
Note that the operator $L_\alpha = A- i \alpha \hat{\Lambda}$ is diagonalized in terms of the Fourier series with respect to the $x$ variable:
\begin{align}
L_\alpha = \oplus_{l\in \Z\setminus\{0\}} \, L_{\alpha,l}\,, \qquad L_{\alpha,l} = A_l - i \alpha l \hat{\Lambda}_l\,.
\end{align}
Hence the estimates of the solution $u$ to \eqref{proof.thm.semigroup.1} are obtained from the estimates for each Fourier mode $\mathcal{P}_l u$, which is given by the semigroup generated by $L_{\alpha,l}$ in $L^2 (\T)$. 
For the estimate of the semigroup $e^{t L_{\alpha,l}}$
it is convenient to use the representation in terms of  the Dunford integral
\begin{align}\label{proof.thm.semigroup.2}
w_l (t) = e^{t L_{\alpha,l}} f_l = \frac{1}{2\pi i} \int_\Gamma e^{t\zeta} (\zeta - L_{\alpha,l})^{-1} f_l \dd \zeta\,.
\end{align}
Here $f_l=(\mathcal{P}_l f) e^{-ilx} \in L^2 (\T)$, and $\Gamma$ is first taken as 
\begin{align}\label{proof.thm.semigroup.3}
\Gamma & =\big \{ \zeta \in \C~|~\Re \zeta = -\frac12\,, ~ |\Im \zeta|\leq 4 |\alpha l|  \, \big \}  \nonumber \\
& \quad  \cup \big \{ \zeta \in \C~|~\Im \zeta = \mp  (\Re \zeta + \frac12) \pm 4 |\alpha l|\,, ~  \Re \zeta\leq -\frac12\,  \big \} \nonumber \\
& =:  \Gamma_{0,-\frac12} +  \Gamma_{\pm, -\frac12}\,,
\end{align}
which is oriented counter-clockwisely.
We note that 
\begin{align}\label{proof.thm.semigroup.4}
\| (\zeta-L_{\alpha,l})^{-1} \|_{L^2(\T)\rightarrow L^2(\T)}\leq \frac{C}{|\Im \zeta|} \,, \qquad |\Im \zeta|\geq 4\alpha |l|
\end{align}
hold with a constant $C$ independent of $\alpha$.
Set $\mathbb{P}_l=I-\mathbb{Q}_l$, where $\mathbb{Q}_l: L^2(\T)\rightarrow Y_l=({\rm Ker}\, \hat{\Lambda}_l)^\bot$ is the orthogonal projection as used in the previous section.
Note that $\mathbb{P}_l=0$ when $|l|\geq 2$.
Our aim is to establish the estimates for $\mathbb{Q}_l e^{t L_{\alpha,l}}$ and $\mathbb{P}_l e^{t L_{\alpha,l}}$.
For the part $\mathbb{Q}_l e^{t L_{\alpha,l}}$ the fast dissipation is expected for large $\alpha$, while for the part $\mathbb{P}_l e^{t L_{\alpha,l}}$ the strong amplification is expected through the interaction term
\begin{align*}
i \alpha l \mathbb{P}_l \hat{\Lambda}_l \mathbb{Q}_l\,, \qquad |l|=1\,,
\end{align*}
which does not vanish due to the lack of the invariance of the space $Y_l=\mathbb{Q}_l L^2 (\T)$ under the action of $\hat{\Lambda}_l$, or in other words,
due to the lack of the symmetry of $\hat{\Lambda}_l$, for this term automatically vanishes if $\hat{\Lambda}_l$ is symmetric.
To estimate $\mathbb{Q}_l e^{t L_{\alpha,l}}$ we observe 
\begin{align*}
\mathbb{Q}_l e^{t L_{\alpha,l}} f & = \frac{1}{2\pi i} \int_\Gamma e^{t\zeta} \mathbb{Q}_l (\zeta-L_{\alpha,l})^{-1} f_l \dd \zeta  =  \frac{1}{2\pi i} \int_\Gamma e^{t\zeta} (\zeta- \mathbb{Q}_l L_{\alpha,l})^{-1} \mathbb{Q}_l f_l \dd \zeta \,.
\end{align*}
The last identity follows from \eqref{eq.rem.thm.abstract.1'}. 
We observe that for each $\lambda\in \R$ the set   
$$\big \{ \, \zeta\in \C~|~|\zeta - i \lambda | < \frac{1}{\| (i\lambda - \mathbb{Q}_l L_{\alpha,l})^{-1} \|_{Y_l\rightarrow Y_l}}  \, \big \}$$
is contained in the resolvent set of $\mathbb{Q}L_{\alpha,l}$ by the standard Neumann series argument,
and in particular, we have 
\begin{align}\label{proof.thm.semigroup.-1}
\| (\zeta - \mathbb{Q}_l L_{\alpha,l})^{-1} \|_{Y_l\rightarrow Y_l} \leq 2 \| (i\lambda - \mathbb{Q}_l L_{\alpha,l})^{-1} \|_{Y_l\rightarrow Y_l}
\end{align}
if $|\zeta - i \lambda | < \frac{1}{2\| (i\lambda - \mathbb{Q}_l L_{\alpha,l})^{-1} \|_{Y_l\rightarrow Y_l}}$. 
Then, in virtue of Theorem \ref{thm.Lambda.Kol.rate.2}, we can shift the integral  
$\int_\Gamma \dd\zeta $ to  $\int_{\Gamma_\alpha}\dd \zeta$ by the Cauchy theorem, where 
\begin{align*}
\int_{\Gamma_\alpha} \dd \zeta & = \int_{\Gamma_{\alpha,+}} \dd\zeta +  \int_{\Gamma_{\alpha,-}} \dd\zeta = 
 \sum_{k=1,2,3} \int_{\Gamma_{\alpha,+,k}} \dd \zeta  +\sum_{k=1,2,3} \int_{\Gamma_{\alpha,-,k}} \dd \zeta  \,,
\end{align*}
where, with the notation of $\tilde \alpha = \alpha l$, 
\begin{align}\label{proof.thm.semigroup.Kol.5}
\begin{split}
\Gamma_{\alpha,+,1} & = \big \{ \, \zeta = s + i \big ( \frac{1}{|\tilde \alpha|} \frac{s^3}{c^3} + |\tilde \alpha| \big )~|~ \quad  -c  |\tilde \alpha|^\frac23 \leq s\leq -c |\tilde \alpha|^\frac12 \big \}\,,\\
\Gamma_{\alpha,+,2} & = \big \{ \, \zeta = -c |\tilde \alpha|^\frac12 + i s ~|~ \quad |\tilde \alpha | - |\tilde \alpha|^\frac12 \leq s \leq |\tilde \alpha| +  |\tilde \alpha|^\frac12 \big \}\,,\\
\Gamma_{\alpha,+,3} & = \big \{ \, \zeta =  s + i (- \frac{s}{c}  + |\tilde \alpha| ) ~|~ \quad s\leq - c|\tilde \alpha|^\frac12\big \}\,,
\end{split}
\end{align}
and each $\Gamma_{\alpha,-,k}$ is a reflection of $\Gamma_{\alpha,+,k}$ with respect to the real axis.
Here $c$ is a positive constant which is independent of $|\alpha| \gg 1$ and also of $l$.
Let us estimate each integral by applying Theorem \ref{thm.Lambda.Kol.rate.2} and \eqref{proof.thm.semigroup.-1}.
As for the integral on $\Gamma_{\alpha,+,1}$ we see  
\begin{align}\label{proof.thm.semigroup.-2}
\| \frac{1}{2\pi} \int_{\Gamma_{\alpha,+,1}} e^{t\zeta} (\zeta - \mathbb{Q}_l L_{\alpha,l})^{-1}  \mathbb{Q}_l f_l \dd \zeta \|_{L^2} & \leq C \int_{-c|\tilde \alpha|^\frac23}^{-c|\tilde \alpha|^\frac12} e^{t s} \frac{\sqrt{1+\frac{9 s^4}{\tilde \alpha^2 c^6}}}{|\tilde \alpha|^\frac23 \big (1-\frac{1}{|\tilde \alpha|} (\frac{s^3}{|\tilde \alpha| c^3} + |\tilde \alpha|) \big )^\frac13} \dd s  \,  \| \mathbb{Q}_l f_l \|_{L^2} \nonumber \\
%& \leq C \int_{-c|\tilde \alpha|^\frac23}^{-c|\tilde \alpha|^\frac12} e^{t s} \frac{\sqrt{1+\frac{9 s^4}{\tilde \alpha^2 c^6}}}{|s|} \dd s  \,  \| \mathbb{Q}_l f \|_{L^2}\nonumber \\
& \leq \frac{C}{|\tilde \alpha|} \int_{-c |\tilde \alpha|^\frac23}^{-c|\tilde \alpha|^\frac12} |s| e^{t s} \dd s  \,  \| \mathbb{Q}_l f _l \|_{L^2} \nonumber \\
& \leq \frac{C}{|\tilde \alpha| t^2} e^{-c' |\tilde \alpha|^\frac12 t} \| \mathbb{Q}_l f_l \|_{L^2}\,,
\end{align}
where $c' = \frac{c}{2}$ and $\tilde \alpha = \alpha l$.
The integral on $\Gamma_{\alpha,+,2}$ is estimated as 
\begin{align}\label{proof.thm.semigroup.-3}
\| \frac{1}{2\pi} \int_{\Gamma_{\alpha,+,2}} e^{t\zeta} (\zeta - \mathbb{Q}_l L_{\alpha,l})^{-1}  \mathbb{Q}_l f _l \dd \zeta \|_{L^2} & \leq \frac{C}{ |\tilde \alpha|^\frac12} \int_{|\tilde \alpha |-|\tilde \alpha|^\frac12}^{|\tilde \alpha| +  |\tilde \alpha|^\frac12} e^{-c  |\tilde \alpha|^\frac12 t} \dd s  \,  \| \mathbb{Q}_l f_l \|_{L^2}  \leq  C  e^{-c|\tilde \alpha|^\frac12 t} \| \mathbb{Q}_l f \|_{L^2} \,.
\end{align}
As for the integral on $\Gamma_{\alpha,+,3}$, we have 
\begin{align}\label{proof.thm.semigroup.-4}
\| \frac{1}{2\pi} \int_{\Gamma_{\alpha,+,3}} e^{t\zeta} (\zeta - \mathbb{Q}_l L_{\alpha,l})^{-1}  \mathbb{Q}_l f_l \dd \zeta \|_{L^2} & \leq \frac{C}{|\tilde \alpha|} \int_{s\leq -c |\tilde \alpha|^\frac12} e^{ts} \frac{\sqrt{1+c^{-2}}}{\frac{1}{|\tilde \alpha|}(-\frac{s}{c} + |\tilde \alpha|) -1}\dd s  \,  \| \mathbb{Q}_l f_l \|_{L^2} \nonumber \\
& \leq C\int_{s\leq -c  |\tilde \alpha|^\frac12} e^{ts} \frac{\dd s}{|s|}  \,  \| \mathbb{Q}_l f_l \|_{L^2} \nonumber \\
& \leq \frac{C}{|\tilde \alpha|^\frac12 t} e^{-c |\tilde \alpha|^\frac12 t}  \| \mathbb{Q}_l f_l \|_{L^2}\,.
\end{align}
The estimates on the curve $\Gamma_{\alpha,-,k}$ are obtained in the same manner.
Hence, by rewriting $c'$ by $c$ for notational convenience, we have arrived at, from $\tilde \alpha = \alpha l$,
\begin{align}\label{proof.thm.semigroup.-5}
\| \mathbb{Q}_l e^{t L_{\alpha,l}} f_l \|_{L^2}  \leq
\begin{cases} 
& \displaystyle C \big ( \frac{1}{|\alpha l| t^2} + \frac{1}{|\alpha l|^\frac12 t} \big )  \, e^{-c |\alpha l|^\frac12 t} \| \mathbb{Q}_l f_l \|_{L^2}\,, \qquad t>0\\
& \displaystyle C e^{-c |\alpha l|^\frac12 t} \| \mathbb{Q}_l f_l \|_{L^2}\,, \qquad t\geq |\alpha l|^{-\frac12}\,.
\end{cases}
\end{align}
On the other hand, the simple energy computation for $\frac{\dd }{\dd t} \langle u_l (t), B_{2,l} u_l (t) \rangle_{L^2(\T)}$ gives the identity
\begin{align*}
\frac{\dd }{\dd t} \langle u_l , B_{2,l} u_l  \rangle_{L^2(\T)} = 2\langle A_l u_l, B_{2,l} u_l \rangle_{L^2(\T)}\,.
\end{align*}
The term in the right-hand side is bounded from above by $-(2l^2-1) )\| \mathbb{Q}_l u_l \|_{L^2(\T)}^2 \leq - (2l^2-1) \langle u_l,B_{2,l} u_l \rangle_{L^2 (\T)}$, where we have used $\langle u_l, B_{2,l} u_l \rangle_{L^2 (\T)} \leq \| \mathbb{Q}_l u \|_{L^2(\T)}^2$. 
Thus we have from the coercive estimate \eqref{coercive.B_2},
\begin{align}\label{proof.thm.semigroup.Kol.7}
\| \mathbb{Q}_l e^{t L_{\alpha,l}} f_l \|_{L^2} \leq 2 e^{-\frac12(2l^2-1) t} \| \mathbb{Q}_l f_l\|_{L^2} \,, \qquad t>0\,.
\end{align}
This estimate is useful for a short time period.
When $|l|\geq 2$ we have obtained the desired semigroup bound since $\mathbb{Q}_l=I$ in this case.
For the estimate of $\mathbb{P}_l e^{t L_{\alpha,l}}$ in the case $|l|=1$ we cannot shift the curve $\Gamma$ as in $\Gamma_\alpha$,
and it has to be computed in a different way.
By the construction of the resolvent in the proof of Theorem \ref{thm.abstract.1}, see \eqref{proof.thm.abstract.1'.3}, we observe that 
\begin{align}\label{proof.thm.semigroup.Kol.8}
\mathbb{P}_l (\zeta-L_{\alpha,l})^{-1} f_l & = - i \alpha l (\zeta-A_l)^{-1}\mathbb{P}_l \hat{\Lambda}_l (\zeta - \mathbb{Q}_l L_{\alpha,l})^{-1} \mathbb{Q}_l f_l + (\zeta-A_l)^{-1}\mathbb{P}_l f_l\,.
\end{align}
Since $\|(\zeta-A_l)^{-1} \|_{L^2 (\T) \rightarrow L^2 (\T)} \leq |\zeta|^{-1}$, we have 
\begin{align}
\| \mathbb{P}_l e^{t L_{\alpha,l}} f_l\|_{L^2} & \leq \alpha |l|\,  \| \frac{1}{2\pi i} \int_{\Gamma} e^{t\zeta} (\zeta-A_l)^{-1} \mathbb{P}_l \hat{\Lambda}_l (\zeta-\mathbb{Q}_l L_{\alpha,l})^{-1} \mathbb{Q}_l f_l \dd \zeta \|_{L^2}  \nonumber \\
& \qquad + \| \frac{1}{2\pi i } \int_{\Gamma} e^{t \zeta}  (\zeta-A_l)^{-1} \mathbb{P}_l f_l \dd \zeta \|_{L^2} \nonumber \\
& = \alpha I_1 + I_2\,,\nonumber 
\end{align} 
and it is not difficult to see 
\begin{align}\label{proof.thm.semigroup.Kol.9}
\| I_2 \|_{X}\leq e^{-t} \| \mathbb{P}_l f_l \|_{L^2}\,, \qquad t>0\,, \quad |l|=1
\end{align}
by the estimate for the semigroup generated by the self-adjoint operator $A_l$ in $L^2(\T)$.
As for $I_1$, we replace $\Gamma$ by $\tilde \Gamma_\alpha$, where
\begin{align*}
\tilde \Gamma_\alpha & = \tilde \Gamma_{-\frac12, \pm}   + \tilde \Gamma_{\alpha,\pm,1} + \Gamma_{\alpha,\pm,2} + \Gamma_{\alpha,\pm,3}\,,
\end{align*}
where $\Gamma_{\alpha,\pm.j}$ with $j=2,3$ are the curves as in \eqref{proof.thm.semigroup.Kol.5}, while $\tilde \Gamma_{-\frac12,+}$ (Rep. $\tilde \Gamma_{-\frac12,-}$)   is the segment connecting $\zeta=-\frac12$ and a point $p_{\alpha,+}$ of $\Gamma_{\alpha,+,1}$ (Resp. $p_{\alpha,-}$ of $\Gamma_{\alpha,-,1}$), and $\tilde \Gamma_{\alpha,\pm,1}$ is the part of $\Gamma_{\alpha,\pm,1}$ which therefore connects $p_{\alpha,\pm}$ and the end point of $\Gamma_{\alpha,\pm,2}$.
We will take $p_{\alpha,\pm}$ as $|\Im p_{\alpha,\pm}| =\frac{|\alpha l|}{2}$, thus, they are away enough from the degenerate case such as $|\Im \zeta| \sim |\alpha l| =|\alpha|$.  
On the curve $\tilde \Gamma_{-\frac12,\pm}$ we have $\| (\zeta -\mathbb{Q}_l L_{\alpha,l})^{-1} \|_{Y_l\rightarrow Y_l} \leq \frac{C}{|\alpha|^\frac23}$ by the choice of $p_{\alpha,\pm}$.
Thus we have 
\begin{align*}
\| \frac{1}{2\pi i} \int_{\tilde \Gamma_{-\frac12,\pm}} e^{t\zeta} (\zeta-A_l)^{-1} \mathbb{P}_l \hat{\Lambda}_l (\zeta-\mathbb{Q}_l L_{\alpha,l})^{-1} \mathbb{Q}_l f_l \dd \zeta \|_{L^2}  & \leq \frac{C}{|\alpha|^\frac23} e^{-\frac{t}{2}}  \int_{\zeta=-\frac12}^{\zeta=p_{\alpha,\pm}}  \frac{|\dd \zeta|}{|\zeta|} \\
& \leq \frac{C|\log \alpha|}{|\alpha|^\frac23} e^{-\frac{t}{2}} \|  \mathbb{Q}_l f_l \|_{L^2}\,.
\end{align*}
As for the integrals on the curve $\tilde \Gamma_{\alpha,\pm,1}$ and $\Gamma_{\alpha,\pm,2}$ we compute as in \eqref{proof.thm.semigroup.-2} and \eqref{proof.thm.semigroup.-3} respectively, but the difference in this case is the presence of the factor $(\zeta- A_l)^{-1} \mathbb{P}_l \hat{\Lambda}_l$, which is bounded by 
\begin{align*}
\| (\zeta- A_l)^{-1} \mathbb{P}_l \hat{\Lambda}_l \|_{L^2\rightarrow L^2} \leq \frac{C}{|\zeta|} \leq \frac{C}{|\alpha|}
\end{align*}
on $\tilde \Gamma_{\alpha,\pm,1}$ and $\Gamma_{\alpha,\pm,2}$.
Moreover, we need to compute the integral so that the singularity at $t=0$ does not appear. 
Hence we have 
\begin{align*}
\| \frac{1}{2\pi i} \int_{\tilde \Gamma_{\alpha,\pm,1}} e^{t\zeta} (\zeta-A_l)^{-1} \mathbb{P}_l \hat{\Lambda}_l (\zeta-\mathbb{Q}_l L_{\alpha,l})^{-1} \mathbb{Q}_l f_l \dd \zeta \|_{L^2}  & \leq \frac{C}{|\alpha|^2} \int_{-c |\alpha|^\frac23}^{-c|\alpha|^\frac12} |s| e^{t s} \dd s  \,  \| \mathbb{Q}_l f _l \|_{L^2} \nonumber \\
& \leq \frac{C}{|\alpha|^\frac23} e^{-c|\alpha|^\frac12 t} \| \mathbb{Q}_l f _l \|_{L^2}\,,
\end{align*}
and similarly,
\begin{align*}
\| \frac{1}{2\pi i} \int_{\Gamma_{\alpha,\pm,2}} e^{t\zeta} (\zeta-A_l)^{-1} \mathbb{P}_l \hat{\Lambda}_l (\zeta-\mathbb{Q}_l L_{\alpha,l})^{-1} \mathbb{Q}_l f_l \dd \zeta \|_{L^2}  & \leq \frac{C}{|\alpha|^\frac32} \int_{|\alpha |-|\alpha|^\frac12}^{|\alpha| + |\alpha|^\frac12} e^{-c|\alpha|^\frac12t} \dd s  \,  \| \mathbb{Q}_l f _l \|_{L^2} \nonumber \\
& \leq \frac{C}{|\alpha|} e^{-c|\alpha|^\frac12 t} \| \mathbb{Q}_l f _l \|_{L^2}\,.
\end{align*}
Finally, the integrals over $\Gamma_{\alpha,\pm,3}$ is computed as in \eqref{proof.thm.semigroup.-4}, and we have 
\begin{align*}
& \| \frac{1}{2\pi} \int_{\Gamma_{\alpha,+,3}} e^{t\zeta} (\zeta-A_l)^{-1} \mathbb{P}_l \hat{\Lambda}_l (\zeta - \mathbb{Q}_l L_{\alpha,l})^{-1}  \mathbb{Q}_l f_l \dd \zeta \|_{L^2} \\
& \qquad \leq \frac{C}{|\alpha l|} \int_{s\leq -c  |\alpha|^\frac12} e^{ts} \frac{\sqrt{1+c^{-2}}}{(|s| +|\alpha|) \big( \frac{1}{|\alpha l|}(-\frac{s}{c} + |\alpha l|) -1 \big )}\dd s  \,  \| \mathbb{Q}_l f_l \|_{L^2} \nonumber \\
& \qquad \leq C\int_{s\leq -c  |\alpha|^\frac12} e^{ts} \frac{1}{(|s|+|\alpha|)|s|} \dd s   \,  \| \mathbb{Q}_l f_l \|_{L^2} \nonumber \\
& \qquad \leq \frac{C|\log \alpha|}{|\alpha|} e^{-c |\alpha|^\frac12 t}  \| \mathbb{Q}_l f_l \|_{L^2}\,.
\end{align*}
Collecting these above,  we have 
\begin{align}\label{proof.thm.semigroup.Kol.10}
|\alpha| I_1 & \leq C|\alpha|^\frac13 |\log \alpha| \,  e^{-\frac{t}{2}} \| \mathbb{Q}_l f_l  \|_{L^2} \,.
\end{align}
Combining \eqref{proof.thm.semigroup.Kol.9} and \eqref{proof.thm.semigroup.Kol.10}, we obtain 
\begin{align}\label{proof.thm.semigroup.Kol.11}
\| \mathbb{P}_l e^{t L_{\alpha,l}} f_l \|_{L^2} \leq C  |\alpha|^\frac13 |\log \alpha |  \, e^{-\frac{t}{2}} \| \mathbb{Q}_l  f_l  \|_{L^2} + e^{-t} \| \mathbb{P}_l f_l \|_{L^2} \,, \qquad t>0\,, \quad |l|=1\,.
\end{align}

Thus we have arrived at the following theorem.
Let us recall that $\mathbb{Q}:L^2_0 (\T^2)\rightarrow Y$ is the orthogonal projection on to $Y$.
\begin{thm}\label{thm.semigroup.Kol} For all sufficiently large $|\alpha|$ the following statement holds.
The semigroup $\{e^{tL_\alpha}\}_{t\geq 0}$ generated by $L_\alpha$ in $L^2_0 (\T^2)$ satisfies for any $l\in \Z\setminus\{0\}$,
\begin{align}\label{est.semigroup.Kol.1}
\| \mathbb{Q} \mathcal{P}_l e^{t L_\alpha} f \|_{L^2(\T^2)} \leq 
\begin{cases}
& \displaystyle 2e^{-\frac12(2l^2-1)t} \| \mathbb{Q} \mathcal{P}_l f \|_{L^2 (\T^2)} \,, \qquad \qquad t>0 \,,\\
& \displaystyle C e^{-c |\alpha l|^\frac12  t} \| \mathbb{Q} \mathcal{P}_l f\|_{L^2(\T^2)}\,, \qquad \qquad t\geq |\alpha l|^{-\frac12}\,,
\end{cases} 
\end{align}
while
\begin{align}\label{est.semigroup.Kol.2}
\begin{split}
\| (I-\mathbb{Q}) \mathcal{P}_l e^{t L_\alpha} f\|_{L^2(\T^2)}  & \leq C  |\alpha|^\frac13 |\log \alpha |  \, e^{-\frac{t}{2}} \| \mathbb{Q} \mathcal{P}_l f \|_{L^2(\T^2)}\\
& \qquad  + e^{-t} \| (I-\mathbb{Q}) \mathcal{P}_l f \|_{L^2 (\T^2)} \,, \qquad t>0\,, \quad |l|=1\,.
\end{split}
\end{align}
Here $C$ and $c$ are independent of $\alpha$, $l$, and $f$.
\end{thm}

\begin{rem}{\rm Recently Wei \cite{W} obtained a refined version of the Gearhart-Pr${\rm \ddot{u}}$ss  theorem for semigroup  in terms of the pseudospectral bound. If one applies this general result of \cite{W} then  the semigroup estimate \eqref{est.semigroup.Kol.1} is a direct consequence of the pseudospecral bound $\displaystyle \sup_{\lambda\in \R} \| (i\lambda + \mathbb{Q}_l L_{\alpha,l})^{-1} \|_{Y_l\rightarrow Y_l}\leq \frac{C}{|\alpha l|^\frac12}$ in Theorem \ref{thm.Lambda.Kol.rate.2}, and thus the proof of Theorem \ref{thm.semigroup.Kol} is much shortened.
}
\end{rem}

Theorem \ref{thm.semigroup.Kol} immediately leads to the estimate of the solution to 
\begin{align}\label{proof.thm.semigroup.1}
\left\{
\begin{aligned}
& \frac{\dd  \omega}{\dd t} - \mathcal{L}^{\nu,a} \omega  =0\,, \qquad t>0\,, \\
& \omega |_{t=0}   = f\in L^2_0 (\T^2)\,.
\end{aligned}\right.
\end{align}
Here $\mathcal{L}^{\nu,a}$ is defined as in \eqref{def.L^nu'} with $a\in \R\setminus\{0\}$ and $0<\nu\ll 1$.
Indeed, by introducing the rescaling $\omega (x,y,t) = w(x,y,\nu t)$ which gives $\alpha=\frac{a}{\nu}$ in Theorem \ref{thm.semigroup.Kol}, we obtain
\begin{cor}\label{cor.thm.semigroup.Kol}
For all sufficiently small $\frac{\nu}{a}>0$ the following statement holds.
The semigroup $\{e^{t \mathcal{L}^{\nu,a}}\}_{t\geq 0}$ generated by $\mathcal{L}^{\nu,a}$ in $L^2_0 (\T^2)$ satisfies for any $l\in \Z\setminus\{0\}$,
\begin{align}\label{est.cor.semigroup.Kol.1}
\| \mathbb{Q} \mathcal{P}_l e^{t \mathcal{L}^{\nu,a}} f \|_{L^2(\T^2)} \leq 
\begin{cases}
& \displaystyle 2 e^{-\frac12(2l^2-1)\nu t} \| \mathbb{Q} \mathcal{P}_l f \|_{L^2 (\T^2)} \,, \qquad \qquad t>0 \,,\\
& \displaystyle C  e^{-c  \sqrt{a |l|\nu}\,  t} \| \mathbb{Q} \mathcal{P}_l f\|_{L^2(\T^2)}\,, \qquad \qquad t\geq \frac{1}{\sqrt{a |l|\nu}}\,,
\end{cases} 
\end{align}
while
\begin{align}\label{est.cor.semigroup.Kol.2}
\begin{split}
\| (I-\mathbb{Q}) \mathcal{P}_l e^{t \mathcal{L}^{\nu,a}} f\|_{L^2(\T^2)}   \leq C  (\frac{a}{\nu})^\frac13 |\log \frac{a}{\nu} |  \, e^{-\frac{\nu t}{2}} \| \mathbb{Q} \mathcal{P}_l f \|_{L^2(\T^2)} & + e^{-\nu t} \| (I-\mathbb{Q}) \mathcal{P}_l f \|_{L^2 (\T^2)} \,,\\
& \qquad  \qquad  t>0\,, \quad |l|=1\,.
\end{split}
\end{align}
Here $C$ and $c$ are independent of $\nu$, $a$, $l$, and $f$.
\end{cor}

\section{Application to the Lamb-Oseen vortices}

In this section we consider the operator related to the Lamb-Oseen vortices.
For details of the derivation of the operators below, the reader is referred to Deng \cite{De2} and Li, Wei, and Zhang \cite{LWZ}. 
Let  $\displaystyle H^1_{1,0}(\R_+) = \{f\in H^1_0 (\R_+)~|~\frac{f(r)}{r}\,, r f(r) \in L^2(\R_+)\}$.
Let $A$  be the realization in $L^2(\R_+)$ of 
\begin{align}\label{def.A.O}
A = \partial_r^2 - \frac{3}{4r^2} - \frac{r^2}{16} + \frac12
\end{align}
with the maximal domain $D(A) = \{f\in H^1_{1,0} (\R_+)~|~A f\in L^2(\R_+)\}$.
It is known that the operator $A$ is conjugate to the two-dimensional Harmonic oscillator restricted to the angular Fourier mode $\pm 1$, and hence, $A$ is a self-adjoint operator in $L^2(\R_+)$ with compact resolvent and $-A$ is strictly positive. 
Moreover, we have the equivalence between the norms such as 
\begin{align}
\| (-A)^\frac12 u \|_{L^2} \sim \| \partial_r u\|_{L^2} + \| \frac{u}{r}\|_{L^2} + \| r u\|_{L^2}\,.
\end{align}
This fact will be frequently used in the analysis below.  
Let $g(r) = e^{-\frac{r^2}{8}}$ and we denote by $Y$ the orthogonal complement space in $L^2(\R_+)$ to the one-dimensional subspace spanned by $r^\frac32 g(r)$, that is, 
\begin{align}
Y=\{r^\frac32 g(r)\}^\bot \qquad {\rm in}~~L^2(\R_+)\,.
\end{align}
Then $Y$ is known to be invariant under the action of $A$.
Let $\rho (r)$ be the function defined by 
\begin{align}\label{def.sigma.O}
\rho (r) = \frac{1-e^{-\frac{r^2}{4}}}{r^2/4}\,.
\end{align}
The direct computation shows that $\rho(0)=1$, $\rho'(0)=0$, $\rho''(0)<0$, and $\rho'(r)<0$ for $r>0$. 
We also introduce the nonlocal operator $Z$ as 
\begin{align}\label{def.Z.O}
Z [f] (r) = r^{-\frac32} \int_0^r f(s) g (s) s^\frac32 \, \dd s\,.
\end{align}
Then $Z$ belongs to the Hilbert-Schmidt class, for the kernel $k(r,s) = r^{-\frac32} g(s) s^\frac32 \chi_{\{0<s<r\}}$ belongs to $L^2(\R_+\times \R_+)$. We define the self-adjoint operators $\hat{\Lambda}_1$, $\hat{\Lambda}_2$, and $\hat{\Lambda}$,  as 
\begin{align}\label{def.Lambda.O}
\hat{\Lambda}_1=M_{\rho}\,, \qquad \hat{\Lambda}_2 = - Z^* Z \,, \qquad \hat{\Lambda} = \hat{\Lambda}_1 + \hat{\Lambda}_2\,,
\end{align}
where $Z^*$ is the adjoint of $Z$ in $L^2(\R_+)$, that is,
\begin{align}\label{def.Z^*.O}
Z^*[f] (r) =  g(r) r^\frac32 \int_r^\infty f(s) s^{-\frac32} \, \dd s\,.
\end{align}
The explicit representation of $\hat{\Lambda}_2$ is given by 
\begin{align*}
\hat{\Lambda}_2 [f] (r) = - \frac{g(r)}{2} \int_0^\infty f (s)  \min\{ \frac{r}{s},\frac{s}{r}\} (rs)^\frac12 g(s) \, \dd s\,.
\end{align*}
We are interested in the resolvent estimate of the operator in $Y$ defined as 
\begin{align*}
L_\alpha = A - i \alpha \hat{\Lambda}\,, \qquad D_Y(L_\alpha)=D_Y(A)=D(A) \cap Y\,.
\end{align*}
We note that, by Proposition \ref{prop.spectral.Lambda.O} (1) below and by the fact that $Y$ is invariant under the action of $A$, the operator $L_\alpha$ is indeed well-defined as the operator from $D(A) \cap Y$ to $Y$.  As observed in \cite{LWZ}, the operator $L_\alpha$ is conjugate to the linearized operator around the Lamb-Oseen vortex with the total circulation $\alpha$ in the self-similar variable (that is, $\mathcal{L}-\alpha \Lambda$ introduced as in \eqref{LambO1}), when it is restricted to the angular Fourier mode $\pm 1$ (and it is known that this case is the most difficult to analyze). 
It will be convenient to see that the inverse of $Z$, denoted  by $W$, is given by 
\begin{align}\label{def.W.O}
W[v] (r) = \frac{1}{g (r) r^\frac32}\big (r^\frac32 v (r) \big )'\,,
\end{align}
and the domain of $W$ is $D(W)=\{v\in H^1_{1,0}(\R_+)~|~W[v]\in L^2 (\R_+)\}$.
The operator $W$ is densely defined and closed in $L^2(\R_+)$.
Let $W^*$ be the adjoint of $W$ in $L^2(\R_+)$. Then the direct computation shows 
\begin{align}\label{def.WW^*.O}
W W^* f =  - \frac{1}{g^2} \Big ( f''  +\frac{r}{2} f' + (\frac{r^2}{16} + \frac14 - \frac{3}{4r^2} ) f \Big ) = -\frac{1}{g^2} \Big ( A f + \frac{r}{2} f' + (\frac{r^2}{8} -\frac14 ) f \Big ) \,. 
\end{align}
Note that $- W W^* \hat{\Lambda}_2=I$.
We have the following estimates for $\hat{\Lambda}_2$:
\begin{align}\label{est.Z^*}
\| \frac{1}{(1+r) g} \partial_r Z^* f \|_{L^2} +\| \frac{1}{rg} Z^* f\|_{L^2} \leq C \| f\|_{L^2}\,, \qquad | \Im \langle A f,  \hat{\Lambda}_2 f \rangle_{L^2} | \leq C  \| M_{\rho'}  f\|_{L^2} \| Z f\|_{L^2}\,.
\end{align} 
The second inequality in \eqref{est.Z^*} follows from $A^*=A$ and 
$$A\hat{\Lambda}_2 f =  -g^2 W W^* \hat{\Lambda}_2 f - \frac{r}{2} \frac{\dd}{\dd r} \hat{\Lambda}_2 f - (\frac{r^2}{8}-\frac14) \hat{\Lambda}_2 f  = g^2 f - \frac{r}{2} \frac{\dd}{\dd r} \hat{\Lambda}_2 f - (\frac{r^2}{8}-\frac14) \hat{\Lambda}_2 f \,,$$
where \eqref{def.WW^*.O} is used.
Indeed, this identity implies 
\begin{align*}
|\Im \langle A f, \hat{\Lambda}_2 f\rangle_{L^2}| = |\Im \langle f, \frac{r}{2}\partial_r Z^* Z f + \frac{r^2}{8} Z^* Z f \rangle_{L^2} | & \leq \| M_{\rho'} f\|_{L^2} \big ( \| M_{\frac{r}{2\rho'}}\partial_r Z^* Z f\|_{L^2} + \| M_{\frac{r^2}{8\rho'}} Z^* Z f\|_{L^2}\big ) \\
& \leq C \| M_{\rho'} f\|_{L^2} \| Z f\|_{L^2}\,. 
\end{align*}
Here we have used the first inequality in \eqref{est.Z^*} in the last line.
As for $\hat{\Lambda}$, the following properties are known.
\begin{prop}\label{prop.spectral.Lambda.O} {\rm (1)} Each $\hat{\Lambda}_j$ is bounded and self-adjoint, and $\hat{\Lambda}_2$ is compact in $L^2(\R_+)$. 

\noindent {\rm (2)} $\sigma (\hat{\Lambda})=[0,1]$, where $\sigma(\hat{\Lambda})$ is the spectrum of $\hat{\Lambda}$ in $L^2(\R_+)$. Moreover, $\{r^\frac32 g(r)\}$ is the eigenspace for the eigenvalue $0$ of $\hat{\Lambda}$ in $L^2(\R_+)$, and $\hat{\Lambda}$ does not possess the eigenvalues except for the eigenvalue $0$.
\end{prop}

Indeed, the statement (1) is clear from the definition of $\hat{\Lambda}_j$ above (the fact $\hat{\Lambda}_2$ is self-adjoint is firstly observed by Gallay-Wayne \cite{GW} in the original formulation of the linearized operator around the Burgers vortex). As for the statement (2), the fact $\sigma(\hat{\Lambda}_1)={\rm Ran}\, (\rho)=[0,1]$ with no eigenvalues is trivial from the definition of $\hat{\Lambda}_1$, and then, since $\hat{\Lambda}_2$ is compact, the essential spectrum of $\hat{\Lambda}$ must coincides with that of $\hat{\Lambda}_1$, i.e., it is the interval $[0,1]$. The structure of the eigenvalues of $\hat{\Lambda}$ follows from \cite{Ma}, which also gives $\sigma (\hat{\Lambda})=[0,1]$.

We also have  
\begin{prop}\label{prop.bound.Lambda.O} $|\Im \langle A f, \hat{\Lambda} f\rangle_{L^2}|  \leq C \| (-A)^\frac12 f \|_{L^2} \| M_{\rho'} f \|_{L^2}$\, for any $f\in D(A)$.
\end{prop}

\noindent {\it Proof.} In virtue of \eqref{est.Z^*} for the estimate of $A\hat{\Lambda}_2$ it suffices to show $$|\Im \langle A f , \hat{\Lambda}_1 f \rangle_{L^2}|\leq  C\| (-A)^\frac12 f \|_{L^2} \| M_{\rho'} f \|_{L^2}\,.$$
To see this, we observe that $A f = (- T^* T + \frac12 ) f$, where $T=\partial_r + \frac{1}{2r} + \frac{r}{4}$ and $T^* = - \partial_r +\frac{1}{2r} + \frac{r}{4}$. Thus we have
\begin{align*}
\Im \langle A f, \hat{\Lambda}_1 f\rangle_{L^2} = - \Im \langle T f, T M_{\rho} f \rangle_{L^2} = - \Im \langle T f, [T, M_{\rho}] f\rangle_{L^2} = - \Im \langle T f, M_{\rho'} f\rangle_{L^2}\,.
\end{align*}   
Then the desired estimate follows from $\|(-A)^\frac12 f\|_{L^2}^2 = \langle -A f , f \rangle_{L^2} = \| T f\|_{L^2}^2 - \frac12 \| f\|_{L^2}^2$ and $\| f\|_{L^2} \leq C \| (-A)^\frac12 f\|_{L^2}$. The proof is complete.

\

Let $h_1(m,\mu)$, $m>0$ and $\mu\in [-\frac12,\frac12]$ be the nonnegative function whose square is defined as
\begin{align}\label{def.h_1.O.1}
h_1 (m,\mu)^2  = 
\begin{cases}
& \displaystyle \frac{1}{m^2\mu^3} \qquad m> \frac{1}{\mu^2}~~{\rm and}~~0<\mu\leq \frac12\,,\\
& \displaystyle \mu \qquad \quad \frac{1}{10\mu} < m \leq \frac{1}{\mu^2}~~{\rm and}~~0<\mu \leq \frac12\,,\\
& \displaystyle \frac{1}{m} \quad \qquad 0 < m\leq \frac{1}{10|\mu|}\,.
\end{cases}
\end{align}
Note that $h_1(m,\mu)^2=\frac1m$ for $-\frac12\leq \mu\leq 0$.
\begin{lem}\label{lem.Lambda.O1} Let $-\frac12 \leq \mu\leq \frac12$. Then there exists $C>0$ such that for $\delta\in (0,1]$ and $u\in D(A)\cap Y$, 
\begin{align}\label{est.lem.Lambda.O1.1}
\delta^2 \| M_{\rho'} u \|_{L^2}^2 + \| Z u\|_{L^2}^2 + \| Z^* Z u \|_{L^2}^2 \leq C \delta^{-2} \| (\mu -\hat{\Lambda}) u \|_{L^2}^2 + C \delta^6 |\mu_+|^3 \| (-A)^\frac12 u \|_{L^2}^2\,,
\end{align}
and 
\begin{align}\label{est.lem.Lambda.O1.2}
\| u\|_{L^2}^2 \leq C \delta^{-2} \| (\mu -\hat{\Lambda}) u \|_{L^2}^2 + C h_1 (\frac{1}{\delta}, \mu)^2 \| (-A)^\frac12 u \|_{L^2}^2\,.
\end{align}
Here $\mu_+=\max\{ \mu,0\}$. 
\end{lem}

\noindent  {\it Proof of \eqref{est.lem.Lambda.O1.1}.} Again we will use the contradiction argument.
Suppose that \eqref{est.lem.Lambda.O1.1} does not hold. Then there exist $\{\delta_n, \mu_n\}_{n\in \N}$, $\delta_n \in (0,1]$, $\mu_n \in (-\frac12,\frac12]$, and $\{u_n\}\subset D(A)$ such that 
\begin{align*}
& \lim_{n\rightarrow \infty} \delta_n=\delta_\infty\in [0,1]\,, \qquad \lim_{n\rightarrow \infty}\mu_n=\mu_\infty \in [-\frac12,\frac12]\,,
\end{align*}
and 
\begin{align}\label{proof.lem.Lambda.O1.1}
\begin{split}
& \delta_n^2 \| M_{\rho'} u_n\|_{L^2}^2+ \| Z u_n \|_{L^2}^2 + \| Z^* Z u_n \|_{L^2}^2  =1\,, \\
& \lim_{n\rightarrow \infty}  \bigg ( \delta_n^{-2}  \|  (\mu_n - \hat{\Lambda}) u_n \|_{L^2}^2 + \delta_n^6 |(\mu_n)_+|^3 \| (-A)^\frac12  u_n\|_{L^2}^2 \bigg ) =0\,.
\end{split}
\end{align}
Set
\begin{align}\label{proof.lem.Lambda.O1.2}
f_n = \delta_n^{-1} (\mu_n   - \hat{\Lambda}) u_n\,, \qquad v_n = Z^* Z u_n\,,
\end{align}
and then, $u_n=W W^* v_n$.
Since $Z^*$ is compact, by \eqref{est.Z^*}, we may assume that $v_n$ converges to $v_\infty$ strongly in $L^2(\R_+)$, $W^*v_n$ converges to $W^*v_\infty$ weakly in $L^2(\R_+)$.
Since $u_n$ is real valued, so is $v_n$, and $v_n$ satisfies $v_n(0)=0$ and 
\begin{align}\label{proof.lem.Lambda.O1.3}
(\mu_n - M_{\rho} )W W^* v_n  = v_n + \delta_n f_n\,, \qquad r>0 \,.
\end{align}
Let $r_n\in [\frac{1}{100},\infty]$ be the critical point, i.e., $r_n=\rho^{-1} (\mu_n)$ for $\mu_n\in (0,\frac12]$ and $r_n=\infty$ for $\mu_n\in [-\frac12,0]$. Then from \eqref{proof.lem.Lambda.O1.3} we have 
\begin{align}\label{proof.lem.Lambda.O1.4}
v_n (r_n) + \delta_n f_n (r_n) =0\,,
\end{align}
which is valid also for the case $r_n=\infty$ by setting $v_n(\infty)=f_n (\infty)=0$.
The following estimates will be frequently used:
\begin{align}
|\rho^{(k)}(r)| & \leq C (1+r)^{-2-k}\,, \qquad \qquad k=0,1,2\,, \label{proof.lem.Lambda.O1.51}\\
|\frac{1}{\mu_n-\rho (r)}| & \leq \frac{Cr_n(1+r \min\{r,r_n\})}{|r_n-r|}\,, \qquad  r\ne r_n\,, \label{proof.lem.Lambda.O1.52}\\
r_n & = \rho^{-1} (\mu_n) \approx 2\mu_n^{-\frac12} \quad \quad {\rm for}~~ 0<\mu_n\ll 1\,.\label{proof.lem.Lambda.O1.53}
\end{align}
When $r_n=\infty$ estimate \eqref{proof.lem.Lambda.O1.52} is interpreted as $C(1+r^2)$.
We define the operation $\infty \cdot 0=0\cdot \infty =0$ and set $\frac{1}{\rho'(\infty)}=\infty$ below.

\noindent {\bf Step 1: $\displaystyle \lim_{n\rightarrow \infty} \frac{v_n (r_n)}{\rho'(r_n)}=  \displaystyle \lim_{n\rightarrow \infty} \frac{\delta_n f_n (r_n)}{\rho'(r_n)}=0$\,.}

\noindent Since $\frac{v_n(\infty)}{\rho'(\infty)} = \infty\cdot 0=0$ and $\frac{\delta_n f_n (\infty)}{\rho'(\infty)}=\infty\cdot 0=0$, it suffices to consider the case $0<r_n<\infty$ for all $n$. 
In this case $\mu_n>0$ for all $n$. 
First we consider the case $\mu_\infty=0$, that is, $\displaystyle \lim_{n\rightarrow \infty} r_n=\infty$. 
In this case we compute as 
\begin{align*}
|v_n (r_n)|^2 \leq 2 \int_{r_n}^\infty |v_n' v_n | \, \dd r \leq 2 \| v_n' \|_{L^2} \| v_n\|_{L^2([r_n,\infty))} & \leq \frac{C}{(1+r_n)^8} \| Z u_n \|_{L^2} \| (1+r)^8 v_n \|_{L^2}\\
& \leq \frac{1}{C(1+r_n)^8}   \| Z u_n \|_{L^2}^2\,.
\end{align*}
Hence, since $\rho'(r) = -2r^{-3} (1+ o(1))$ for $r\gg 1$, we have $\displaystyle \lim_{n\rightarrow \infty} |\frac{v_n (r_n)}{\rho'(r_n)}| \leq \lim_{n\rightarrow \infty} \frac{C}{1+r_n}=0$.
Then we also have $\displaystyle \lim_{n\rightarrow \infty} \frac{\delta_n f_n (r_n)}{\rho'(r_n)}=0$ by using \eqref{proof.lem.Lambda.O1.4}. Next we consider the case $\mu_\infty\in (0,\frac12]$, that is, $\sup_n r_n<\infty$ and $\sup_n \frac{1}{|\rho'(r_n)|}<\infty$. 
Let $\kappa_1>0$ be fixed but arbitrary small number. 
We find $\tilde r_n\in [r_n, r_n + \kappa_1^2 \delta_n^2]$ such that $\kappa_1^2 \delta_n^2 |f_n (\tilde r_n)|^2 \leq \| f_n\|_{L^2([r_n,r_n+\kappa_1^2 \delta_n^2]}^2 \leq \| f_n \|_{L^2}^2$, and thus, $\delta_n |f_n (\tilde r_n)|\leq \frac{1}{\kappa_1} \| f_n \|_{L^2}$. On the other hand, we have 
\begin{align*}
|\delta_n f_n (r_n) -\delta_n f_n (\tilde r_n)| & \leq C \kappa_1 \delta_n \| \partial_r \delta_n f_n \|_{L^2([r_n,\tilde r_n])}
\end{align*}
and  \eqref{proof.lem.Lambda.O1.3} yields
\begin{align}\label{proof.lem.Lambda.O1.6}
\| \partial_r \delta_n f_n \|_{L^2([r_n,\tilde r_n])} & \leq C \big ( \| (\mu_n-M_\rho) \partial_r u_n \|_{L^2([r_n,\tilde r_n])} + \| M_{\rho'} u_n \|_{L^2} + \| \partial_r v_n \|_{L^2} \big ) \nonumber  \\
& \leq C \kappa_1^2 \delta_n^2 \| \partial_r u_n \|_{L^2} + C + C\delta_n\,, 
\end{align}
since $|\mu_n-\rho'(r) |\leq C \kappa_1^2\delta_n^2$ for $r\in [r_n,r_n+\kappa_1^2 \delta_n^2]$ and the normalized condition \eqref{proof.lem.Lambda.O1.1}. Then \eqref{proof.lem.Lambda.O1.1} implies 
\begin{align*}
\limsup_{n\rightarrow \infty} |\delta_n f_n (r_n)| & \leq 
\limsup_{n\rightarrow \infty} (|\delta_n f_n (r_n) - \delta_n f_n (\tilde r_n)| + |\delta_n f_n (\tilde r_n) |) \\
& \leq \limsup_{n\rightarrow \infty} \Big ( \frac{\|f_n\|_{L^2}}{\kappa_1} + C\delta_n^3 \| (-A)^\frac12 u_n \|_{L^2} + C \kappa_1 \Big ) \leq C \kappa_1\,.
\end{align*}
Since $\kappa_1>0$ is arbitrary, we have $\displaystyle \lim_{n\rightarrow \infty} \delta_n f_n (r_n)=0$, which also gives $\displaystyle \lim_{n\rightarrow \infty} v_n (r_n)=0$ in virtue of \eqref{proof.lem.Lambda.O1.4}.
This proves the claim in the case $\sup_n r_n<\infty$.

\

\noindent {\bf Step 2: $\displaystyle \limsup_{n\rightarrow \infty} \delta_n \| M_{\rho'} u_n \|_{L^2}\leq C \limsup_{n\rightarrow \infty} \delta_n \| Z^* Z u_n \|_{L^2}$\,.}

\noindent Let $\kappa_2>0$ be fixed but arbitrary small number. Let $r_n<\infty$ (thus, $\mu_n>0$).
Then, by setting $u_n(r)=0$ if $r<0$,
\begin{align*}
\int_{r_n-\frac{\delta_n^2}{\kappa_2^2}}^{r_n+\frac{\delta_n^2}{\kappa_2}} |M_{\rho'} u_n|^2\, \dd r \leq \frac{C\delta_n^2}{\kappa_2^2} \| M_{\rho'} u_n \|_{L^\infty ([r_n-\frac{\delta_n^2}{\kappa_2^2},\infty))}^2 &\leq \frac{C\delta_n^2}{\kappa_2} \| \partial_r (M_{\rho'} u_n ) \|_{L^2 ([r_n-\frac{\delta_n^2}{\kappa_2^2},\infty))} \| M_{\rho'} u_n \|_{L^2 ([r_n-\frac{\delta_n^2}{\kappa_2^2},\infty))} \\
& \leq \frac{C\delta_n}{\kappa_2} \| M_{\rho'} \partial_r u_n + M_{\rho''} u_n \|_{L^2([r_n-\frac{\delta_n^2}{\kappa_2^2},\infty))} \\
& \leq \frac{C\delta_n}{\kappa_2 (1+r_n)^3} \| (-A)^\frac12 u_n \|_{L^2} \leq C \delta_n \mu_n^\frac32 \| (-A)^\frac12 u_n \|_{L^2}\,. 
\end{align*}
Here we have used $\|\rho'\|_{L^\infty ([r_n-\frac{\delta_n^2}{\kappa_2^2},\infty))} + \| \rho'' \|_{L^\infty([r_n-\frac{\delta_n^2}{\kappa_2^2},\infty))} \leq C (1+r_n)^{-3}$. Note that the above estimate is valid also for the case $\mu_n\in (-\frac12,0]$, for $r_n=\infty$ in this case.
On the other hand, we have from \eqref{proof.lem.Lambda.O1.3}, when $r_n<\infty$,
\begin{align*}
\| M_{\rho'} u_n\|_{L^2([r_n-\frac{\delta_n^2}{\kappa_2^2}, r_n+\frac{\delta_n^2}{\kappa_2^2}]^c)}
& \leq \| \frac{\rho'}{\mu_n-\rho} v_n \|_{L^2([r_n-\frac{\delta_n^2}{\kappa_2^2}, r_n+\frac{\delta_n^2}{\kappa_2^2}]^c)} + \| \frac{\rho'}{\mu_n-\rho} \delta_n f_n \|_{L^2([r_n-\frac{\delta_n^2}{\kappa_2^2}, r_n+\frac{\delta_n^2}{\kappa_2^2}]^c)}\\
&  \leq \| \frac{\rho'}{\mu_n-\rho}  \|_{L^2([r_n-\frac{\delta_n^2}{\kappa_2^2}, r_n+\frac{\delta_n^2}{\kappa_2^2}]^c)} \|v_n \|_{L^\infty} +  \| \frac{\rho'}{\mu_n-\rho}  \|_{L^\infty ([r_n-\frac{\delta_n^2}{\kappa_2^2}, r_n+\frac{\delta_n^2}{\kappa_2^2}]^c)}  \| \delta_n f_n \|_{L^2} \\
& \leq \frac{C\kappa_2}{\delta_n} \| v_n \|_{L^\infty} + \frac{C\kappa_2^2}{\delta_n} \| f_n \|_{L^2}\,.
\end{align*}
When $r_n=\infty$ (and thus, $\mu_n\in [-\frac12,0]$), we see $[r_n-\frac{\delta_n^2}{\kappa_2^2}, r_n+\frac{\delta_n^2}{\kappa_2^2}]^c=[0,\infty)$ and $\frac{|\rho'|}{|\mu_n-\rho|}\leq \frac{|\rho'|}{\rho}$, and thus, 
\begin{align*}
\| M_{\rho'} u_n\|_{L^2}
&  \leq \| \frac{\rho'}{\rho}  \|_{L^\infty} \|v_n \|_{L^2} +  \| \frac{\rho'}{\rho}  \|_{L^\infty}  \| \delta_n f_n \|_{L^2} \leq C \| v_n \|_{L^2} + C \delta_n \| f_n \|_{L^2}\,.
\end{align*}
Since $\|v_n \|_{L^\infty}\leq C$ by \eqref{proof.lem.Lambda.O1.1} and the interpolation inequality, collecting these above, we have 
\begin{align*}
\delta_n \| M_{\rho'} u_n \|_{L^2}\leq C \big (  \delta_n^3 |(\mu_n)_+|^\frac32 \| (-A)^\frac12 u_n \|_{L^2} \big )^\frac12 + C \kappa_2 + C \| f_n \|_{L^2} + C \delta_n \| v_n \|_{L^2}\,.
\end{align*}
This shows from \eqref{proof.lem.Lambda.O1.1} that $\displaystyle \limsup_{n\rightarrow \infty} \delta_n \| M_{\rho'} u_n \|_{L^2}\leq C \kappa_2 + C\limsup_{n\rightarrow \infty} \delta_n \| v_n \|_{L^2}$, and the claim is proved since $\kappa_2>0$ is arbitrary.

\

\noindent {\bf Step 3: Estimate of $\| Z u_n \|_{L^2}$\,.}

\noindent Let us recall from \eqref{proof.lem.Lambda.O1.3} that 
\begin{align}\label{proof.lem.Lambda.O1.7}
\langle W^* v_n, W^* \varphi \rangle_{L^2} = J[\varphi]:= \langle \frac{v_n + \delta_n f_n}{\mu_n-\rho}, \varphi\rangle_{L^2}\,, \qquad \varphi\in D(W^*)\,.
\end{align}
Here the right-hand side is well-defined in virtue of \eqref{proof.lem.Lambda.O1.4}.
Let $\kappa_3>0$ be fixed but arbitrary small number.
We decompose $J$ as $J_1+J_2$, where 
\begin{align*}
J_1[\varphi] = \int_{r_n-\kappa_3^2 \delta_n^2}^{r_n+\kappa_3^2\delta_n^2}\frac{v_n + \delta_n f_n}{\mu_n-\rho} \cdot \varphi\, \dd r\,, \qquad J_2 [\varphi] = J[\varphi]-J_1[\varphi]\,.
\end{align*}
Let us estimate $J_1$. In virtue of \eqref{proof.lem.Lambda.O1.4}, we have 
\begin{align*}
| J_1[\varphi]| & \leq  \int_{r_n-\kappa_3^2 \delta_n^2}^{r_n+\kappa_3^2\delta_n^2}\frac{\rho' |v_n -v_n (r_n)+ \delta_n (f_n-f_n (r_n))|}{|\mu_n-\rho|} \, \dd r \, \| M_{\frac{1}{\rho'}} \varphi\|_{L^\infty ([\frac{r_n}{2},\infty))} \\
& \leq C \kappa_3 \delta_n \| \partial_r v_n \|_{L^2} \| \varphi\|_{D(W^*)}  + C \kappa_3 \delta_n^2 \| \partial_r f_n \|_{L^2([r_n-\kappa_3^2\delta_n^2, r_n+\kappa_3\delta_n^2])} \| \varphi\|_{D(W^*)} \,.
\end{align*}
Here $\|\varphi \|_{D(W^*)} = \| W^* \varphi\|_{L^2} + \| \varphi \|_{L^2}$, and we have used $\inf_n r_n >0$ (by the condition $\mu_n\leq \frac12$) and $\| M_{\frac{1}{\sigma'}} \varphi \|_{L^\infty ([\frac{r_n}{2},\infty))} \leq C \| \varphi \|_{D(W^*)}$. As in \eqref{proof.lem.Lambda.O1.6}, it follows from \eqref{proof.lem.Lambda.O1.3} that 
\begin{align*}
\delta_n \| \partial_r f_n \|_{L^2([r_n-\kappa_3^2\delta_n^2, r_n+\kappa_3\delta_n^2])} \leq C \kappa_3^2 \delta_n^2 (1+r_n)^{-3} \| \partial_r u_n \|_{L^2} + C \| M_{\rho'} u_n \|_{L^2} + C \| \partial_r v_n \|_{L^2}\,,
\end{align*}
and thus, 
\begin{align}\label{proof.lem.Lambda.O1.8}
|J_1[\varphi]| & \leq C \kappa_3 \| \varphi \|_{D(W^*)} + C \Big ( \delta_n^3 |(\mu_n)_+|^\frac32 \| (-A)^\frac12 u_n \|_{L^2} + \kappa_3 \Big ) \| \varphi \|_{D(W^*)} \nonumber \\
& \leq  C \Big (  \delta_n^3 |(\mu_n)_+|^\frac32 \| (-A)^\frac12 u_n \|_{L^2} + \kappa_3 \Big ) \| \varphi \|_{D(W^*)} \,.
\end{align}
As for $J_2$, we decompose it as $J_{2,1} + J_{2,2}+J_{2,3}$, where 
\begin{align*}
J_{2,1}[\varphi] & = \int_{[r_n-\kappa_3^2 \delta_n^2, \, r_n + \kappa_3^2\delta_n^2]^c} \frac{1}{\mu_n-\rho}  \big ( \frac{v_n}{\rho'}-\frac{v_n (r_n)}{\rho'(r_n)}\big )  \varphi \, \rho'  \dd r \,, \\
J_{2,2}[\varphi] & = \frac{v_n (r_n)}{\rho'(r_n)} \int_{[r_n-\kappa_3^2 \delta_n^2, \, r_n + \kappa_3^2\delta_n^2]^c}  \frac{1}{\mu_n-\rho} \varphi\,  \rho' \dd r\,, \\
J_{2,3}[\varphi] & =  \int_{[r_n-\kappa_3^2 \delta_n^2, \, r_n + \kappa_3^2\delta_n^2]^c} \frac{\delta_n f_n}{\mu_n-\rho}  \varphi \, \rho'  \dd r \,.
\end{align*}
The term $J_{2,1}$ is estimated as, from $v_n=Z^*Zu_n$, $\frac{|\rho'|}{|\mu_n-\rho|}\leq \frac{Cr_n}{|r_n-r|(1+r)}$, and the Hardy inequality,
\begin{align*}
|J_{2,1}[\varphi]| \leq C \| (1+r) \partial_r (\frac{Z^* Z u_n}{\rho'}) \|_{L^2} \| \varphi \|_{L^2} \leq C \| Z u_n \|_{L^2} \| \varphi \|_{L^2} \quad ({\rm by~\eqref{est.Z^*}})\,.
\end{align*} 
Note that this estimate is valid also including the case $r_n=\infty$.
Next, when $\kappa_3>0$ is sufficiently small, we can check that when $r_n<\infty$,
\begin{align}\label{proof.lem.Lambda.O1.9}
\big | \Big ( \int_{r_n-1}^{r_n-\kappa_3^2 \delta_n^2} + \int_{r_n + \kappa_3^2\delta_n^2}^{r_n+1} \Big ) \frac{1}{\mu_n-\rho}   \rho' \dd r \big | \leq C\,,
\end{align}
by considering the change of the variable $\rho (r)=s$. Hence, when $r_n<\infty$ we have 
\begin{align*}
|J_{2,2}[\varphi]| & \leq \frac{|v_n (r_n)|}{|\rho'(r_n)|}  \Big ( \big | \int_0^{r_n-1} \cdots \, \dd r \big | + \big | \big ( \int_{r_n-1}^{r_n-\kappa_3^2 \delta_n^2}  + \int_{r_n + \kappa_3^2\delta_n^2}^{r_n+1} \big ) \cdots \, \dd r\big | + \big | \int_{r_n+1}^\infty \cdots \, \dd r\big | \Big ) \\
& \leq \frac{|v_n (r_n)|}{|\rho'(r_n)|}  \Big ( \| \frac{\rho'}{\mu_n-\rho}\|_{L^2([0,r_n-1])} \| \varphi \|_{L^2} + |\varphi (r_n) | \big | \big ( \int_{r_n-1}^{r_n-\kappa_3^2 \delta_n^2}  + \int_{r_n + \kappa_3^2\delta_n^2}^{r_n+1} \big ) \frac{\rho'}{\mu_n-\rho} \, \dd r\big | \\
& \qquad + \big | \big ( \int_{r_n-1}^{r_n-\kappa_3^2 \delta_n^2}  + \int_{r_n + \kappa_3^2\delta_n^2}^{r_n+1} \big ) \frac{\rho'}{\mu_n-\rho} (\varphi-\varphi (r_n) ) \, \dd r\big | + \| \frac{\rho'}{\mu_n-\rho}\|_{L^2([r_n+1,\infty))} \| \varphi \|_{L^2} \Big )\\
& \leq C  \frac{|v_n (r_n)|}{|\rho'(r_n)|} \|   \varphi \|_{H^1} \leq C  \frac{|v_n (r_n)|}{|\rho'(r_n)|} \|  \varphi \|_{D(W^*)}\,.
\end{align*}
When $r_n=\infty$ we have $J_{2,2}[\varphi]=0$ since $\frac{v_n(r_n)}{\rho'(r_n)}=\infty \cdot 0=0$.
Finally it is not difficult to show 
\begin{align*}
|J_{2,3}[\varphi]|\leq C \| \frac{\rho'}{\mu_n-\rho} \|_{L^2 ([r_n-\kappa_3^2 \delta_n^2, \, r_n + \kappa_3^2\delta_n^2]^c)} \| \delta_n f_n \|_{L^2} \| \varphi \|_{L^\infty} \leq \frac{C}{\kappa_3\delta_n} \| \delta_n f_n \|_{L^2} \| \varphi \|_{H^1} \leq \frac{C}{\kappa_3} \| f_n \|_{L^2} \| \varphi \|_{D(W^*)}
\end{align*} 
for $\varphi\in D(W^*)$.
Hence, we have 
\begin{align}\label{proof.lem.Lambda.O1.10}
|J_2[\varphi]|\leq C \| Z u_n \|_{L^2} \| \varphi \|_{L^2} + C \Big ( \frac{|v_n (r_n)|}{|\rho'(r_n)|} + \frac{C}{\kappa_3} \| f_n \|_{L^2} \Big ) \| \varphi \|_{D(W^*)}\,, \qquad \varphi \in D(W^*)\,.
\end{align}
Estimates \eqref{proof.lem.Lambda.O1.9} and \eqref{proof.lem.Lambda.O1.10} with \eqref{proof.lem.Lambda.O1.1} yield for $\varphi\in D(W^*)$,
\begin{align}\label{proof.lem.Lambda.O1.11}
\big |\langle Z u_n, W^* \varphi \rangle_{L^2} \big | = |J[\varphi]|\leq C \| \varphi \|_{L^2} +  C \Big (  \delta_n^3 |(\mu_n)_+|^\frac32 \| (-A)^\frac12 u_n \|_{L^2} + \kappa_3 +  \frac{|v_n (r_n)|}{|\rho'(r_n)|} + \frac{\|f_n\|_{L~2}}{\kappa_3}  \Big ) \| \varphi \|_{D(W^*)}\,,
\end{align}
that is, from \eqref{proof.lem.Lambda.O1.1} and Step 2, 
\begin{align}\label{proof.lem.Lambda.O1.12}
\big | \langle W^* v_\infty, W^*\varphi\rangle_{L^2} \big | =\big | \langle W^* v_\infty, W^* \varphi\rangle_{L^2} \big | & = \limsup_{n\rightarrow \infty} \big | \langle Z u_n, W^*  \varphi \rangle_{L^2} \big |  \leq  C \| \varphi \|_{L^2} \,, \qquad \varphi\in D(W^*) \,.
\end{align}
Moreover,  by taking $\varphi=Z^*Z u_n$ in \eqref{proof.lem.Lambda.O1.11}, we conclude that
\begin{align*}
\| Z u_n \|_{L^2} \leq  C\| Z^* Z u_n \|_{L^2} + C \Big (  \delta_n^3 |(\mu_n)_+|^\frac32 \| (-A)^\frac12 u_n \|_{L^2} + \kappa_3 +  \frac{|v_n (r_n)|}{|\rho'(r_n)|} + \frac{\|f_n\|_{L~2}}{\kappa_3}  \Big ) \,.
\end{align*}
which gives again from \eqref{proof.lem.Lambda.O1.1} and Step 1, 
\begin{align}\label{proof.lem.Lambda.O1.13}
\limsup_{n\rightarrow \infty} \| Z u_n\|_{L^2} \leq C \limsup_{n\rightarrow \infty} \| Z^* Z u_n \|_{L^2}\,,
\end{align}
by taking $\kappa_3\rightarrow 0$ after $n\rightarrow \infty$.

\

\noindent {\bf Step 4: Completion of the proof of \eqref{est.lem.Lambda.O1.1}.}

\noindent  Suppose that $\displaystyle \limsup_{n\rightarrow \infty} \| Z^* Z u_n \|_{L^2}=0$. Then Step 2 and \eqref{proof.lem.Lambda.O1.13} imply that $\displaystyle \limsup_{n\rightarrow \infty} \delta_n \| M_{\rho'} u_n \|_{L^2} = \limsup_{n\rightarrow \infty} \|  Z u_n \|_{L^2}=0$, which contradicts with the normalized condition \eqref{proof.lem.Lambda.O1.1}. Hence we may assume that $\inf_n \| Z^* Z u_n \|_{L^2}>0$ by taking a suitable subsequence if necessary. This implies that the limit $\displaystyle v_\infty = \lim_{n\rightarrow \infty} Z^* Z u_n\in L^2 (\R_+)$, $W^* v_\infty \in L^2 (\R_+)$, is  nontrivial. Moreover, \eqref{proof.lem.Lambda.O1.12} shows that $W^* v_\infty \in D(W)$. Then it is not difficult to show from $(\mu_n-M_{\rho}) W W^* v_n = v_n + \delta_n f_n$ that $u_\infty = W W^* v_\infty \in L^2 (\R_+)$ satisfies $(\mu_\infty - M_{\rho} ) u_\infty = v_\infty= Z^* Z u_\infty$ in $L^2 (\R_+)$, that is, $\mu_\infty$ must be an eigenvalue of $\hat{\Lambda}$ in $L^2(\R_+)$ and $u_\infty$ is the associated eigenfunction. To achieve the contradiction it remains to show that $u_\infty \in Y$. This is proved as follows: since $W^* (r^\frac32 g)=0$ by the definition of $W^*$,
\begin{align*}
\langle u_\infty,  r^\frac32 g \rangle_{L^2} = \langle W W^* v_\infty, r^\frac32 g\rangle _{L^2} =  \langle W^* v_\infty, W^* (r^\frac32 g) \rangle_{L^2} =0\,.
\end{align*}
The proof of \eqref{est.lem.Lambda.O1.1} is complete.

\

\noindent {\it Proof of \eqref{est.lem.Lambda.O1.2}.} Note that \eqref{est.lem.Lambda.O1.2} is equivalent with $\| u\|_{L^2}^2 \leq C \delta^{-4} \| (\mu -\hat{\Lambda}) u \|_{L^2}^2 + C h_1 (\frac{1}{\delta^2}, \mu)^2 \| (-A)^\frac12 u \|_{L^2}^2$ for any $\delta\in (0,1]$, and this is equivalent with 
\begin{align}\label{proof.lem.Lambda.O1.14}
\delta^2 \| u\|_{L^2}^2 \leq C \delta^{-2} \| (\mu -\hat{\Lambda}) u \|_{L^2}^2 + C \delta^2 h_1 (\frac{1}{\delta^2}, \mu)^2 \| (-A)^\frac12 u \|_{L^2}^2\,, \qquad \delta\in (0,1]\,.
\end{align}
Hence we shall prove \eqref{proof.lem.Lambda.O1.14} by contradiction argument. Suppose that \eqref{proof.lem.Lambda.O1.14} does not hold. 
Then there exist $\{\delta_n, \mu_n\}_{n\in \N}$, $\delta_n\in (0,1]$, $\mu_n \in (-\frac12,\frac12]$, and $\{u_n\}\subset D(A)$ such that 
\begin{align*}
& \lim_{n\rightarrow \infty} \delta_n=\delta_\infty\in [0,1]\,, \qquad \lim_{n\rightarrow \infty}\mu_n=\mu_\infty \in [-\frac12,\frac12]\,,
\end{align*}
and 
\begin{align}\label{proof.lem.Lambda.O1.15}
\begin{split}
\delta_n^2 \| u_n\|_{L^2}^2=1\,,  \quad  \lim_{n\rightarrow \infty}  \bigg ( \delta_n^{-2}  \|  (\mu_n - \hat{\Lambda}) u_n \|_{L^2}^2 + \delta_n^2 h_1 (\frac{1}{\delta_n^2},\mu_n)^2 \| (-A)^\frac12  u_n\|_{L^2}^2 \bigg ) =0\,.
\end{split}
\end{align}
We set 
\begin{align}\label{proof.lem.Lambda.O1.16}
f_n = \delta_n^{-1} (\mu_n   - \hat{\Lambda} ) u_n\,, \qquad v_n = Z^* Z u_n\,,
\end{align}
By the definition of $h_1^2$, we see 
\begin{align}\label{proof.lem.Lambda.O1.17}
\delta_n^6 |(\mu_n)_+|^3 \| (-A)^\frac12 u_n\|_{L^2}^2 \leq C \delta_n^2 h_1(\frac{1}{\delta_n^2},\mu_n)^2 \| (-A)^\frac12 u_n\|_{L^2}^2 \rightarrow 0 \quad (n\rightarrow \infty)\,.
\end{align}
Hence, \eqref{est.lem.Lambda.O1.1} implies that 
\begin{align}\label{proof.lem.Lambda.O1.18}
\lim_{n\rightarrow \infty} \Big ( \delta_n\| M_{\rho'} u_n \|_{L^2} + \| Z u_n \|_{L^2} + \| Z^* Z u_n \|_{L^2} \Big ) =0\,.
\end{align}
If $\frac{1}{\delta_n^2}\leq \frac{1}{10|\mu_n|}$, then we see $h_1(\frac{1}{\delta_n^2},\mu_n)=\delta_n$ and 
\begin{align*}
\delta_n \| u_n\|_{L^2 ([\frac{1}{2\delta_n},\infty))} \leq C \delta_n^2 \| r u_n\|_{L^2} \leq C \delta_n h_1 (\frac{1}{\delta_n^2},\mu_n) \| (-A)^\frac12 u_n \|_{L^2}\,,
\end{align*}
while from $u_n=\frac{v_n+\delta_n f_n}{\mu_n-\rho}$ and $\frac{1}{2\delta_n}\leq \frac{r_n}{2}$,  
\begin{align*}
\delta_n \| u_n \|_{L^2 ([0,\frac{1}{2\delta_n}])} & \leq \delta_n \| \frac{v_n}{\mu_n-\rho} \|_{L^2 ([0,\frac{1}{2\delta_n}])} +\delta_n \| \frac{\delta_n f_n}{\mu_n-\rho} \|_{L^2 ([0,\frac{1}{2\delta_n}])} \\
& \leq C \delta_n \| (1+r)^2 v_n\|_{L^2} + \delta_n \|\frac{1}{\mu_n-\rho}\|_{L^\infty ([0,\frac{1}{2\delta_n}])}  \| \delta_n f_n \|_{L^2} \\
& \leq C \delta_n \| (1+r)^2 v_n\|_{L^2} + C \delta_n\cdot \delta_n^{-2} \| \delta_n f_n \|_{L^2}  \leq C \delta_n \| Z u_n \|_{L^2} + C \| f_n\|_{L^2}\,.
\end{align*}
Thus, we conclude when $\frac{1}{\delta_n^2}\leq \frac{1}{10|\mu_n|}$,
\begin{align}\label{proof.lem.Lambda.O1.19}
\delta_n \| u_n\|_{L^2} \leq C \delta_n h_1 (\frac{1}{\delta_n^2},\mu_n) \| (-A)^\frac12 u_n \|_{L^2} + C \delta_n \| Z u_n \|_{L^2} + C \| f_n\|_{L^2}\,.
\end{align}
If  $\mu_n>0$ and $\frac{1}{10\mu_n} \leq \frac{1}{\delta_n^2}\leq \frac{1}{\mu_n^2}$, then $r_n<\infty$ and we see $h_1(\frac{1}{\delta_n^2},\mu_n)=\mu_n^\frac12$ and 
\begin{align*}
\delta_n \| u_n\|_{L^2 ([\frac{r_n}{2},\infty))} \leq \frac{C\delta_n}{r_n} \| r u_n\|_{L^2} \leq C\delta_n \mu_n^\frac12\| (-A)^\frac12 u_n \|_{L^2} \leq C \delta_n h_1 (\frac{1}{\delta_n^2},\mu_n) \| (-A)^\frac12 u_n \|_{L^2}\,,
\end{align*}
while from $u_n=\frac{v_n+\delta_n f_n}{\mu_n-\rho}$,  
\begin{align*}
\delta_n \| u_n \|_{L^2 ([0,\frac{r_n}{2}])} & \leq \delta_n \| \frac{v_n}{\mu_n-\rho} \|_{L^2 ([0,\frac{r_n}{2}])} +\delta_n \| \frac{\delta_n f_n}{\mu_n-\rho} \|_{L^2 ([0,\frac{r_n}{2}])} \\
& \leq C \delta_n \| (1+r)^2 v_n\|_{L^2} + \delta_n \|\frac{1}{\mu_n-\rho}\|_{L^\infty ([0,\frac{r_n}{2}])}  \| \delta_n f_n \|_{L^2} \\
& \leq C \delta_n \| (1+r)^2 v_n\|_{L^2} + C \delta_n\cdot r_n^{2} \| \delta_n f_n \|_{L^2}  \leq C \delta_n \| Z u_n \|_{L^2} + \frac{C \delta_n^2}{\mu_n} \| f_n\|_{L^2}\,.
\end{align*}
Thus, we conclude when $\mu_n>0$ and $\frac{1}{10\mu_n}\leq \frac{1}{\delta_n^2}\leq \frac{1}{\mu_n^2}$,
\begin{align}\label{proof.lem.Lambda.O1.20}
\delta_n \| u_n\|_{L^2} \leq C \delta_n h_1 (\frac{1}{\delta_n^2},\mu_n) \| (-A)^\frac12 u_n \|_{L^2} + C \delta_n \| Z u_n \|_{L^2} + C \| f_n\|_{L^2}\,.
\end{align}
Next we consider the case $\mu_n>0$ and $\frac{1}{\delta_n^2}\geq \frac{1}{\mu_n^2}$.
In this case $r_n<\infty$ and $h_1(\frac{1}{\delta_n^2},\mu_n) = \frac{\delta_n^2}{\mu_n^\frac32}$.
We observe that 
\begin{align*}
\delta_n^2 \int_{r_n-\frac{\delta_n^2}{\mu_n^\frac32}}^{r_n+\frac{\delta_n^2}{\mu_n^\frac32}} u_n^2 \, \dd r\leq \frac{C\delta_n^4}{\mu_n^\frac32} \|u_n\|_{L^\infty}^2 \leq  \frac{C\delta_n^4}{\mu_n^\frac32} \| \partial_r u_n\|_{L^2} \| u_n\|_{L^2}
& \leq \frac{C\delta_n^3}{\mu_n^\frac32} \| (-A)^\frac12 u_n \|_{L^2} \\
& \leq C\delta_n h_1 (\frac{1}{\delta_n^2},\mu_n) \| (-A)^\frac12 u_n \|_{L^2}\,.
\end{align*}
On the other hand, we have from $u_n=\frac{v_n+\delta_n f_n}{\mu_n-\rho}$, 
\begin{align*}
\delta_n \|u_n\|_{L^2([r_n-\frac{\delta_n^2}{\mu_n^{3/2}},\, r_n+\frac{\delta_n^2}{\mu_n^{3/2}}]^c)} 
& \leq 
\delta_n \|\frac{v_n}{\mu_n-\rho}\|_{L^2([r_n-\frac{\delta_n^2}{\mu_n^{3/2}}, \, r_n+\frac{\delta_n^2}{\mu_n^{3/2}}]^c)} + \delta_n \| \frac{\delta_n f_n}{\mu_n-\rho} \|_{L^2([r_n-\frac{\delta_n^2}{\mu_n^{3/2}}, \, r_n+\frac{\delta_n^2}{\mu_n^{3/2}}]^c)} \\
& \leq \delta_n \|\frac{\rho'}{\mu_n-\rho}\|_{L^2([r_n-\frac{\delta_n^2}{\mu_n^{3/2}}, \, r_n+\frac{\delta_n^2}{\mu_n^{3/2}}]^c)} \| M_{\frac{1}{\rho'}} v_n \|_{L^\infty}  \\
& \quad +  \delta_n \|\frac{1}{\mu_n-\rho}\|_{L^\infty ([r_n-\frac{\delta_n^2}{\mu_n^{3/2}}, \, r_n+\frac{\delta_n^2}{\mu_n^{3/2}}]^c)}\|\delta_n f_n\|_{L^2} \\
& \leq C \delta_n \cdot \frac{\mu_n^\frac43}{\delta_n} \| Z u_n\|_{L^2} + C \delta_n \cdot (1+r_n)^3 \frac{\mu_ n^\frac32}{\delta_n^2} \| \delta_n f_n \|_{L^2} \\
& \leq C \| Z u_n \|_{L^2} + C \| f_n \|_{L^2}\,.
\end{align*}
Here we have used the fact that $r_n^3 \mu_n^\frac32\leq C$ since $r_n\approx \mu_n^{-\frac12}$ for $0<\mu_ n\ll 1$. Thus, we have when $\frac{1}{\delta_n^2}\geq \frac{1}{\mu_n^2}$, 
\begin{align}\label{proof.lem.Lambda.O1.21}
\delta_n \| u_n \|_{L^2} \leq C\delta_n h_1 (\frac{1}{\delta_n^2},\mu_n) \| (-A)^\frac12 u_n \|_{L^2} + C \| Z u_n \|_{L^2} + C \| f_n\|_{L^2}\,.
\end{align}
Collecting \eqref{proof.lem.Lambda.O1.19}, \eqref{proof.lem.Lambda.O1.20}, and \eqref{proof.lem.Lambda.O1.21}, we obtain by applying \eqref{proof.lem.Lambda.O1.1} and \eqref{proof.lem.Lambda.O1.18},
\begin{align*}
\lim_{n\rightarrow \infty} \delta_n \| u_n \|_{L^2} =0\,.
\end{align*}
This contradicts with the normalized condition \eqref{proof.lem.Lambda.O1.15}. The proof of \eqref{est.lem.Lambda.O1.2} is complete.

\

The analysis of the case $|\mu|\geq \frac12$ is similar to the case of the Kolmogorov flow.
First we observe that 
\begin{lem}\label{lem.Lambda.O2} There exist $\kappa\in (0,1)$ and $C>0$ such that the following statements hold for all $\delta\in (0,1]$.
If  $\mu\in \R$ satisfies $1-\kappa \delta^2 \leq \mu\leq 1+\kappa \delta^2$ then 
\begin{align}\label{est.lem.Lambda.O2.1}
\delta^2 \| u \|_{L^2}^2  + \| M_{\rho'} u \|_{L^2}^2 \leq C \big ( \delta^{-2} \| (\mu-\hat{\Lambda} ) u \|_{L^2}^2 + \delta^4 \| (-A)^\frac12 u \|_{L^2}^2 \big )\,, \qquad u\in D (A)\,,
\end{align}
while if $\mu>1$ then 
\begin{align}\label{est.lem.Lambda.O2.2}
\begin{split}
& (\mu-1)^2 \|  u \|_{L^2}^2  + (\mu-1) \| M_{\rho'}  u\|_{L^2}^2  \leq C  \| (\mu-\hat{\Lambda} ) u \|_{L^2}^2 \,,  \qquad u\in  D (A)\,.
\end{split}
\end{align}
Moreover, if $\mu<0$ then 
\begin{align}\label{est.lem.Lambda.O2.3}
\mu^2 \| u\|_{L^2}^2 +  \| M_{\rho'} u\|_{L^2}^2  \leq C \| (\mu-\hat{\Lambda})u\|_{L^2}^2\,, \qquad u\in D(A)\cap Y\,.
\end{align} 
\end{lem}

\noindent {\it Proof.} The bounds $\|u\|_{L^2} \leq \frac{1}{\mu-1} \| (\mu-\hat{\Lambda}) u \|_{L^2}$ when $\mu>1$ and $\|u\|_{L^2} \leq \frac{1}{|\mu|} \| (\mu-\hat{\Lambda}) u \|_{L^2}$ when $\mu<0$, stated in \eqref{est.lem.Lambda.O2.1} and \eqref{est.lem.Lambda.O2.2} respectively,  is a direct consequence of the fact that $\hat{\Lambda}$ is a self-adjoint operator in $L^2 (\R_+)$ with the spectrum $\sigma (\hat{\Lambda})=[0,1]$. Then \eqref{est.lem.Lambda.O2.3} is proved, for the desired  estimate of $\|M_{\rho'} u\|_{L^2}$ in the case $\mu\in [-\frac12, 0]$ is already shown in Lemma \ref{lem.Lambda.O1}, while the estimate for the case $\mu\leq -\frac12$ follows from the estimate of $\| u\|_{L^2}$.
To show the other estimates we set $f=(\mu-\hat{\Lambda})u$. 
Then we have 
\begin{align}\label{proof.lem.Lambda.O2.1}
\int_0^\infty (1-\rho) |u|^2 \, \dd r +  \| Z u\|_{L^2}^2 = (1-\mu) \| u\|_{L^2}^2 +  \langle f, u \rangle_{L^2}\,.
\end{align}
Since there exists $C_0>0$ such that $(\rho')^2\leq C_0 (1-\rho)$ in $[0,\infty)$, \eqref{proof.lem.Lambda.O2.1} implies 
\begin{align}\label{proof.lem.Lambda.O2.2}
\| M_{\rho'} u\|_{L^2}^2 \leq C_0 \int_0^\infty (1-\rho) |u|^2 \, \dd r \leq C (1-\mu) \| u\|_{L^2}^2 + C \langle f, u\rangle_{L^2}\,.
\end{align}
This proves \eqref{est.lem.Lambda.O2.2} for the case $\mu>1$. Moreover, if $\mu\geq 1-\kappa\delta^2$ then \eqref{proof.lem.Lambda.O2.2} gives $\|M_{\rho'} u\|_{L^2}^2 \leq C \delta^2 \|u\|_{L^2}^2 + C \delta^{-2} \|f\|_{L^2}^2$, and thus, 
it suffices to consider the estimate of $\|u\|_{L^2}^2$ for small enough$\delta$ to complete the proof of \eqref{est.lem.Lambda.O2.1}. We see
\begin{align*}
\int_0^\infty (1-\rho) |u|^2 \, \dd r \geq \int_\delta^\infty (1-\rho) |u|^2 \, \dd r \geq \frac{\delta^2}{C} \int_\delta^\infty |u|^2\, \dd r & \geq \frac{\delta^2}{C} \big ( \| u\|_{L^2}^2 - \int_0^\delta |u|^2 \, \dd r \big ) \\
& \geq \frac{\delta^2}{C} \big (\|u\|_{L^2}^2 - C \delta \| u\|_{L^\infty}^2\big )\,.
\end{align*}
This gives from \eqref{proof.lem.Lambda.O2.1}, by taking $\kappa>0$ small enough,
\begin{align*}
\delta^2 \| u\|_{L^2}^2 \leq C \kappa\delta^2 \|u\|_{L^2}^2 + C \|f\|_{L^2} \|u\|_{L^2} + C \delta^3 \|u\|_{L^\infty}^2 \leq C\delta^{-2} \|f\|_{L^2}^2 + C \delta^4 \| (-A)^\frac12 u\|_{L^2}^2\,.
\end{align*}
Estimate \eqref{est.lem.Lambda.O2.1} is proved. The proof is complete.

\

The proof of the following lemma is very parallel to the proof of Lemma \ref{lem.Lambda.Kol.rate.2}.
\begin{lem}\label{lem.Lambda.O3}
Let $\kappa\in (0,1)$ be the number in Lemma \ref{lem.Lambda.O2}.
There exists $C>0$  such that if $\delta\in (0,1]$,  and $\mu\in \R$ with $\frac12\leq \mu < 1-\kappa \delta^2$, then 
\begin{align}\label{est.lem.Lambda.O3.1}
\delta^2 \| u \|_{L^2}^2 + \| Z u \|_{L^2}^2 + \frac{1}{1-\mu} \| Z^* Z u\|_{L^2}^2 \leq C \bigg ( \delta^{-2} \| (\mu-\hat{\Lambda}) u \|_{L^2}^2 + \frac{\delta^6}{1-\mu} \| (-A)^\frac12 u\|_{L^2}^2  \bigg ) \,, \quad u\in D (A) \,,
\end{align}
and 
\begin{align}\label{est.lem.Lambda.O3.2}
\| M_{\rho'} u \|_{L^2}^2 \leq C \bigg ( \delta^{-2} \| (\mu-\hat{\Lambda}) u \|_{L^2}^2 + \delta^2 (1-\mu ) \, \| (-A)^\frac12 u\|_{L^2}^2 \bigg ) \,, \quad u\in D (A) \,.
\end{align}
\end{lem}

\noindent {\it Proof of \eqref{est.lem.Lambda.O3.1}.} We only state the outline of the proof, for the argument is parallel to the proof of Lemma \ref{lem.Lambda.Kol.rate.2}.
Suppose that \eqref{est.lem.Lambda.O3.1} does not hold. 
Then there exist $\{\delta_n, \mu_n\}_{n\in \N}$, $\delta_n \in (0,1]$, $\mu_n\in  [\frac12, 1-\kappa \delta_n^2)$, and $\{u_n\}\subset D (A)$ such that 
\begin{align*}
\lim_{n\rightarrow \infty} \delta_n=\delta_\infty\in [0,1]\,,
\qquad  \lim_{n\rightarrow \infty}\mu_n=\mu_\infty \in [\frac12, 1-\kappa \delta_\infty^2]\,,
\end{align*}
and 
\begin{align}\label{proof.prop.Lambda.O3.1}
\begin{split}
& \delta_n^2 \| u_n \|_{L^2}^2 + \| Z u_n \|_{L^2}^2 + \frac{1}{1-\mu_n} \| Z^* Z u_n \|_{L^2}^2 =1\,, \\
& \lim_{n\rightarrow \infty}  \bigg ( \delta_n^{-2}  \|  (\mu_n - \hat{\Lambda}) u_n \|_{L^2}^2 + \frac{\delta_n^6}{1-\mu_n} \| (-A)^\frac12   u_n \|_{L^2}^2 \bigg ) =0\,.
\end{split}
\end{align}
Set
\begin{align}\label{proof.prop.Lambda.O3.2}
f_n = \delta_n^{-1} (\mu_n  - \hat{\Lambda}) u_n\,, \qquad v_n = Z^* Z u_n\,, \qquad r_n = \rho^{-1} (\mu_n)\in (0,10)\,,
\end{align}
and then $v_n$ satisfies 
\begin{align}\label{proof.prop.Lambda.O3.3}
(\mu_n - M_{\rho}) W^* W v_n =  v_n + \delta_n f_n\,.
\end{align}
Note that 
\begin{align}\label{proof.prop.Lambda.O3.4}
v_n (r_n) + \delta_n f_n (r_n) = 0
\end{align}
holds.
We may assume that $\delta_\infty=0$ (otherwise the contradiction is easily achieved).
We may also assume from \eqref{proof.prop.Lambda.O3.1} that $v_n$ converges to a limit $v_\infty$ strongly in $L^2(\R_+)$ and $W^* v_n$ converges to $W^*v_\infty$ weakly in $L^2 (\R_+)$. 
The direct computation implies that 
\begin{align}
\frac{1}{C|1-\mu_n|^\frac12} \leq |\rho'(r_n)| \leq \frac{C}{|1-\mu_n|^\frac12}\,,
\end{align}
for $\rho'(r)<0$ and $\rho'(r) \approx -\frac{r}{4}$ near $r=0$. 
As in the proof of Lemma \ref{lem.Lambda.Kol.rate.2}, we can show the following claims.

\

\noindent {\bf Step 1: $\displaystyle \lim_{n\rightarrow \infty} \delta_n \|u_n \|_{L^2}=0$.} \qquad \noindent {\bf Step 2: $\displaystyle \lim_{n\rightarrow \infty} \frac{v_n (r_n)}{|\rho' (r_n)|^\frac12} = \lim_{n\rightarrow \infty} \frac{\delta_n f_n (r_n)}{|\rho' (r_n)|^\frac12} =0$.}

\noindent {\bf Step 3: $\displaystyle \lim_{n\rightarrow \infty} \frac{\| v_n \|_{L^2 ([r_n,\infty))}}{|\rho'(r_n)|} =0$.} \qquad \noindent {\bf Step 4: Estimate of $\| W^* v_n\|_{L^2}$.}

\noindent In Step 4 we verify the estimates 
\begin{align}\label{proof.prop.Lambda.O3.5}
\limsup_{n\rightarrow \infty} \| W^* v_n\|_{L^2} \leq \frac{C}{\kappa_4} \limsup_{n\rightarrow \infty} \frac{\| v_n\|_{L^2([0,r_n])}}{|\rho'(r_n)|} + C \kappa_4^\frac12
\end{align}
for any sufficiently small $\kappa_4>0$, and we can also show when $\frac12\leq \mu_\infty <1$,
\begin{align}\label{proof.prop.Lambda.O3.6}
|\langle W^* v_\infty, W^* \varphi\rangle_{L^2}\big | \leq C\| \varphi \|_{L^2}\,, \qquad \varphi\in D(W^*)\,,
\end{align}
The details of the proof of the above steps are omitted here.
In virtue of \eqref{proof.prop.Lambda.O3.5} and Step 1 we may assume that $\inf_n \frac{\| v_n\|_{L^2([0,r_n])}}{|\rho'(r_n)|}>0$ (by taking a suitable subsequence if necessary), otherwise we achieve the contradiction to the normalized condition \eqref{proof.prop.Lambda.O3.1}.  
If $\frac12 \leq \mu_\infty <1$ then the limit $v_\infty$ is nontrivial, and \eqref{proof.prop.Lambda.O3.6} implies that $W^* v_\infty\in D(W)$, that is, $u_\infty= W W^* v_\infty\in L^2(\R_+)$ is an eigenfunction to the eigenvalue $\mu_\infty$ of $\hat{\Lambda}$ in $L^2 (\R_+)$, which contradicts with the absence of the eigenvalues in $[\frac12,1]$.
It remains to consider the case $\mu_\infty=1$, for which we need a rescaling process.
Set 
\begin{align*}
w_n (s) = \frac{1}{|1-\mu_n|^\frac14} v_n ( |1-\mu_n|^\frac12 s)\,, \qquad s\in [0,10]\,, \qquad s_n = \frac{r_n}{|1-\mu_n|^\frac12}\,.
\end{align*}
Then we have 
\begin{align}\label{proof.prop.Lambda.O3.5}
\begin{split}
\| \partial_s w_n \|_{L^2 ([0,10])}^2 + \| w_n \|_{L^2([0,10])}^2 & \leq \| \partial_r v_n \|_{L^2}^2 + \frac{\| v_n \|_{L^2}^2}{1-\mu_n} \leq 1\,,\\
\inf_n \| w_n \|_{L^2 ([0,c_n])}^2 & = \inf_n \frac{\| v_n \|_{L^2([0,r_n])}^2}{1-\mu_n}>0\,.
\end{split}
\end{align}
Hence we may assume that $w_n$ converges to a limit $w_\infty$ weakly in $H^1 ([0,10])$ and strongly in $L^2 ([0,10])$. The lower bound in \eqref{proof.prop.Lambda.O3.5} implies that $w_\infty$ is nontrivial.
Since $\rho (r) = 1 + \frac12 (-\frac{r^2}{4}) + \frac16 (-\frac{r^2}{4})^2 + \cdots$, we have from $\mu_n = \rho (r_n)$ that $1-\mu_n= \frac{r_n^2}{8} \big ( 1-\frac{r_n^2}{12} + O(r_n^4)\big )$ for $0<1-\mu_n\ll 1$, which gives 
\begin{align*}
r_n = 2\sqrt{2} (1-\mu_n)^\frac12 \big ( 1+ O (r_n^2) \big )\,,
\end{align*}
and hence, 
\begin{align}\label{proof.prop.Lambda.O3.6}
s_n = 2\sqrt{2} \big ( 1 + O (r_n^2) \big ) ~~ \rightarrow  ~~2\sqrt{2} \quad (n\rightarrow \infty)\,.
\end{align}
By Step 2 and Step 3 we also have 
\begin{align}\label{proof.prop.Lambda.O3.7}
w_\infty (s) =0 \qquad s\in [2\sqrt{2},10]\,.
\end{align}
To obtain the equation for $w_\infty$ we observe that 
\begin{align*}
WW^* v_n = - \frac{1}{gr^\frac32} \frac{\dd}{\dd r} \Big ( r^3 \frac{\dd }{\dd r} (\frac{v_n}{gr^\frac32}) \Big ) = - \frac{1}{g^2} \Big ( v_n''  +\frac{r}{2} v_n' + (\frac{r^2}{16} + \frac14 - \frac{3}{4r^2} ) v_n \Big )\,. 
\end{align*}
Hence, $(\mu_n-\rho (r)) W W^* v_n = v_n + \delta_n f_n$ is written as 
\begin{align*}
\big (1-\mu_n - \frac{r^2}{8} + O(r^4)\big ) \Big ( v_n'' + \frac{r}{2} v_n'  +  (\frac{r^2}{16} + \frac14 - \frac{3}{4r^2} ) v_n \Big  ) = g^2 \big ( v_n + \delta_n f_n \big )
\end{align*}
and then, $v_n (r) = |1-\mu_n|^\frac14 w_n (\frac{r}{|1-\mu_n|^\frac12})$ shows 
\begin{align}\label{proof.prop.Lambda.O3.8}
\big ( 1-\frac{s^2}{8} + q_n (s) \big ) \Big ( w_n'' (s) + \frac{|1-\mu_ n| s}{2} w_n'  + ( -\frac{3}{4s^2} + p_n (s) ) w_n \Big ) = \big (1+b_n (s)  \big ) (w_n + \tilde f_n )\,,
\end{align}
where $1-\frac{s_n^2}{8} + q_n (s_n)=0$,  
$\displaystyle \lim_{n\rightarrow \infty} \big ( \| q_n \|_{C^1 ([0,10])} + \| p_n \|_{C^1([0,10])} + \| b_n \|_{C^1([0,10])} \big )=0$, and 
\begin{align}\label{proof.prop.Lambda.O3.9}
\tilde f_n = \frac{\delta_n}{|1-\mu_n|^\frac14} f_n (|1-\mu_n|^\frac12 s)\,,
\end{align}
which satisfies 
\begin{align}\label{proof.prop.Lambda.O3.10}
\| \tilde f_n \|_{L^2([0,10])} \leq  \frac{\delta_n}{|1-\mu_n|^\frac12} \| f_n \|_{L^2} \rightarrow 0\qquad (n\rightarrow \infty)\,.
\end{align}
Let $\kappa_5>0$ be fixed but arbitrary small number.
Then, by arguing as Step 6 in the proof of Lemma \ref{lem.Lambda.Kol.rate.2}, we can show 
\begin{align}\label{proof.prop.Lambda.O3.11}
\big | \langle \partial_s w_\infty, \partial_s \varphi \rangle_{L^2 ([\kappa_5,10])} \big | \leq C_{\kappa_5} \| \varphi \|_{L^2([\kappa_5,10])}\,, \qquad \varphi\in C_0^\infty ((\kappa_5,10))\,.
\end{align}
Here $C_{\kappa_5}$ depends only on $\kappa_5$. This implies $w_\infty \in H^2 ((\kappa_5,10))$. Since $\kappa_5>0$ is arbitrary, we conclude that $w_\infty\in H^2_{loc} ((0,10))\cap H^1 ([0,10])$.
We can also check from \eqref{proof.prop.Lambda.O3.8} that  $w_\infty$ satisfies 
\begin{align}\label{proof.prop.Lambda.O3.12}
(1-\frac{s^2}{8}) \big ( w_\infty'' -\frac{3}{4s^2} w_\infty \big ) = w_\infty\,, \qquad s\in (0,10)\setminus \{2\sqrt{2}\}\,.
\end{align}
The regularity $w_\infty\in H^2_{loc} ((0,10))\cap H^1 ([0,10])$ with \eqref{proof.prop.Lambda.O3.7} yields $w_\infty (2\sqrt{2})=\partial_s w_\infty(2\sqrt{2})=0$, and it is easy to show that any solution to \eqref{proof.prop.Lambda.O3.12} in $H^2_{loc} ((0,10))\cap H^1 ([0,10])$ satisfying the condition $w_\infty (2\sqrt{2})=\partial_s w_\infty(2\sqrt{2})=0$ is trivial.
This is a contradiction. The proof of \eqref{est.lem.Lambda.O3.1} is complete.

\

\noindent {\it Proof of \eqref{est.lem.Lambda.O3.2}.} The argument is parallel to the proof of \eqref{est.lem.Lambda.Kol.rate.2'} in Lemma \ref{lem.Lambda.Kol.rate.2}.
Suppose that \eqref{est.lem.Lambda.O3.2} does not hold. Then there exist $\{\delta_n, \mu_n\}_{n\in \N}$, $\delta_n\in (0,1]$, $\mu_n \in  [\frac12, 1-\kappa \delta_n^2)$, and $\{u_n\}\subset D (A)$ such that 
\begin{align*}
\lim_{n\rightarrow \infty} \delta_n=\delta_\infty\in [0,1]\,,
\qquad  \lim_{n\rightarrow \infty}\mu_n=\mu_\infty \in [\frac12, 1-\kappa \delta_\infty^2]\,,
\end{align*}
and 
\begin{align}\label{proof.prop.Lambda.O3.1'}
\begin{split}
\| M_{\rho'} u_n \|_{L^2}^2 =1\,, \qquad  \lim_{n\rightarrow \infty}  \bigg ( \delta_n^{-2}  \|  (\mu_n - \hat{\Lambda}) u_n \|_{L^2}^2 + \delta_n^2 (1-\mu_n) \| (-A)^\frac12   u_n \|_{L^2}^2 \bigg ) =0\,.
\end{split}
\end{align}
Set
\begin{align}\label{proof.prop.Lambda.O3.2'}
f_n = \delta_n^{-1} (\mu_n  - \hat{\Lambda}) u_n\,, \qquad v_n = Z^* Z u_n\,, \qquad r_n = \rho^{-1} (\mu_n)\in (0,10)\,,
\end{align}
and then $v_n$ satisfies 
\begin{align}\label{proof.prop.Lambda.O3.3'}
(\mu_n - M_{\rho}) W^* W v_n =  v_n + \delta_n f_n\,.
\end{align}
Note that 
\begin{align}\label{proof.prop.Lambda.O3.4'}
v_n (r_n) + \delta_n f_n (r_n) = 0
\end{align}
holds. Since $\frac{\delta_n^6}{1-\mu_n} \| (-A)^\frac12 u_n \|_{L^2}^2 \leq C\delta_n^2 (1-\mu_n) \| (-A)^\frac12 u_n \|_{L^2}^2 \rightarrow 0 ~~(n\rightarrow \infty)$, we have from \eqref{est.lem.Lambda.O3.1} that 
\begin{align}\label{proof.prop.Lambda.O3.5'}
\lim_{n\rightarrow \infty} \big ( \delta_n^2 \| u_n \|_{L^2}^2 + \| Z u_n \|_{L^2}^2 + \| Z^* Z u_n \|_{L^2}^2 \big ) =0\,.
\end{align} 
Then, as in \eqref{proof.lem.Lambda.Kol.rate.2-2.-7}, we can show 
\begin{align*}
\limsup_{n\rightarrow \infty} \| M_{\rho'} u_n\|_{L^2} \leq C \limsup_{n\rightarrow \infty} \frac{|v_n (r_n)|}{\delta_n^\frac12}\,.
\end{align*}
The fact $\displaystyle \limsup_{n\rightarrow \infty} \frac{|v_n (r_n)|}{\delta_n^\frac12}=0$ is proved in the same manner as in the case of Lemma \ref{lem.Lambda.Kol.rate.2}, by investigating $\delta_n f_n (r_n)$ (here, recall \eqref{proof.prop.Lambda.O3.4'}). We omit the details. The proof is complete.

\

Let $\kappa\in (0,1)$ be the number in Lemma \ref{lem.Lambda.O2}.
Taking \eqref{def.h_1.O.1} into account, we refine $h_1(m,\mu)^2$ for $m\geq 100$ and $\mu\in \R$ as follows.
\begin{align}\label{def.h_1.O.2}
h_1^2 (m, \mu)  = 
\begin{cases}
& \displaystyle 0\,, \qquad \qquad \qquad \mu>1 + \frac{\kappa}{m}\,,\\
& \displaystyle \frac{1}{m}\,, \qquad \qquad 1-\frac{\kappa}{m} < \mu \leq 1+ \frac{\kappa}{m}\,,\\
& \displaystyle \frac{1}{m^2(1-\mu)}\,, \qquad \quad \frac12 < \mu \leq 1-\frac{\kappa}{m}\,,\\ 
& \displaystyle \frac{1}{m^2\mu^3} \qquad  \qquad \quad \frac{1}{m^\frac12} < \mu\leq \frac12\,,\\
& \displaystyle \mu \qquad \qquad  \qquad \frac{1}{10m} < \mu\leq \frac{1}{m^\frac12}\,,\\
& \displaystyle \frac{1}{m} \qquad \qquad   ~~-\frac{1}{10m}\leq \mu \leq \frac{1}{10m}\,, \\
& \displaystyle 0\,, \qquad \qquad \qquad \mu<-\frac{1}{10m}\,.
\end{cases}
\end{align}
We also define $h_2(m,\mu)^2$ for $m\geq 100$ and $\mu\in \R$ as  
\begin{align}\label{def.h_2.O.1}
h_2^2 (m, \mu)  = 
\begin{cases}
& \displaystyle 0\,, \qquad \qquad \qquad \mu>1 + \frac{\kappa}{m^2}\,,\\
& \displaystyle \frac{1}{m^4}\,, \qquad \qquad 1-\frac{\kappa}{m^2} < \mu \leq 1+ \frac{\kappa}{m^2}\,,\\
& \displaystyle \frac{1-\mu}{m^2}\,, \qquad \quad \frac12 < \mu \leq 1-\frac{\kappa}{m^2}\,,\\ 
& \displaystyle \frac{\mu^3}{m^2} \qquad  \qquad \quad 0 < \mu\leq \frac12\,,\\
& \displaystyle 0\,, \qquad \qquad \qquad \mu\leq 0\,.
\end{cases}
\end{align}

Note that $\displaystyle \lim_{m\rightarrow \infty} \sup_{\mu\in \R} h_j (m,\mu)=0$ holds.
Lemmas \ref{lem.Lambda.O1}, \ref{lem.Lambda.O2}, and \ref{lem.Lambda.O3} yield
\begin{prop} Let $m\geq 100$ and $\mu\in \R$. Let $h_j(m,\mu)$, $j=1,2$, be the nonnegative function defined by \eqref{def.h_1.O.2} and \eqref{def.h_2.O.1}. Then there exists a positive constant $C$ independent of $m$ and $\mu$ such that for any $u\in D(A)\cap Y$, 
\begin{align}
\| u\|_{L^2}^2 & \leq C m^2 \| (\mu-\hat{\Lambda}) u\|_{L^2}^2 + C h_1^2 (m,\mu) \| (-A)^\frac12 u\|_{L^2}^2\,\\
\| M_{\rho'} u\|_{L^2}^2 & \leq C m^2 \| (\mu-\hat{\Lambda}) u\|_{L^2}^2 + C h_2^2 (m,\mu) \| (-A)^\frac12 u \|_{L^2}^2\,.
\end{align} 
\end{prop}

\

To obtain the resolvent estimate by applying Theorem \ref{thm.abstract.2} we need to evaluate the function 
\begin{align*}
F(\alpha,\frac{\lambda}{\alpha}) = \inf_{m_1,m_2\geq 100} \Big ( \frac{m_1}{|\alpha|} + \frac{m_1^2m_2^2}{\alpha^2} + \frac{m_1^2 h_2 (m_2,\frac{\lambda}{\alpha})}{|\alpha|} + h_1^2 (m_1, \frac{\lambda}{\alpha}) \Big )\,.
\end{align*} 

Set $\mu=\frac{\lambda}{\alpha}$.

\noindent {\bf Case 1: $\mu>1+ \frac{\kappa}{|\alpha|^\frac12}$.} Take $m_1=m_2^2=\frac{2\kappa}{\mu-1}$, for which $1+\frac{\kappa}{m_1}=1+\frac{\kappa}{m_2^2}=1+\frac{\mu-1}{2}<\mu$ holds. Then $h_1(m_1,\mu)^2=h_2(m_2,\mu)^2=0$ and we have $F(\alpha,\mu)\leq \frac{C}{|\alpha|(\mu-1)}$.

\noindent {\bf Case 2: $1-\frac{\kappa}{|\alpha|^\frac12} < \mu \leq 1 + \frac{\kappa}{|\alpha|^\frac12}$.} Take $m_1=m_2^2= |\alpha|^\frac12$, for which $1+\frac{\kappa}{m_1}=1+\frac{\kappa}{m_2^2} = 1 + \frac{\kappa}{|\alpha|^\frac12}\geq \mu>1-\frac{\kappa}{m_1}=1-\frac{\kappa}{m_2^2}$ holds. Then $h_1(m_1,\mu)^2 = \frac{1}{m_1}$ and $h_2(m_2,\mu)^2=\frac{1}{m_2^4}$, which gives $F(\alpha,\mu)\leq \frac{C}{|\alpha|^\frac12}$.

\noindent {\bf Case 3: $\frac12 < \mu \leq 1-\frac{\kappa}{|\alpha|^\frac12}$.} Take $m_1=(\frac{|\alpha|}{1-\mu} )^\frac13$ and $m_2 = \kappa^{-\frac13} \big ( |\alpha| (1-\mu)^\frac12 \big )^\frac13$, for which $\frac{\kappa}{m_1} = \frac{\kappa (1-\mu)^\frac13}{|\alpha|^\frac13}= (\frac{\kappa}{|\alpha|^\frac12})^\frac23 \kappa^\frac13 (1-\mu)^\frac13\leq \kappa^\frac13 (1-\mu)<1-\mu$, and $\frac{\kappa}{m_2^2}=\frac{\kappa^\frac43}{|\alpha|^\frac23 (1-\mu)^\frac13} = (\frac{\kappa}{|\alpha|^\frac12})^\frac43 \frac{1}{(1-\mu)^\frac13}\leq 1-\mu$. Then $h_1(m_1,\mu)^2 = \frac{1}{m_1^2(1-\mu)}$ and $h_2 (m_2,\mu)^2=\frac{1-\mu}{m_2^2}$, which gives $F(\alpha,\mu) \leq \frac{C}{|\alpha|^\frac23 (1-\mu)^\frac13}$.

\noindent {\bf Case 4: $\frac{1}{|\alpha|^\frac13} < \mu\leq \frac12$.} Take $m_1=\frac{|\alpha|^\frac13}{\mu}$ and $m_2=|\alpha|^\frac13 \mu^\frac12$, for which $\frac{1}{m_1^\frac12} = \frac{\mu^\frac12}{|\alpha|^\frac16}<\mu$. Then $h_1 (m_1,\mu)^2 = \frac{1}{m_1^2\mu^3}$ and $h_2(m_2,\mu)^2=\frac{\mu^3}{m_2^3}$, which gives $F(\alpha,\mu)\leq \frac{C}{|\alpha|^\frac23 \mu}$.

\noindent {\bf Case 5: $\frac{1}{|\alpha|^\frac12} < \mu \leq \frac{1}{|\alpha|^\frac13}$.} Take $m_1=|\alpha| \mu$ and $m_2=|\alpha|^\frac13 \mu^\frac12$, for which $\frac{1}{10m_1}=\frac{1}{10|\alpha|\mu}<\frac{\mu}{10}<\mu$. Then $h_1(m_1,\mu)^2=\mu$ and $h_2(m_2,\mu)^2=\frac{\mu^3}{m_2^2}$, which gives $F(\alpha,\mu)\leq C\mu$.

\noindent {\bf Case 6: $-\frac{1}{|\alpha|^\frac12} < \mu\leq \frac{1}{|\alpha|^\frac12}$.} Take $m_1= \frac{|\alpha|^\frac12}{10}$ and $m_2=|\alpha|^\frac13 |\mu|^\frac12$, for which $|\mu| \leq \frac{1}{|\alpha|^\frac12} =\frac{1}{10m_1}$. Then $h_1(m_1,\mu)^2 = \frac{1}{m_1}$ and $h_2 (m_2,\mu) \leq \frac{|\mu|^3}{m_2^2}$, which gives $F(\alpha,\mu)\leq \frac{C}{|\alpha|^\frac12}$.

\noindent {\bf Case 7: $\mu\leq -\frac{1}{|\alpha|^\frac12}$.}  Take $m_1=\frac{1}{|\mu|}$ and $m_2=100$, for which $\mu=-|\mu|=-\frac{1}{m_1}<-\frac{1}{10m_1}$. Then $h_1(m_1,\mu)^2 =h_2(m_2,\mu)^2=0$, which gives $F(\alpha,\mu) \leq \frac{C}{|\alpha|\, |\mu|}$.

\

Summarizing these above, we obtain 
\begin{thm}\label{thm.Oseen} The exist positive numbers $C$ and $\alpha_0$ such that the following resolvent estimate holds for all $\lambda\in \R$ and for all $\alpha$ with $|\alpha|\geq \alpha_0$.
\begin{align}\label{est.thm.Oseen}
\| (i\lambda + L_\alpha)^{-1} \|_{Y\rightarrow Y} \leq C
\begin{cases}
& \displaystyle \frac{1}{|\alpha|(\frac{\lambda}{\alpha}-1)} \qquad \text{if} ~~ \frac{\lambda}{\alpha}>1+\frac{1}{|\alpha|^\frac12}\,,\\ 
& \displaystyle \frac{1}{|\alpha|^\frac12} \qquad \qquad  \text{if} ~~ 1-\frac{1}{|\alpha|^\frac12} <\frac{\lambda}{\alpha} \leq 1 + \frac{1}{|\alpha|^\frac12}\,,\\
& \displaystyle \frac{1}{|\alpha|^\frac23 (1-|\frac{\lambda}{\alpha}|)^\frac13} \qquad \text{if} ~~ \frac12 < \frac{\lambda}{|\alpha|} \leq 1- \frac{1}{|\alpha|^\frac12}\,,\\
& \displaystyle \frac{1}{|\alpha|^\frac23 \frac{\lambda}{\alpha}} \qquad \qquad \text{if} ~~ \frac{1}{|\alpha|^\frac13} < \frac{\lambda}{\alpha} \leq \frac12\,,\\
& \displaystyle \frac{\lambda}{\alpha} \qquad \qquad \quad  ~~ \text{if} ~~ \frac{1}{|\alpha|^\frac12} <  \frac{\lambda}{\alpha} \leq  \frac{1}{|\alpha|^\frac13}\,,\\
& \displaystyle \frac{1}{|\alpha|^\frac12} \qquad \qquad \text{if} ~~ -\frac{1}{|\alpha|^\frac12} < \frac{\lambda}{\alpha} \leq \frac{1}{|\alpha|^\frac12}\,,\\
& \displaystyle \frac{1}{|\lambda|}  \qquad \qquad  \quad \text{if} ~~ \frac{\lambda}{\alpha} \leq -\frac{1}{|\alpha|^\frac12}\,.
\end{cases}
\end{align}
In particular, we have $\displaystyle \sup_{\lambda\in \R} \|(i\lambda + L_\alpha)^{-1}\|_{Y\rightarrow Y} \leq \frac{C}{|\alpha|^\frac13}$. 
\end{thm}
The bound $\displaystyle \sup_{\lambda\in \R} \|(i\lambda + L_\alpha)^{-1}\|_{Y\rightarrow Y} \leq \frac{C}{|\alpha|^\frac13}$ is firstly shown in \cite{LWZ} by constructing the wave operator, and in \cite{LWZ} the optimality of the rate $O(|\alpha|^{-\frac13})$ is also proved. Theorem \ref{thm.Oseen} gives a different proof of their result without using the wave operator. Although \eqref{est.thm.Oseen} looks complicated, the dependence on $\alpha$ in each regime is compatible with the optimal result \cite{GaGaNi} for the case when $\hat{\Lambda}$ does not contain a nonlocal part. The rate $O(|\alpha|^{-\frac13})$ appears in the regime $\frac{\lambda}{\alpha} \sim O (|\alpha|^{-\frac13})$ and is related to the behavior of $\rho(r) \approx \frac{4}{r^2}$ for $r\gg 1$.

\section*{Acknowledgement}
The first author is partially supported by NSERC Discovery grant \# 371637-2014,
and also acknowledges the kind hospitality of the New York University in Abu Dhabi.
The second author is partially supported by JSPS Program for Advancing Strategic International Networks
to Accelerate the Circulation of Talented Researchers, 'Development of Concentrated Mathematical Center Linking to Wisdom of  the Next Generation', which is organized by Mathematical Institute of Tohoku University. 
The third author is partially supported by the NSF grant DMS-1716466.


\begin{thebibliography}{}

\bibitem{AB} A. L. Afendikov and K. I. Babenko;
\newblock Bifurcation in the presence of a symmetry group and loss of stability of some plane flows of a viscous fluid. Soviet Math. Dokl., {\bf 33} (1986), 742-747.


\bibitem{BW} M. Beck and C. E. Wayne;
\newblock Metastability and rapid convergence to quasi-stationary bar states for the two-dimensional Navier-Stokes equations. Proc. Roy. Soc. Edinburgh Sect. A  {\bf 143}  (2013),  no. 5, 905-927. 


\bibitem{BeGerMa1} J. Bedrossian, P. Germain, and N. Masmoudi;
Dynamics near the subcritical transition of the 3D Couette flow I: Below threshold case, 2015. arXiv 1506.03720. (to appear in Memoire of the AMS).


\bibitem{BeGerMa2} J. Bedrossian, P. Germain, and N. Masmoudi;
Dynamics near the subcritical transition of the 3D Couette flow II: Above threshold case, 2015. arXiv 1506.03721.


\bibitem{BeGerMa3} J. Bedrossian, P. Germain, and N. Masmoudi;
On the stability threshold for the 3D Couette flow in Sobolev regularity. Ann. of Math. (2)  {\bf 185} (2017) no. 2, 541-608.


\bibitem{BeMaVi} J. Bedrossian, N. Masmoudi, and V. Vicol;
Enhanced dissipation and inviscid damping in the inviscid limit of the Navier-Stokes equations near the two dimensional Couette flow. Arch. Ration. Mech. Anal.  {\bf 219} (2016), no. 3, 1087-1159.


\bibitem{BS} F. Bouchet and E. Simonnet;
\newblock Random Changes of Flow Topology in Two-Dimensional and Geophysical Turbulence. Phys. Rev. Lett. {\bf 102} (2009), 09450.


\bibitem{CKRZ} P. Constantin, A. Kiselev, L. Ryzhik, and A. Zlato${\rm \check{s}}$;
\newblock Diffusion and mixing in fluid flow. 
Ann. of Math. (2), {\bf 168} (2008), 643-674.


\bibitem{De1} W. Deng;
\newblock Resolvent estimates for a two-dimensional non-self-adjoint operator.
Commun. Pure Appl. Anal. {\bf 12} (2013), 547-596.

\bibitem{De2} W. Deng;
\newblock Pseudospectrum for Oseen vortices operators.
Int. Math. Res. Not. IMRN {\bf 2013} (2013), 1935-1999.


\bibitem{GaGaNi} I. Gallagher, Th. Gallay, and F. Nier;
\newblock Special asymptotics for large skew-symmetric perturbations of the harmonic oscillator. Int. Math. Res. Not. {\bf 12} (2009), 2147-2199.

\bibitem{Ga2017} Th. Gallay;
\newblock Enhanced dissipation and axisymmetrization of two-dimensional viscous vortices.
Arch. Ration. Mech. Anal. {\bf 230} (2018), 939-975. 


\bibitem{GaMa} Th. Gallay and Y. Maekawa;
\newblock Existence and Stability of Viscous Vortices,  in Handbook of Mathematical Analysis in Mechanics of Viscous Fluids, Springer International Publishing Switzerland 2016. 
DOI 10.1007/978-3-319-10151-4$\_$13-1


\bibitem{GW} Th. Gallay and G. E. Wayne;
\newblock Global stability of vortex solutions of the two-dimensional Navier-Stokes equation. Commun. Math. Phys. {\bf 255} (2005), 97-129.


\bibitem{Iu} V. I. Iudovich;
\newblock Example of the generation of a secondary statonary or periodic flow when there is loss of stability of the laminar flow of a viscous incompressible fluid. J. Appl. Math. Mech., {\bf 29} (1965), 527-544.


\bibitem{Ka} T. Kato;
\newblock {\it Perturbation Theory for Linear Operators, 2nd edn.} Springer, London (1976).


\bibitem{LWZ} T. Li, D. Wei, and Z. Zhang;
\newblock Pseudospectral and spectral bounds for the Oseen vortices operator. Preprint. https://arxiv.org/abs/1701.06269


\bibitem{LWZ2} T. Li, D. Wei, and Z. Zhang;
Pseudospectral bound and transition threshold for the $3$D Kolmogorov flow.
arXiv:1801.05645.


\bibitem{LX} Z. Lin, and M. Xu;
\newblock Metastability of Kolmogorov flows and inviscid damping of shear flows. Preprint. https://arxiv.org/abs/1707.00278.


\bibitem{Ma} Y. Maekawa;
\newblock Spectral properties of the linearization at the Burgers vortex in the high rotation limit.
Journal of Mathematical Fluid Mechanics, {\bf 13} (2011) 515-532. 


\bibitem{Mar} C. Marchioro;
\newblock An example of absence of turbulence for any Reynolds number.
Commun. Math. Phys., {\bf 105} (1986), 99-106.


\bibitem{MaMi} M. Matsuda and S. Miyatake;
\newblock Bifurcation analysis of Kolmogorov flows.
T${\rm \hat{o}}$hoku Math. J., {\bf 36} (1984), 623-646.


\bibitem{MSMOM} W. H. Matthaeus, W. T. Stribling, D. Martinez, S. Oughton, and D. Montgomery;
\newblock Decaying, two dimensional, Navier-Stokes turbulence at very long times. Physica D {\bf 51}, 531.


\bibitem{MeSi} L. D. Meshalkin and Y. G. Sinai;
\newblock Investigation of the stability of a stationary solution to a system of equations for the plane movement of a incompressible viscous liquid.
J. Appl. Math. Mech., {\bf 25} (1962), 1700-1705.


\bibitem{OS} H. Okamoto and M. Sh${\rm \bar{o}}$ji;
\newblock Bifurcation diagram in Kolmogorov's problem of viscous incompressible fluid on $2$-D Tori.
Japan J. Indus. Appl. Math., {\bf 10} (1993), 191-218.


\bibitem{Vi} C. Villani;
\newblock {\it Hypocoercivity.} Memoirs of the Americal Mathematical Society (Providence, RI: American Mathematical Society, 2009).
 

\bibitem{W} D. Wei;
Diffusion and mixing in fluid flow via the resolvent estimate. arXiv:1811.11904.


\bibitem{WZZ} D. Wei, Z. Zhang, W. Zhao;
\newblock Linear inviscid damping and enhanced dissipation for the Kolmogorov flow. Preprint, 
arXiv:1711.01822.


\bibitem{WZZ2} D. Wei, Z. Zhang, W. Zhao;
Linear inviscid damping and vorticity depletion for shear flows.
Ann. PDE {\bf 5} (2019), no. 1, Art. 3, 101 pp.


\bibitem{Ya} M. Yamada;
\newblock Nonlinear stability theory of spatially periodic parallel flows.
J. Phys. Soc. Japan, {\bf 55} (1986), 3073-3079.


\bibitem{YMC} Z. Yin, D. C. Montgomery, and H. J. H. Clercx;
\newblock Alternative statistical-mechanical descriptions of decaying two-dimensional turbulence in terms of ``patches'' and ``points''. Phys. Fluids {\bf 15} (2009), 1937-1953.


\bibitem{Z} A. Zlato${\rm \check{s}}$;
Diffusion in fluid flow: dissipation enhancement by flows in $2$D. 
Comm. Partial Differential Equations {\bf 35} (2010), 496-534. 




\end{thebibliography}
\end{document}